\documentclass[preprint,3p,times]{elsarticle}

\usepackage[section]{placeins}
\usepackage{multirow}
\usepackage{makecell}
\usepackage{longtable}
\usepackage{booktabs}
\usepackage{supertabular}
\usepackage[table]{xcolor}

\usepackage[bf]{caption2}

\usepackage{amssymb}
\usepackage{amsmath}
\newtheorem{theorem}{Theorem}
\newtheorem{lemma}{Lemma}

\newtheorem{definition}{Definition}

\newtheorem{example}{Example}

\journal{Elsevier}

\begin{document}

\begin{frontmatter}

\title{A new mapped WENO scheme using order-preserving mapping}

\author[a]{Ruo Li}
\ead{rli@math.pku.edu.cn}

\author[b,c]{Wei Zhong\corref{cor1}}
\ead{zhongwei2016@pku.edu.cn}

\cortext[cor1]{Corresponding author}

\address[a]{CAPT, LMAM and School of Mathematical Sciences, Peking
University, Beijing 100871, China}

\address[b]{School of Mathematical Sciences, Peking University,
Beijing 100871, China}

\address[c]{Northwest Institute of Nuclear Technology, Xi'an 
710024, China}

\begin{abstract}

  Existing mapped WENO schemes can hardly prevent spurious
  oscillations while preserving high resolutions at long output
  times. We reveal in this paper the essential reason for such
  phenomena. It is actually caused by that the mapping function in
  these schemes can not preserve the order of the nonlinear weights 
  of the stencils. The nonlinear weights may be increased for 
  non-smooth stencils and be decreased for smooth stencils. It is 
  then indicated to require the set of mapping functions to be
  \textit{Order-Preserving} in mapped WENO schemes. Therefore, we
  propose a new mapped WENO scheme with a set of mapping functions to
  be order-preserving which exhibits a remarkable advantage over the
  mapped WENO schemes in references. For long output time 
  simulations, the new scheme has the capacity to attain high 
  resolutions and avoid spurious oscillations near discontinuities 
  meanwhile.

\end{abstract}



\end{frontmatter}

\section{Introduction}
\label{secIntroduction}
Many numerical methods have been studied to solve the hyperbolic 
problems which may generate discontinuities as time evolves in its 
solution even if the initial condition is smooth, especially for 
nonlinear cases. As the discontinuities often cause spurious 
oscillations in numerical calculations, it is very difficult to 
design high order numerical schemes. Thus, the essentially 
non-oscillatory (ENO) schemes \cite{ENO1987JCP71, ENO1987V24, 
ENO1986, ENO1987JCP83} and weighted ENO (WENO) schemes 
\cite{ENO-Shu1988, ENO-Shu1989, WENO-LiuXD,WENO-JS, WENOoverview} 
have developed quite successfully in recent decades, and they are 
very popular methods to solve the hyperbolic conservation laws 
because of their success to the ENO property. The goal of this paper 
is to devise a new version of the fifth-order finite volume WENO 
scheme for solving the following hyperbolic conservation laws
\begin{equation}\left\{
\begin{array}{l}
\dfrac{\partial \mathbf{u}}{\partial t} +
\displaystyle\sum\limits_{\alpha = 1}^{d} \dfrac{\partial
\mathbf{f}_{\alpha}(\mathbf{u})}{\partial x_{\alpha}} = 0,
\quad x_{\alpha} \in \mathbb{R}, t > 0, \\
\mathbf{u}(\mathbf{x},0) = \mathbf{u}_{0}(\mathbf{x}),
\end{array}\right.
\label{governingEquation}
\end{equation}
with proper boundary conditions. Here, $\mathbf{u} = (u_{1}, \cdots, 
u_{m}) \in \mathbb{R}^{m}$ are the conserved variables and 
$\mathbf{f}_{\alpha}: \mathbb{R}^{m} \rightarrow \mathbb{R}^{m}$, 
$\alpha = 1,2,\cdots,d$ are the Cartesian components of flux.  

The first WENO scheme was developed by Liu et al. \cite{WENO-LiuXD} 
in the finite volume version. It converts an $r$th-order ENO scheme 
\cite{ENO1987JCP71,ENO1987V24,ENO1986} into an $(r+1)$th-order WENO 
scheme through a convex combination of all candidate substencils 
instead of only using the optimal smooth candidate stencil by the 
original ENO scheme. By introducing a different definition of the 
smoothness indicators used to measure the smoothness of the 
numerical solutions on substencils, the classic $(2r-1)$th-order 
WENO-JS scheme was proposed by Jiang and Shu \cite{WENO-JS}. The 
WENO-JS scheme has been successfully used in a wide number of 
applications. The WENO methodology within the general framework of 
smoothness indicators and nonlinear weights proposed in the WENO-JS 
scheme \cite{WENO-JS} is still in development to improve the 
convergence rate in smooth regions and reduce the dissipation near 
the discontinuities \cite{WENO-M,WENO-IM, WENO-PM, WENO-Z, WENO-Z01, 
WENO-Z02, WENO-Accuracy, WENO-GlobalAccuracy}.

Henrick et al. \cite{WENO-M} pointed out that the classic WENO-JS 
scheme was less than fifth-order for many cases such as at or near 
critical points of order $n_{\mathrm{cp}} = 1$ in the smooth 
regions. Here, we refer to $n_{\mathrm{cp}}$ as the order of the 
critical point; e.g., $n_{\mathrm{cp}} = 1$ corresponds to 
$f'=0, f'' \neq 0$ and $n_{\mathrm{cp}} = 2$ corresponds to 
$f'=0, f'' = 0, f''' \neq 0$, etc. The necessary and sufficient 
conditions on the nonlinear weights for optimality of the 
convergence rate of the fifth-order WENO schemes were derived by 
Henrick et al. in \cite{WENO-M}. These conditions were reduced to a 
simple sufficient condition \cite{WENO-Z} which could be easily 
extended to the $(2r-1)$th-order WENO schemes \cite{WENO-IM}. Then, 
by designing a mapping function that satisfies the sufficient 
condition to achieve the optimal order of accuracy, the original 
mapped WENO scheme, named WENO-M, was devised by Henrick et al. 
\cite{WENO-M}.

Recently, Feng et al. \cite{WENO-PM} noted that, when the WENO-M
scheme was used for solving the problems with discontinuities for 
long output times, its mapping function may amplify the effect from 
the non-smooth stencils leading to a potential loss of accuracy near
discontinuities. To address this issue, they proposed a piecewise
polynomial mapping function with two additional requirements, that 
is, $g'(0) = 0$ and $g'(1)=0$ ($g(x)$ denotes the mapping function), 
to the original criteria in \cite{WENO-M}. However, the resultant
WENO-PM$k$ scheme \cite{WENO-PM} may generate the non-physical
oscillations near the discontinuities as shown in Fig. 8 of
\cite{WENO-IM} and Figs. 3-8 of \cite{WENO-RM260}. Later, Feng et
al. \cite{WENO-IM} devised an improved mapping method which is
referred to as WENO-IM($k, A$) where $k$ is a positive even integer
and $A$ a positive real number.  The broader class of the improved
mapping functions of the family of the WENO-IM($k, A$) schemes
includes the mapping function of the WENO-M scheme as a special case
by taking $k=2$ and $A=1$. Feng et al. indicated that by taking $k=2$
and $A=0.1$ in the fifth-order WENO-IM($k, A$) scheme, far better
numerical solutions with less dissipation and higher resolution could
be obtained than that of the fifth-order WENO-M scheme. However, the
possible over-amplification of the contributions from non-smooth
stencils exists for the WENO-IM($k, A$) scheme as the first 
derivative of its mapping function satisfies
$\big( g^{\mathrm{IM}} \big)'_{s}(0; k, A)= 1 +
\frac{1}{Ad^{k-1}_{s}}$, and this excessive amplification of the
weights of non-smooth stencils causes spurious oscillations to occur
in the solution which may even render the algorithm unstable
\cite{WENO-AIM}, especially for higher order cases. The numerical
experiments in \cite{WENO-RM260} showed that the seventh- and ninth-
order WENO-IM(2, 0.1) schemes generated evident spurious oscillations
near discontinuities when the output time is large. In addition, our
calculations as shown in Figs.
\ref{fig:BiCWP:800}-\ref{fig:BiCWP:6400} of this paper indicate that,
even for the fifth-order WENO-IM(2, 0.1) scheme, the spurious
oscillations are also produced when the grid number increases or a
different initial condition is used.

Many other improved mapped WENO schemes, such as the WENO-PPM$n$
\cite{WENO-PPM5}, WENO-RM($mn0$) \cite{WENO-RM260}, WENO-MAIM$i$
\cite{WENO-MAIMi}, WENO-ACM \cite{WENO-ACM} schemes and et al., have
been successfully developed to improve the performances of the 
classic WENO-JS scheme in some ways, like achieving optimal 
convergence orders near critical points in smooth regions, having 
lower numerical dissipations and achieving higher resolutions near 
discontinuities, and we refer to the references for more details. 
However, as mentioned in previously published works of literature 
\cite{WENO-IM,WENO-RM260}, most of the existing improved mapped WENO 
schemes could not prevent the generation of the spurious 
oscillations near discontinuities, especially for long output time 
simulations.

Taking a long output time simulation of the linear advection problem
with discontinuities as an example, we make a further study of the
nonlinear weights of the existing mapped WENO schemes developed in
\cite{WENO-PM,WENO-IM,WENO-MAIMi} and the MIP-WENO-ACM$k$ scheme,
which is a generalized form of the WENO-ACM scheme \cite{WENO-ACM}
(see details in subsection \ref{subsec:MIP-WENO-ACMk}). We find that,
in many points, the order of the nonlinear weights for the 
substencils of the same global stencil has been changed in the 
mapping process of all these considered mapped WENO schemes. The 
order change of the nonlinear weights is caused by weights 
increasing of non-smooth substencils and weights decreasing of 
smooth substencils. Through theoretical analysis and extensive 
numerical tests, we reveal that this is the essential cause of the 
potential loss of accuracy of the WENO-M scheme and the spurious 
oscillation generation of the existing improved mapped WENO schemes.

Indicated by the observation above, we give the definition of the
\textit{order-preserving} (see Definition \ref{def:OPM} below) 
mapping and suggest it as an additional criterion in the design
of the mapping function. Then we propose a new mapped WENO scheme,
referred to as MOP-WENO-ACM$k$ below, with its mapping function 
satisfying the additional criterion. This new version of the mapped 
WENO scheme achieves the optimal convergence order of accuracy even 
at critical points. It also has low numerical dissipation and high 
resolution but does not generate spurious oscillation near 
discontinuities even if the output time is large.

At first, a series of accuracy tests show the capacity of the 
proposed scheme to achieve the optimal convergence order in smooth 
regions with first-order critical points and its advantage in long 
output time simulations of the problems with very high-order 
critical points. Some linear advection examples with long output 
times are then presented to demonstrate that the proposed scheme can 
obtain high resolution and does not generate spurious oscillation 
near discontinuities. At last, some benchmark problems modeled via 
the two-dimensional Euler equations are computed by various 
considered WENO schemes to compare with the proposed scheme. It is 
clear that the proposed scheme exhibits significant advantages in 
preventing spurious oscillations.

The organization of this paper is as follows. Preliminaries to
understand the finite volume method and the procedures of the WENO-JS
\cite{WENO-JS}, WENO-M \cite{WENO-M}, WENO-PM6 \cite{WENO-PM} and 
WENO-IM($k, A$) \cite{WENO-IM} schemes are reviewed in Section 2. In 
Section 3, we present a detailed analysis on how the nonlinear 
weights are mapped in some existing mapped WENO schemes 
and the consequences of these mappings on the numerical solutions. 
In Section 4, we propose a set of mapping functions that is 
order-preserving, as well as its properties, and apply it to 
construct the new mapped WENO scheme. Numerical experiments are 
presented in Section 5 to illustrate the advantages of the proposed 
WENO scheme. Finally, some concluding remarks are made in Section 6 
to close this paper.


\section{Brief review of the fifth-order WENO schemes}
\label{secMappedWENO}

\subsection{Finite volume methodology}
We consider the finite volume method for the following 
one-dimensional scalar hyperbolic conservation laws
\begin{equation}\left\{
\begin{array}{l}
\dfrac{\partial u}{\partial t}+\dfrac{\partial f(u)}{\partial x}=0,
\quad x_{l} < x < x_{r}, t > 0, \\
u(x,0) = u_{0}(x).
\end{array}\right.
\label{eq:1D-hyperbolicLaw}
\end{equation}
For brevity in the description, we assume that the computational 
domain is discretized into uniform cells $I_{j} = [x_{j-1/2}, 
x_{j+1/2}], j = 1,\cdots,N$ with the uniform cell size $\Delta x = 
\frac{x_{r} - x_{l}}{N}$, and the associated cell centers and cell 
boundaries are denoted by $x_{j} = x_{l} + (j - 1/2)\Delta x$ and 
$x_{j \pm 1/2} = x_{j} \pm \Delta x/2$ respectively. Let 
$\bar{u}(x_{j}, t)=\dfrac{1}{\Delta x}\int_{x_{j-1/2}}^{x_{j+1/2}}u
(\xi,t)\mathrm{d}\xi$ be the cell average of $I_{j}$, then by 
integrating Eq.(\ref{eq:1D-hyperbolicLaw}) over the control volumes 
$I_{j}$ and employing some simple mathematical manipulations, we 
approximate Eq.(\ref{eq:1D-hyperbolicLaw}) by the following finite 
volume conservative formulation
\begin{equation}
\dfrac{\mathrm{d}\bar{u}_{j}(t)}{\mathrm{d}t} \approx -\dfrac{1}
{\Delta x}\bigg( \hat{f}(u_{j+1/2}^{-},u_{j+1/2}^{+}) - 
\hat{f}(u_{j-1/2}^{-},u_{j-1/2}^{+}) \bigg),
\label{eq:discretizedFunction}
\end{equation}
where $\bar{u}_{j}(t)$ is the numerical approximation to the cell 
average $\bar{u}(x_{j}, t)$, and the numerical flux $\hat{f}(u_{j \pm
1/2}^{-},u_{j \pm 1/2}^{+})$ where $u_{j \pm 1/2}^{\pm}$ refer to 
the limits of $u$ is a replacement of the physical flux function 
$f(u)$ at the cell boundaries $x_{j\pm1/2}$. The values of 
$u_{j\pm1/2}^{\pm}$ can be obtained by the technique of 
reconstruction like some WENO reconstructions narrated later. In 
this paper, we choose the global Lax-Friedrichs flux $\hat{f}(a,b) = 
\frac{1}{2}[f(a) + f(b) -\alpha(b - a)]$, where $\alpha = \max_{u} 
\lvert f'(u) \rvert$ is a constant and the maximum is taken over the 
whole range of $u$. For the hyperbolic conservation laws system, a 
local characteristic decomposition is commonly used and more details 
can be found in \cite{WENO-JS}. In two-dimensional Cartesian meshes, 
two classes of finite volume WENO schemes are studied in detail in 
\cite{FVMaccuracyProofs03}, and we take the one denoted as class A 
in this paper.

The ordinary differential equation (ODE) system 
Eq.(\ref{eq:discretizedFunction}) can be solved using a suitable 
time discretization, and the following explicit, third-order, strong 
stability preserving (SSP) Runge-Kutta method \cite{ENO-Shu1988,
SSPRK1998,SSPRK2001} is employed in our calculations
\begin{equation*}
\begin{array}{l}
\begin{aligned}
&u^{(1)} = u^{n} + \Delta t \mathcal{L}(u^{n}), \\
&u^{(2)} = \dfrac{3}{4} u^{n} + \dfrac{1}{4} u^{(1)} + \dfrac{1}{4} 
\Delta t \mathcal{L}(u^{(1)}), \\
&u^{n + 1} = \dfrac{1}{3} u^{n} + \dfrac{2}{3} \Delta t \mathcal{L}(u
^{(2)}),
\end{aligned}
\end{array}
\end{equation*}
where 
\begin{equation*}
\mathcal{L}(u_{j}):=-\dfrac{1}{\Delta x}\bigg( \hat{f}(u_{j+1/2}^{-},
u_{j+1/2}^{+}) - \hat{f}(u_{j-1/2}^{-},u_{j-1/2}^{+}) \bigg),
\end{equation*}
and $u^{(1)}$, $u^{(2)}$ are the intermediate stages, $u^{n}$ 
is the value of $u$ at time level $t^{n} = n\Delta t$, and $\Delta t$
is the time step satisfying some proper CFL condition. As mentioned 
earlier, the WENO reconstructions will be applied to obtain 
$\mathcal{L}(u_{j})$.

\subsection{WENO-JS}
We review the process of the fifth-order WENO-JS reconstruction 
\cite{WENO-JS}. As the right-biased reconstruction $u_{j + 1/2}^{+}$ 
can easily be obtained by mirror symmetry with respect to the 
location $x_{j + 1/2}$ of that for the left-biased one 
$u_{j + 1/2}^{-}$, we describe only the reconstruction procedure of 
$u_{j + 1/2}^{-}$. For simplicity of notation, we do not use the 
subscript ``-'' in the following content. 

For constructing the values of $u_{j + 1/2}$ from known cell average 
values, the fifth-order WENO-JS scheme uses a 5-point global stencil 
$S^{5}=\{I_{j-2}, I_{j-1}, I_{j}, I_{j+1}, I_{j+2}\}$. Normally,
$S^{5}$ is subdivided into three 3-point substencils $S_{s} = 
\{I_{j+s-2}, I_{j+s-1}, I_{j+s}\}, s = 0, 1, 2$. Explicitly, the 
third-order approximations of $u(x_{j+1/2},t)$ associated with these 
substencils are given by
\begin{equation}
\begin{array}{l}
\begin{aligned}
&u_{j+1/2}^{0} = \dfrac{1}{6}(2\bar{u}_{j-2} - 7\bar{u}_{j-1}
+ 11\bar{u}_{j}), \\
&u_{j+1/2}^{1} = \dfrac{1}{6}(-\bar{u}_{j-1} + 5\bar{u}_{j}
+ 2\bar{u}_{j+1}), \\
&u_{j+1/2}^{2} = \dfrac{1}{6}(2\bar{u}_{j} + 5\bar{u}_{j+1}
- 2\bar{u}_{j+2}).
\end{aligned}
\end{array}
\label{eq:approx_ENO}
\end{equation}

The above three third-order approximations are combined in a 
weighted average to define the fifth-order WENO approximation of 
$u(x_{j+1/2},t)$, 
\begin{equation}
u_{j + 1/2} = \sum\limits_{s = 0}^{2}\omega_{s}u_{j + 1/2}^{s},
\label{eq:approx_WENO}
\end{equation}
where $\omega_{s}$ is the nonlinear weight of the substencil $S_{s}$
. In the classic WENO-JS scheme, $\omega_{s}$ is calculated as
\begin{equation} 
\omega_{s}^{\mathrm{JS}} = \dfrac{\alpha_{s}^{\mathrm{JS}}}{\sum_{l =
 0}^{2} \alpha_{l}^{\mathrm{JS}}}, \alpha_{s}^{\mathrm{JS}} = \dfrac{
 d_{s}}{(\epsilon + \beta_{s})^{2}}, \quad s = 0,1,2,
\label{eq:weights:WENO-JS}
\end{equation} 
where $d_{0} = 0.1, d_{1} = 0.6, d_{2} = 0.3$ are ideal weights of 
$\omega_{s}$ satisfying $\sum\limits_{s=0}^{2} d_{s} u^{s}_{j+1/2} =
u(x_{j + 1/2}, t) + O(\Delta x^{5})$ in smooth regions, $\epsilon$ 
is a small positive number introduced to prevent the denominator from
becoming zero, and the parameters $\beta_{s}$ are the smoothness 
indicators for the third-order approximations $u_{j+1/2}^{s}$ and 
their explicit forms defined by Jiang and Shu \cite{WENO-JS} is 
given as
\begin{equation*}
\begin{array}{l}
\begin{aligned}
\beta_{0} &= \dfrac{13}{12}\big(\bar{u}_{j - 2} - 2\bar{u}_{j - 1} + 
\bar{u}_{j} \big)^{2} + \dfrac{1}{4}\big( \bar{u}_{j - 2} - 4\bar{u}_
{j - 1} + 3\bar{u}_{j} \big)^{2}, \\
\beta_{1} &= \dfrac{13}{12}\big(\bar{u}_{j - 1} - 2\bar{u}_{j} + \bar
{u}_{j + 1} \big)^{2} + \dfrac{1}{4}\big( \bar{u}_{j - 1} - \bar{u}_{
j + 1} \big)^{2}, \\
\beta_{2} &= \dfrac{13}{12}\big(\bar{u}_{j} - 2\bar{u}_{j + 1} + \bar
{u}_{j + 2} \big)^{2} + \dfrac{1}{4}\big( 3\bar{u}_{j} - 4\bar{u}_{j 
+ 1} + \bar{u}_{j + 2} \big)^{2}.
\end{aligned}
\end{array}
\end{equation*}

The fifth-order WENO-JS scheme is able to achieve optimal order of 
accuracy in smooth regions without critical points. However, it 
loses accuracy and its order of accuracy decreases to third-order or 
even less at critical points. More details can be found in 
\cite{WENO-M}.

\subsection{WENO-M}
\label{subsecWENO-M}
It has been indicated that \cite{WENO-M,WENO-Z,WENO-IM,MWENO-P,
WENO-eta} a sufficient condition that ensures the fifth-order WENO 
schemes retaining optimal order of convergence is simply given by
\begin{equation}
\omega^{\pm}_{s} - d_{s} = O(\Delta x^{3}), \quad s = 0, 1, 2.
\label{eq:sufficient}
\end{equation}
The condition Eq.(\ref{eq:sufficient}) may not hold in the case of 
smooth extrema or at critical points when the fifth-order WENO-JS 
scheme is used. Henrick et al. \cite{WENO-M} proposed a fix to this 
deficiency in their WENO-M scheme by introducing a mapping function 
that makes $\omega_{s}$ approximating the ideal weights $d_{s}$ with 
increased accuracy. The mapping function of the nonlinear weights 
$\omega \in [0, 1]$ is given by
\begin{equation}
\big( g^{\mathrm{M}} \big)_{s}(\omega) = \dfrac{ \omega \big( d_{s} +
d_{s}^2 - 3d_{s}\omega + \omega^{2} \big) }{ d_{s}^{2} + (1 - 2d_
{s})\omega }, \quad \quad s = 0, 1, 2.
\label{mappingFunctionWENO-M}
\end{equation}
One can verify that $\big( g^{\mathrm{M}}\big)_{s}(\omega)$ meets the requirement in Eq.(\ref{eq:sufficient}), and clearly, 
this mapping function is a non-decreasing monotone function on 
$[0, 1]$ with finite slopes which satisfies the following properties.
\begin{lemma} 
The mapping function $\big( g^{\mathrm{M}} \big)_{s}(\omega)$
defined by Eq.(\ref{mappingFunctionWENO-M}) satisfies: \\

C1. $0 \leq \big(g^{\mathrm{M}}\big)_{s}(\omega)\leq 1, \big(
g^{\mathrm{M}}\big)_{s}(0)=0, \big(g^{\mathrm{M}}\big)_{s}(1) = 1$;

C2. $\big( g^{\mathrm{M}} \big)_{s}(d_{s}) = d_{s}$;

C3. $\big( g^{\mathrm{M}} \big)_{s}'(d_{s}) = \big(g^{\mathrm{M}}
\big)_{s}''(d_{s}) = 0$. 
\label{lemmaWENO-Mproperties}
\end{lemma}

With the mapping function defined by 
Eq.(\ref{mappingFunctionWENO-M}), the nonlinear weights of the 
WENO-M scheme are defined as
\begin{equation*}
\omega_{s}^{\mathrm{M}} = \dfrac{\alpha _{s}^{\mathrm{M}}}{\sum_{l = 
0}^{2} \alpha _{l}^{\mathrm{M}}}, \alpha_{s}^{\mathrm{M}} = \big( g^{
\mathrm{M}} \big)_{s}(\omega^{\mathrm{JS}}_{s}), \quad s = 0,1,2,
\end{equation*}
where $\omega_{s}^{\mathrm{JS}}$ are calculated by 
Eq.(\ref{eq:weights:WENO-JS}).

In \cite{WENO-M}, it has been analyzed and proved in detail that the
WENO-M scheme can retain the optimal order of accuracy in
smooth regions even near the first-order critical points.

\subsection{WENO-PM$k$}
\label{subsecWENO-PM}
Recently, Feng et al. \cite{WENO-PM} noticed that the mapping 
function $\big( g^{\mathrm{M}} \big)_{s}(\omega)$ in 
Eq.(\ref{mappingFunctionWENO-M}) amplifies the effect from 
the non-smooth stencils by a factor of $(1+1/d_{s})$ as its first 
derivative satisfies $\big( g^{\mathrm{M}} \big)'_{s}(0) = 1+1/d_{s}$
. They argued that this may cause the potential loss of accuracy 
near the discontinuities or the parts with sharp gradients. To 
address this issue, Feng et al. \cite{WENO-PM} add two requirements, 
that is, $g'_{s}(0) = 0$ and $g'_{s}(1) = 0$, to the original 
criteria as shown in Lemma \ref{lemmaWENO-Mproperties}. To meet 
these criteria they proposed a new mapping by the following 
piecewise polynomial function
\begin{equation}
\big( g^{\mathrm{PM}} \big)_{s}(\omega) = c_{1}(\omega - d_{s})^{k+1}
(\omega + c_{2}) + d_{s}, \quad \quad k \geq 2, s = 0,1,2,
\label{mappingFunctionWENO-PM}
\end{equation}
where $c_{1}, c_{2}$ are constants with specified parameters $k$ and 
$d_{s}$, taking the following forms
\begin{equation*}
\begin{array}{ll}
c_{1} = \left\{
\begin{array}{ll}
\begin{aligned}
&(-1)^{k}\dfrac{k+1}{d_{s}^{k+1}}, & 0 \leq \omega \leq d_{s}, \\
&-\dfrac{k+1}{(1-d_{s})^{k+1}},    & d_{s} < \omega \leq 1,
\end{aligned}
\end{array}\right.
&
c_{2} = \left\{
\begin{array}{ll}
\begin{aligned}
&\dfrac{d_{s}}{k+1},        &  0 \leq \omega \leq d_{s}, \\
&\dfrac{d_{s}-(k+2)}{k+1},  & d_{s} < \omega \leq 1. 
\end{aligned}
\end{array}\right.
\end{array}
\end{equation*}

\begin{lemma} 
The mapping function $\big( g^{\mathrm{PM}} \big)_{s}(\omega)$
defined by Eq.(\ref{mappingFunctionWENO-PM}) satisfies: \\

C1. $\big( g^{\mathrm{PM}} \big)'_{s}(\omega) \geq 0, \omega \in [0, 
1]$;

C2. $\big( g^{\mathrm{PM}} \big)_{s}(0) = 0, \big( g^{\mathrm{PM}} 
\big)_{s}(1) = 1, \big( g^{\mathrm{PM}} \big)_{s}(d_{s}) = d_{s}$;

C3. $\big( g^{\mathrm{PM}} \big)'_{s}(d_{s}) = \cdots = \big( g^{
\mathrm{PM}} \big)^{(k)}_{s}(d_{s}) = 0$;

C4. $\big( g^{\mathrm{PM}} \big)'_{s}(0) = \big( g^{\mathrm{PM}} 
\big)'_{s}(1) = 0$.
\label{lemmaWENO-PMproperties}
\end{lemma}

Similarly, with the mapping function defined by 
Eq.(\ref{mappingFunctionWENO-PM}) where the parameter $k$ is taken 
to be $6$ as recommended in \cite{WENO-PM}, the WENO-PM6 scheme is 
proposed by computing the nonlinear weights as
\begin{equation*}
\omega_{s}^{\mathrm{PM}6}=\dfrac{\alpha _{s}^{\mathrm{PM}6}}{\sum_{l 
= 0}^{2} \alpha _{l}^{\mathrm{PM}6}}, \alpha_{s}^{\mathrm{PM}6} = 
\big(g^{\mathrm{PM}6} \big)_{s}(\omega^{\mathrm{JS}}_{s}), 
\quad s = 0,1,2.
\end{equation*}

It has been shown by numerical experiments \cite{WENO-PM,WENO-RM260} 
that the two additional requirements are effective and the 
resolution near discontinuities of the WENO-PM6 scheme is 
significantly higher than those of the WENO-JS and WENO-M schemes, 
especially for long output times. We refer to \cite{WENO-PM} for 
more details.

\subsection{WENO-IM($k, A$)}
\label{subsecWENO-IM}
Feng et al. \cite{WENO-IM} has proposed the WENO-IM($k, A$) scheme 
by rewriting the mapping function of the WENO-M scheme as shown in 
Eq.(\ref{mappingFunctionWENO-M}). The broader class of improved 
mapping functions $\big( g^{\mathrm{IM}} \big)_{s}(\omega; k, A)$ is 
defined by
\begin{equation}
\big( g^{\mathrm{IM}} \big)_{s}(\omega; k, A) = d_{s} + \dfrac{\big( 
\omega - d_{s} \big)^{k + 1}A}{\big( \omega - d_{s} \big)^{k}A + 
\omega(1 - \omega)}, \quad A > 0, k = 2n, n \in \mathbb{N}^{+}, 
s = 0,1,2.
\label{mappingFunctionWENO-IM}
\end{equation}
Then, the associated nonlinear weights are given by
\begin{equation*}
\omega_{s}^{\mathrm{IM}} = \dfrac{\alpha _{s}^{\mathrm{IM}}}{\sum_{l 
= 0}^{2} \alpha _{l}^{\mathrm{IM}}}, \alpha_{s}^{\mathrm{IM}} = \big(
g^{\mathrm{IM}} \big)_{s}(\omega^{\mathrm{JS}}_{s}; k, A), 
\quad s = 0,1,2.
\end{equation*}
It is trivial to show that $\big( g^{\mathrm{M}} \big)_{s}(\omega)$ 
belongs to the $\big( g^{\mathrm{IM}} \big)_{s}(\omega; k, A)$ 
family of functions as $\big( g^{\mathrm{M}} \big)_{s}(\omega) = 
\big( g^{\mathrm{IM}} \big)_{s}(\omega; 2, 1)$. Actually, the 
selection of parameters $k$ and $A$ has been discussed carefully in 
\cite{WENO-IM}, and $k = 2, A = 0.1$ was recommended.

\begin{lemma} 
The mapping function $\big( g^{\mathrm{IM}} \big)_{s}(\omega; k, A)$ 
defined by Eq.(\ref{mappingFunctionWENO-IM}) satisfies: \\

C1. $\big( g^{\mathrm{IM}} \big)'_{s}(\omega; k, A) \geq 0, \omega 
\in [0, 1]$;

C2. $\big( g^{\mathrm{IM}} \big)_{s}(0; k, A) = 0, 
\big( g^{\mathrm{IM}} \big)_{s}(1; k, A) = 1$; 

C2. $\big( g^{\mathrm{IM}} \big)_{s}(d_{s}; k, A) = d_{s}$; 

C3. $\big( g^{\mathrm{IM}} \big)'_{s}(d_{s}; k, A) = \cdots = \big( g
^{\mathrm{IM}} \big)^{(k)}_{s}(d_{s}; k, A) = 0, \big( g^{\mathrm{IM}
} \big)^{(k + 1)}_{s}(d_{s}; k, A) \neq 0$.
\label{lemmaWENO-IMproperties}
\end{lemma}

We refer to \cite{WENO-IM} for the detailed proof of Lemma 
\ref{lemmaWENO-IMproperties}.



\section{Analysis of the nonlinear weights of the existing mapped 
WENO schemes}
\label{analysis_nonlinear_weights}
\subsection{Monotone increasing piecewise mapping function and the 
generalized WENO-ACM schemes}\label{subsec:MIP-WENO-ACMk}

In \cite{WENO-ACM}, the present authors have proposed the 
fifth-order WENO-ACM scheme. It has been demonstrated that taking 
narrower transition intervals (standing for the intervals over which 
the mapping results are in a transition from $0$ to $d_{s}$ or from 
$d_{s}$ to $1$) of the mapping function does not bring any adverse 
effects on the resolutions and convergence orders, and the 
associated scheme still performs very well even if the transition 
intervals are infinitely close to $0$. Therefore, we can set the 
transition intervals to be $0$ leading to a simpler form of the 
mapping function as follows
\begin{equation}
\big( g^{\mathrm{MIP-ACM}} \big)_{s}(\omega) = \left\{
\begin{array}{ll}
0, & \omega \in \Omega_{1} = [0, \mathrm{CFS}_{s}), \\
d_{s}, & \omega \in \Omega_{2} = [\mathrm{CFS}_{s},\overline{
\mathrm{CFS}}_{s}],\\
1, & \omega \in \Omega_{3} = (\overline{\mathrm{CFS}}_{s},1],
\end{array}
\right.
\label{eq:mappingFunctionMIP-ACM}
\end{equation}
where $\mathrm{CFS}_{s}$ is the same as that in 
\cite{WENO-MAIMi,WENO-ACM} satisfying $\mathrm{CFS}_{s}\in (0,d_{s})$
, and $\overline{\mathrm{CFS}}_{s}=1- \frac{1-d_{s}}{d_{s}}\times 
\mathrm{CFS}_{s}$ with $\overline{\mathrm{CFS}}_{s} \in (d_{s}, 1)$. 
Clearly, $\big( g^{\mathrm{MIP-ACM}} \big)_{s}(\omega)$ is a 
discontinuous function with two jump discontinuities at $\omega = 
\mathrm{CFS}_{s}$ and $\omega = \overline{\mathrm{CFS}}_{s}$ on the 
interval $[0,1]$, while differentiable mapping functions on the 
interval $[0,1]$ were required in previously published mapped WENO 
schemes \cite{WENO-M,WENO-PM,WENO-IM,WENO-RM260,WENO-PPM5,
WENO-RM-Vevek2018,WENO-AIM,WENO-MAIMi,WENO-ACM}. However, after 
extensive numerical tests, we find that a continuous mapping 
function is not essential in the design of the mapped WENO scheme. 
Actually, in the evaluation at $\omega_{s}^{\mathrm{JS}}$ of the 
Taylor series approximations of the mapping function about the 
optimal weights $d_{s}$, which plays the core role in the 
convergence analysis of the mapped WENO schemes (originally proposed 
by Henrick in the statement of page 556 in \cite{WENO-M}), one needs 
only the mapping function to be differentiable near the neighborhood 
of $\omega =d_{s}$ but not over the whole range of $\omega \in[0,1]$
. Therefore, we innovatively propose the definition of the
\textit{monotone increasing piecewise} mapping function. 
\begin{definition}(monotone increasing piecewise mapping function)
Let $\Omega = [0,1]$, and assume that $\Omega$ is divided into a 
sequence of nonoverlapping intervals $\Omega_{i}, i=1,2,\cdots,M$, 
that is, $\Omega = \Omega_{1}\cup\Omega_{2}\cup\cdots\cup\Omega_{M}$ 
and $\Omega_{i}\cap\Omega_{j}=\varnothing$, for $\forall i,j = 1,2,
\cdots, M$ and $i\neq j$. Let $\overline{\Omega}_{i} = \{\omega 
\lvert \omega \in \Omega_{i}$ and $\omega \notin \partial \Omega_{i}
\}$, and suppose that $\big( g^{\mathrm{MIP-X}} \big)_{s}(\omega)$ 
is a mapping function on the interval $[0,1]$, then $\big( 
g^{\mathrm{MIP-X}} \big)_{s}(\omega)$ is called a \textbf{monotone
increasing piecewise} mapping function, if it satisfies the following
conditions: (C1) for $\forall \omega \in \overline{\Omega}_{i}, i= 1,
\cdots, M$, $\big( g^{\mathrm{MIP-X}} \big)_{s}(\omega)$ is 
differentiable and $\big(g^{\mathrm{MIP-X}}\big)_{s}'(\omega)\geq 0$
; (C2) for $\forall \omega_{i}, \omega_{j} \in \Omega$, if
$\omega_{i} \geq \omega_{j}$, then $\big(g^{\mathrm{MIP-X}} \big)_{s}
(\omega_{i}) \geq \big(g^{\mathrm{MIP-X}}\big)_{s}(\omega_{j})$.
\label{def:MIP-mapping}
\end{definition}

It is trivial to verify that $\big( g^{\mathrm{MIP-ACM}} \big)_{s}(
\omega)$ defined by Eq.(\ref{eq:mappingFunctionMIP-ACM}) is a 
monotone increasing piecewise mapping function and it satisfies the 
following properties.

\begin{lemma}
The mapping function $\big(g^{\mathrm{MIP-ACM}}\big)_{s}(\omega)$ 
defined by Eq.(\ref{eq:mappingFunctionMIP-ACM}) satisfies the
following properties: \\

C1. for $\forall \omega \in \overline{\Omega}_{i}, i=1,2,3$,
$\big( g^{\mathrm{MIP-ACM}}\big)_{s}'(\omega)\geq 0$; 

C2. for $\forall \omega \in \Omega$, $0\leq\big( g^{\mathrm{MIP-ACM}}
\big)_{s}(\omega)\leq 1$;

C3. $d_{s} \in \Omega_{2}$, $\big( g^{\mathrm{MIP-ACM}}\big)_{s}(
d_{s}) = d_{s}, \big( g^{\mathrm{MIP-ACM}}\big)'_{s}(d_{s}) = 
\big( g^{\mathrm{MIP-ACM}}\big)''_{s}(d_{s}) = \cdots = 0$;

C4. $\big( g^{\mathrm{MIP-ACM}}\big)_{s}(0) = 0, 
\big(g^{\mathrm{MIP-ACM}}\big)_{s}(1) = 1, 
\big( g^{\mathrm{MIP-ACM}}\big)'_{s}(0^{+}) = \big( 
g^{\mathrm{MIP-ACM}}\big)'_{s}(1^{-}) = 0$.  
\label{lemma:mappingFunction:WENO-ACM}
\end{lemma}

Naturally, we can derive a generalized version of the mapping 
function $\big( g^{\mathrm{MIP-ACM}} \big)_{s}(\omega)$ as follows
\begin{equation}
\big( g^{\mathrm{MIP-ACM}k} \big)_{s}(\omega) = \left\{
\begin{array}{ll}
k_{s} \omega, & \omega \in \Omega_{1}, \\
d_{s}, & \omega \in \Omega_{2}, \\
1 - k_{s} (1 - \omega), & \omega \in \Omega_{3},
\end{array}
\right.
\label{eq:mappingFunctionMIP-ACMk}
\end{equation}
where $\Omega_{1},\Omega_{2},\Omega_{3}$ are the same as in 
Eq.(\ref{eq:mappingFunctionMIP-ACM}) and $k_{s} \in \Big[0, 
\frac{d_{s}}{\mathrm{CFS}_{s}}\Big]$. Obviously, if $k_{s}$ is taken 
to be $0$, $\big( g^{\mathrm{MIP-ACM}k} \big)_{s}(\omega)$ exactly 
turns into $\big( g^{\mathrm{MIP-ACM}} \big)_{s}(\omega)$. Thus, we 
need only discuss $\big( g^{\mathrm{MIP-ACM}k} \big)_{s}(\omega)$.
Similarly, it is easy to know that $\big(g^{\mathrm{MIP-ACM}k} 
\big)_{s}(\omega)$ is a monotone increasing piecewise mapping
function and it satisfies the following properties.

\begin{lemma}
The mapping function $\big(g^{\mathrm{MIP-ACM}k}\big)_{s}(\omega)$ 
defined by Eq.(\ref{eq:mappingFunctionMIP-ACMk}) satisfies the
following properties: \\

C1. for $\forall \omega \in \overline{\Omega}_{i}, i=1,2,3, \big( 
g^{\mathrm{MIP-ACM}k}\big)_{s}'(\omega)\geq 0$; 

C2. for $\forall \omega \in \Omega, 0\leq\big( g^{\mathrm{MIP-ACM}k}
\big)_{s}(\omega)\leq 1$;

C3. $d_{s} \in \Omega_{2}$, $\big( g^{\mathrm{MIP-ACM}k}\big)_{s}(
d_{s}) = d_{s}, \big( g^{\mathrm{MIP-ACM}k}\big)'_{s}(d_{s}) = 
\big( g^{\mathrm{MIP-ACM}k}\big)''_{s}(d_{s}) = \cdots = 0$;

C4. $\big( g^{\mathrm{MIP-ACM}k}\big)_{s}(0) = 0, 
\big( g^{\mathrm{MIP-ACM}k}\big)_{s}(1) = 1, 
\big( g^{\mathrm{MIP-ACM}k}\big)'_{s}(0^{+}) = 
\big( g^{\mathrm{MIP-ACM}k}\big)'_{s}(1^{-}) = k_{s}$.  
\label{lemma:mappingFunction:WENO-ACMk}
\end{lemma}

As the proofs of Lemma \ref{lemma:mappingFunction:WENO-ACM} and Lemma
\ref{lemma:mappingFunction:WENO-ACMk} are very easy, we do not state
them here and we can observe these properties intuitively from the
$\big(g^{\mathrm{MIP-ACM}k}\big)_{s}(\omega)\sim \omega$ curves as 
shown in Fig. \ref{fig:gOmega:analysis:2} below.

Now, we give the monotone increasing piecewise
approximate-constant-mapped WENO scheme, denoted as MIP-WENO-ACM$k$
, with the mapped weights as follows
\begin{equation}
\omega_{s}^{\mathrm{MIP-ACM}k} = \dfrac{\alpha _{s}^{\mathrm{MIP-ACM}
k}}{\sum_{l = 0}^{2} \alpha _{l}^{\mathrm{MIP-ACM}k}}, \alpha_{s}^{
\mathrm{MIP-ACM}k}=\big( g^{\mathrm{MIP-ACM}k} \big)_{s}(\omega^{
\mathrm{JS}}_{s}).
\end{equation}

We present Theorem \ref{theorem:convergenceRates}, which will show 
that the MIP-WENO-ACM$k$ scheme can recover the optimal convergence 
orders for different values of $n_{\mathrm{cp}}$ in smooth regions.

\begin{theorem}
When $\mathrm{CFS}_{s} \ll d_{s}$, for $\forall n_{\mathrm{cp}} < 
r - 1$, the $(2r-1)$th-order MIP-WENO-ACM$k$ scheme can achieve the 
optimal convergence rates of accuracy if the new mapping function
$\big( g^{\mathrm{MIP-ACM}k} \big)_{s}(\omega)$ is applied to the 
weights of the $(2r-1)$th-order WENO-JS scheme.
\label{theorem:convergenceRates:MIP}
\end{theorem}

We can prove Theorem \ref{theorem:convergenceRates:MIP} by employing 
the Taylor series analysis and using Lemma 
\ref{lemma:mappingFunction:WENO-ACMk} of this paper and Lemma 1 and 
Lemma 2 in the statement of page 456 to 457 in \cite{WENO-IM}, and 
the detailed proof process is almost identical to the one in 
\cite{WENO-M}.

\subsection{Discussion about the effects of $g'(0)$ of the mapped 
WENO schemes on resolutions and spurious oscillations}
\label{subsec:importantDiscussion}
In order to study the effects of $g'(0)$ of the mapped WENO schemes 
on resolutions and spurious oscillations in simulating the problems 
with discontinuities, we calculate the one-dimensional linear 
advection equation 
\begin{equation}
u_{t} + u_{x} = 0,
\label{eq:LAE} 
\end{equation}
with the following initial condition \cite{WENO-JS}
\begin{equation}
\begin{array}{l}
u(x, 0) = \left\{
\begin{array}{ll}
\dfrac{1}{6}\big[ G(x, \beta, z - \hat{\delta}) + 4G(x, \beta, z) + G
(x, \beta, z + \hat{\delta}) \big], & x \in [-0.8, -0.6], \\
1, & x \in [-0.4, -0.2], \\
1 - \big\lvert 10(x - 0.1) \big\rvert, & x \in [0.0, 0.2], \\
\dfrac{1}{6}\big[ F(x, \alpha, a - \hat{\delta}) + 4F(x, \alpha, a) +
 F(x, \alpha, a + \hat{\delta}) \big], & x \in [0.4, 0.6], \\
0, & \mathrm{otherwise},
\end{array}\right. 
\end{array}
\label{eq:LAE:SLP}
\end{equation}
where $G(x, \beta, z) = \mathrm{e}^{-\beta (x - z)^{2}}, F(x, \alpha
, a) = \sqrt{\max \big(1 - \alpha ^{2}(x - a)^{2}, 0 \big)}$, and 
the constants are $z = -0.7, \hat{\delta} = 0.005, \beta = \dfrac{
\log 2}{36\hat{\delta} ^{2}}, a = 0.5$ and $\alpha = 10$. The 
periodic boundary condition is used in the two directions and the 
CFL number is set to be $0.1$. This problem consists of a Gaussian,
a square wave, a sharp triangle and a semi-ellipse. For brevity in 
the presentation, we call this \textit{\textbf{L}inear 
\textbf{P}roblem} SLP as it is presented by \textit{\textbf{S}hu} et 
al. in \cite{WENO-JS}. 

The following two groups of mapped WENO schemes with various values 
of $\big(g^{\mathrm{X}}\big)'(0)$ (X stands for some specific mapped 
WENO scheme) are employed in the discussion.

\subsubsection{Study on the WENO-PM6 and WENO-IM($k,A$) schemes}
In this subsection, we focus on the performances of the WENO-PM6 
scheme \cite{WENO-PM} and the WENO-IM($k,A$) schemes \cite{WENO-IM}
with $k=2$ and $A=0.1, 0.5$ on solving SLP. A uniform mesh size of 
$N = 400$ is used and the output time is set to be $t = 200$. 

\begin{figure}[ht]
\centering
\includegraphics[height=0.27\textwidth]
{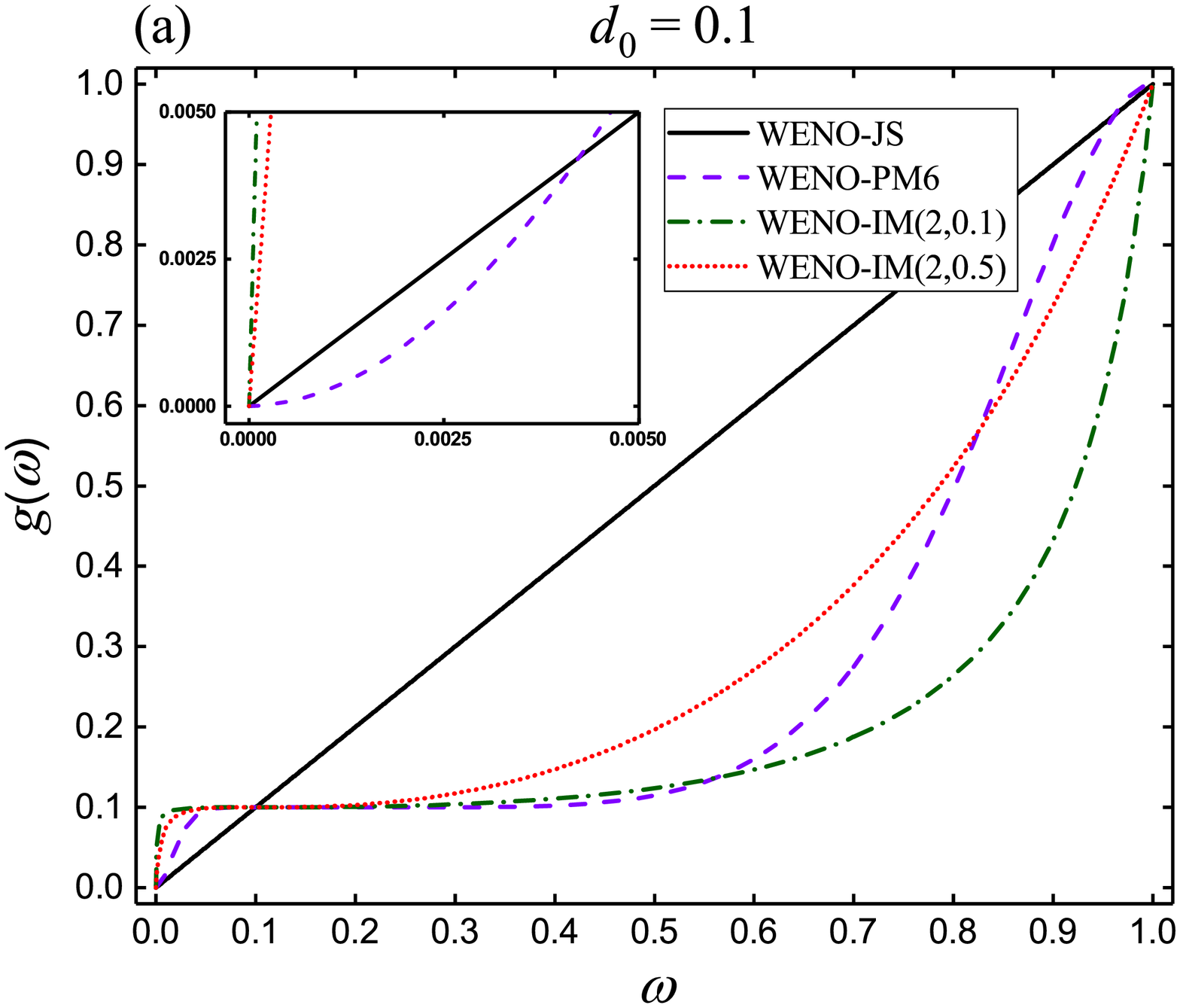}\quad
\includegraphics[height=0.27\textwidth]
{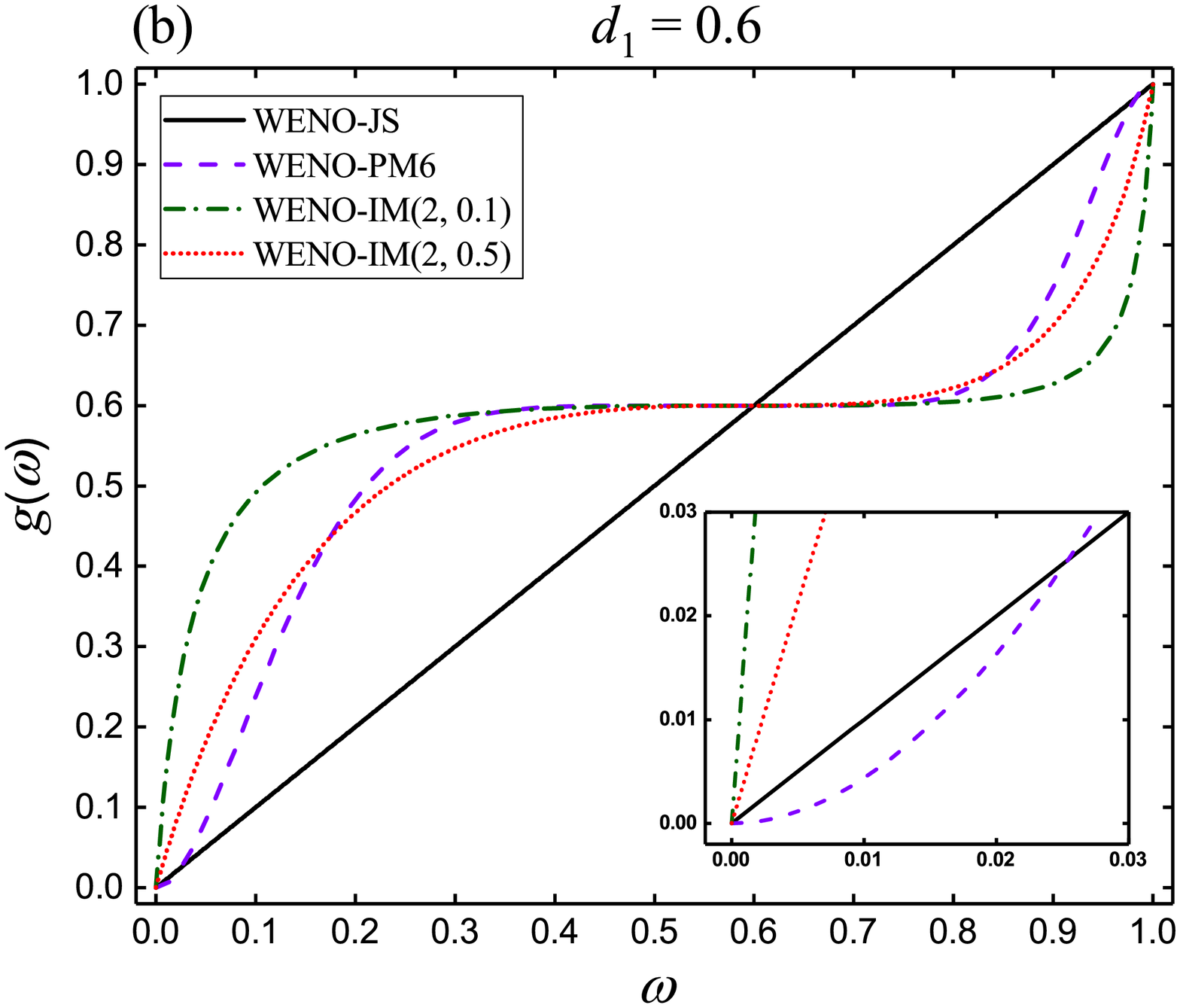}
\includegraphics[height=0.27\textwidth]
{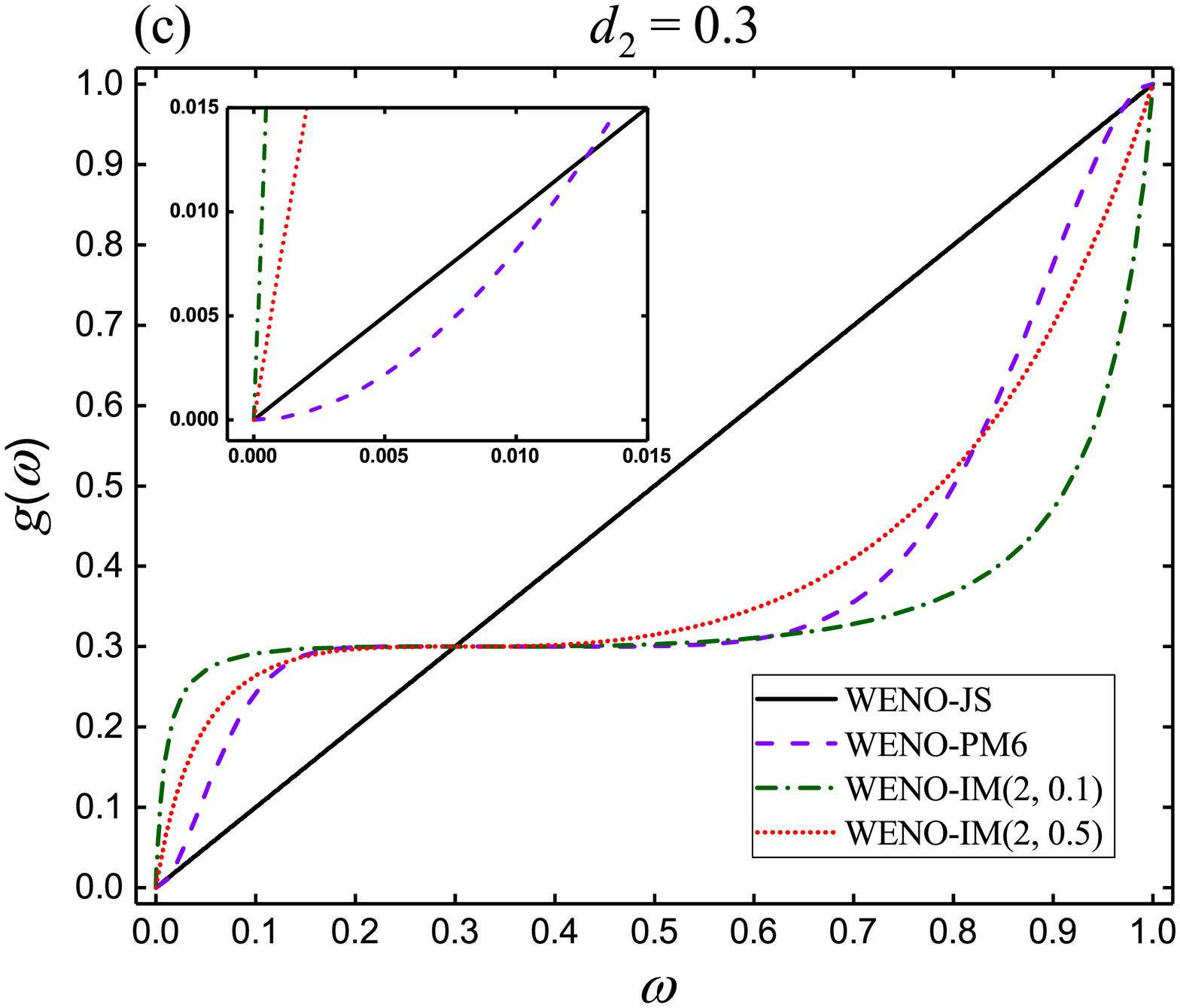}
\caption{The mapping functions of the WENO-PM6, WENO-IM($2,0.1$) and 
WENO-IM($2,0.5$) schemes, $d_{0}=0.1, d_{1} = 0.6, d_{2} = 0.3$.}
\label{fig:gOmega:analysis:1}
\end{figure}

From Fig. \ref{fig:gOmega:analysis:1}, we can easily see that the 
values of $\big(g^{\mathrm{X}}\big)'(0)$ satisfy 
\begin{equation}
\big(g^{\mathrm{PM}6}\big)'(0) < \big(g^{\mathrm{IM}(2,0.5)}\big)'(0)
< \big(g^{\mathrm{IM}(2,0.1)} \big)'(0).
\label{eq:gOmega_0:analysis:1}
\end{equation}

\begin{figure}[ht]
\centering
\includegraphics[height=0.5\textwidth]
{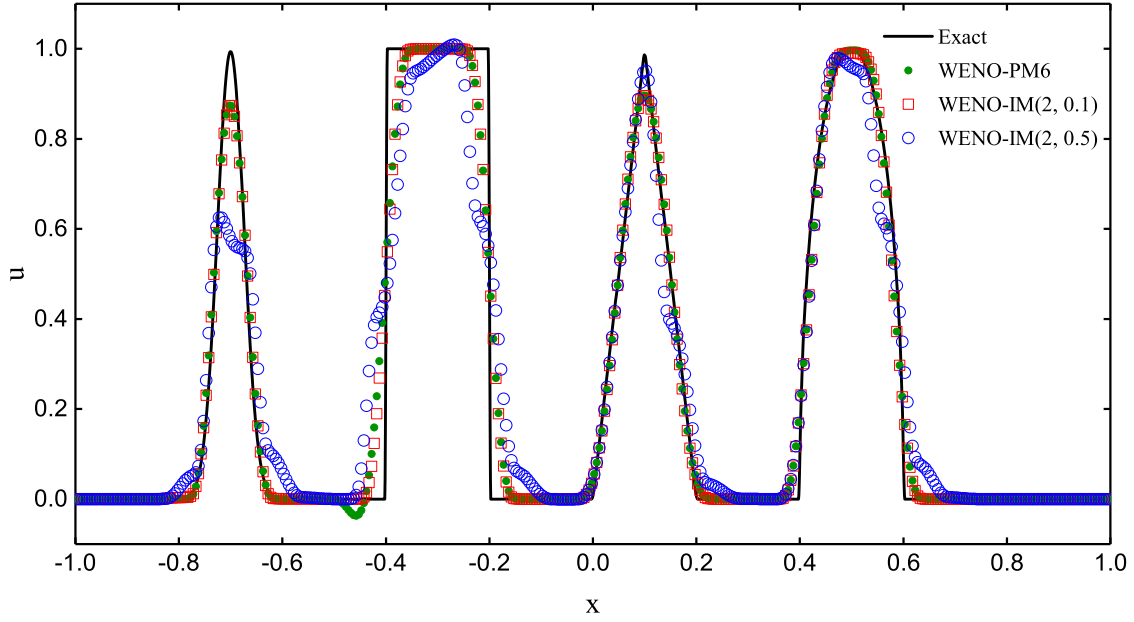}
\caption{Performance of the fifth-order WENO-PM6, WENO-IM($2,0.1$) 
and WENO-IM($2,0.5$) schemes for SLP with $N=400$ at long output 
time $t=200$.}
\label{fig:x-u:analysis:1}
\end{figure}

Fig. \ref{fig:x-u:analysis:1} shows the calculating results, and 
Table \ref{table:errors:analysis:1} shows the $L_{1}, L_{2}$, 
$L_{\infty}$ errors and the order of these errors (in brackets in 
descending manner), i.e., in the second column, $1$ indicates the 
largest $L_{1}$ error and $2$ indicates the second largest one, etc. 
From Fig. \ref{fig:x-u:analysis:1} and Table 
\ref{table:errors:analysis:1}, we can observe that: 
(1) the WENO-IM($2,0.1$) scheme, whose 
$\big(g^{\mathrm{IM}(2,0.1)}\big)'(0)$ is the largest, presents the 
smallest spurious oscillation (actually, there is no spurious 
oscillation in present computing conditions) and gives the smallest 
$L_{1}, L_{2}$ and $L_{\infty}$ errors; (2) the WENO-PM6 scheme, 
whose $\big(g^{\mathrm{PM6}}\big)'(0)$ is the smallest and satisfies
$\big(g^{\mathrm{PM6}}\big)'(0) = 0$, presents the largest spurious 
oscillation and gives the second largest $L_{1}, L_{2}$ and 
$L_{\infty}$ errors; (3) the WENO-IM($2,0.5$) scheme shows the 
lowest resolutions and presents evident spurious oscillation at the 
top of the square wave, and it gives the largest $L_{1}, L_{2}$ and 
$L_{\infty}$ errors, while its $\big(g^{\mathrm{IM}(2,0.5)}\big)'(0)$
is neither largest nor smallest.

\begin{table}[ht]
\centering
\caption{The $L_{1}, L_{2}, L_{\infty}$ errors for SLP with $N=400$ 
at long output time $t=200$, computed by the WENO-PM6, 
WENO-IM($2,0.1$) and WENO-IM($2,0.5$) schemes.}
\label{table:errors:analysis:1}
\begin{tabular*}{\hsize}
{@{}@{\extracolsep{\fill}}llll@{}}
\hline
Schemes  & $L_{1}$ error  & $L_{2}$ error   & $L_{\infty}$ error  \\ 
\hline
WENO-PM6         & 5.69929e-02(2) & 1.06646e-01(2) & 4.80453e-01(2)\\
WENO-IM($2,0.1$) & 5.46038e-02(3) & 1.04363e-01(3) & 4.52252e-01(3)\\
WENO-IM($2,0.5$) & 1.23950e-01(1) & 1.75554e-01(1) & 5.24783e-01(1)\\
\hline
\end{tabular*}
\end{table}

\subsubsection{Study on the WENO-MAIM2 and MIP-WENO-ACM$k$ schemes}
In this subsection, we focus on the performances of the WENO-MAIM2 
\cite{WENO-MAIMi} and MIP-WENO-ACM$k$ schemes on solving SLP. We 
still use a uniform mesh size of $N = 400$ and choose the output 
time $t = 200$. 

\begin{figure}[ht]
\centering
\includegraphics[height=0.27\textwidth]
{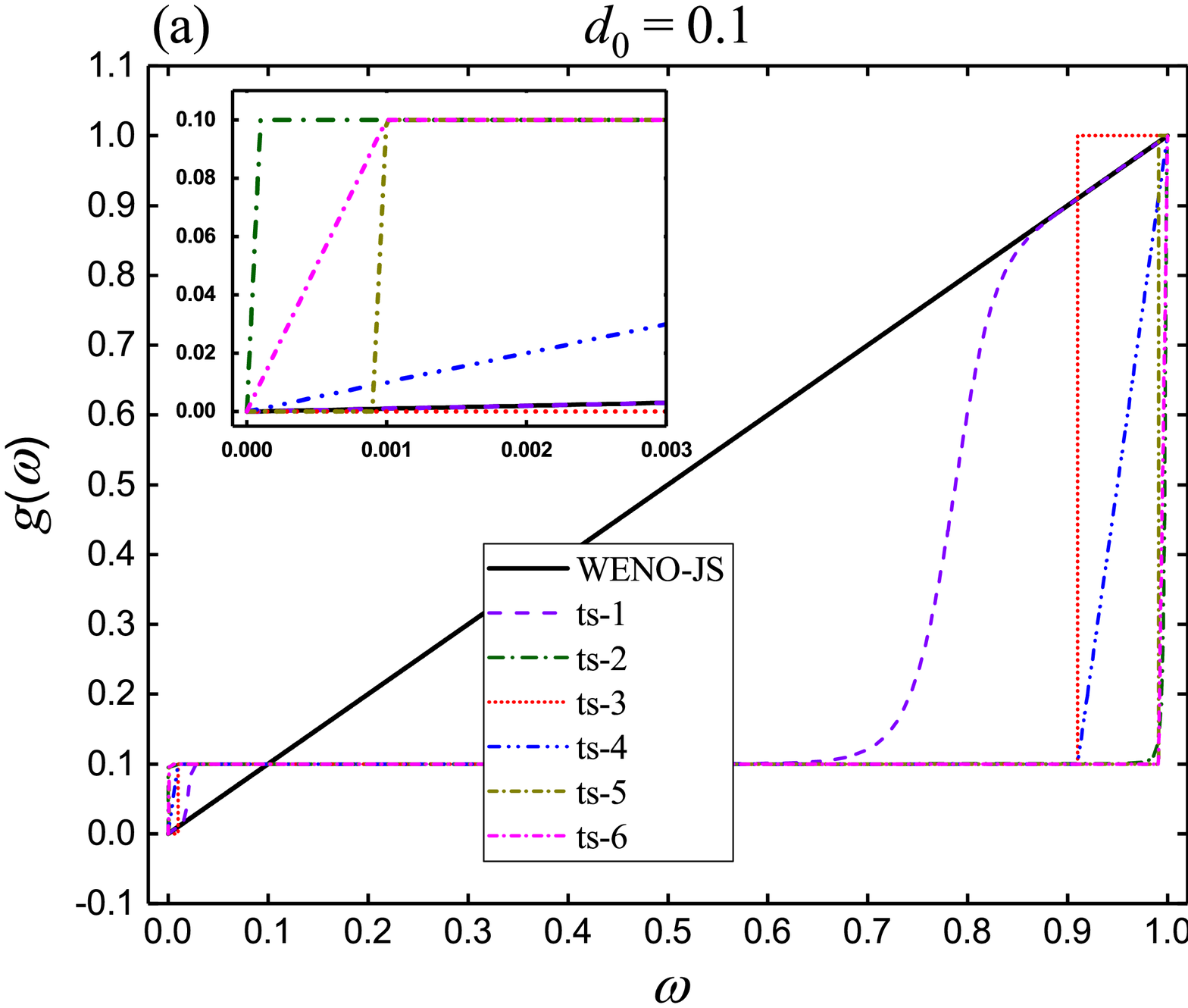}\quad
\includegraphics[height=0.27\textwidth]
{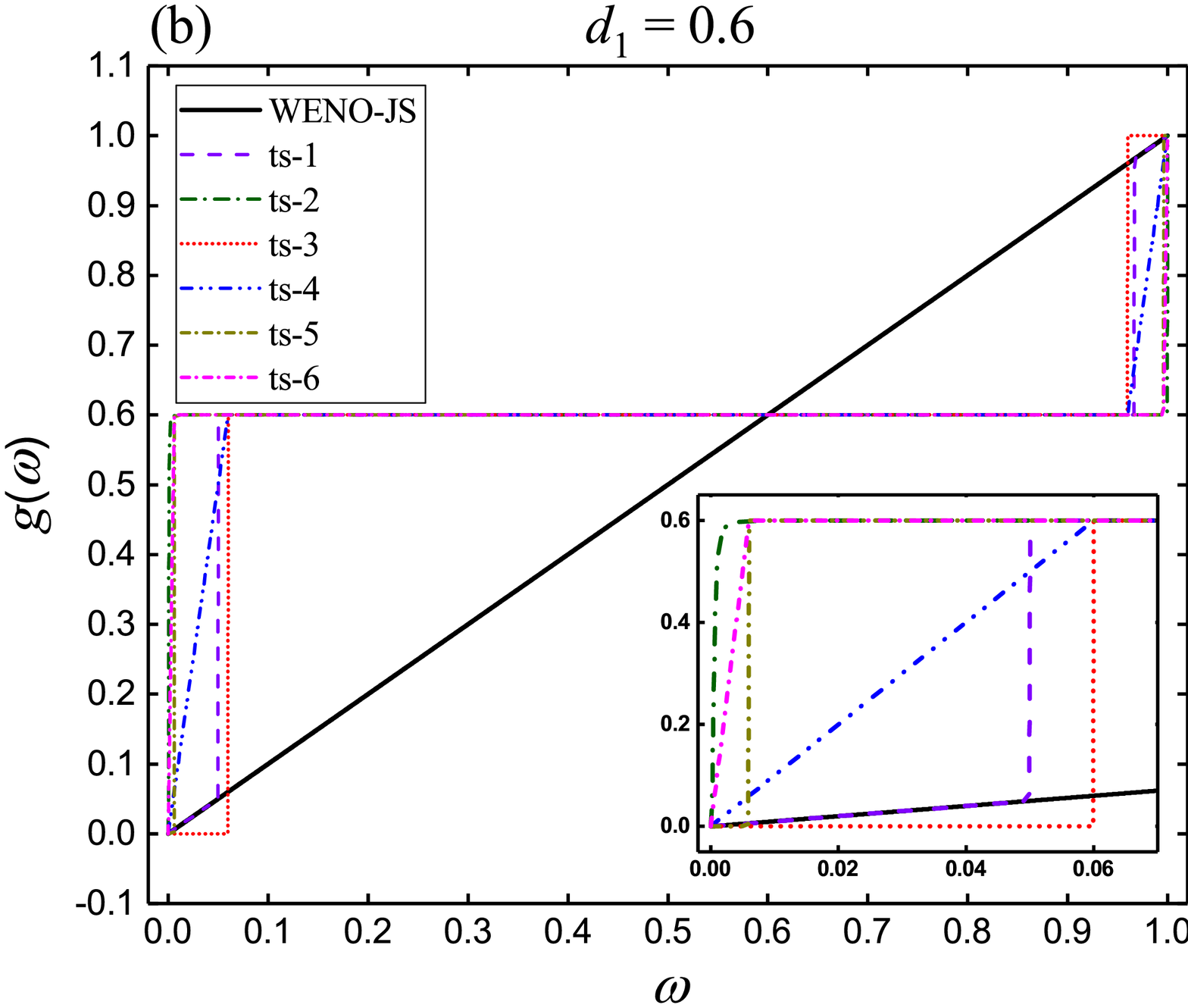}
\includegraphics[height=0.27\textwidth]
{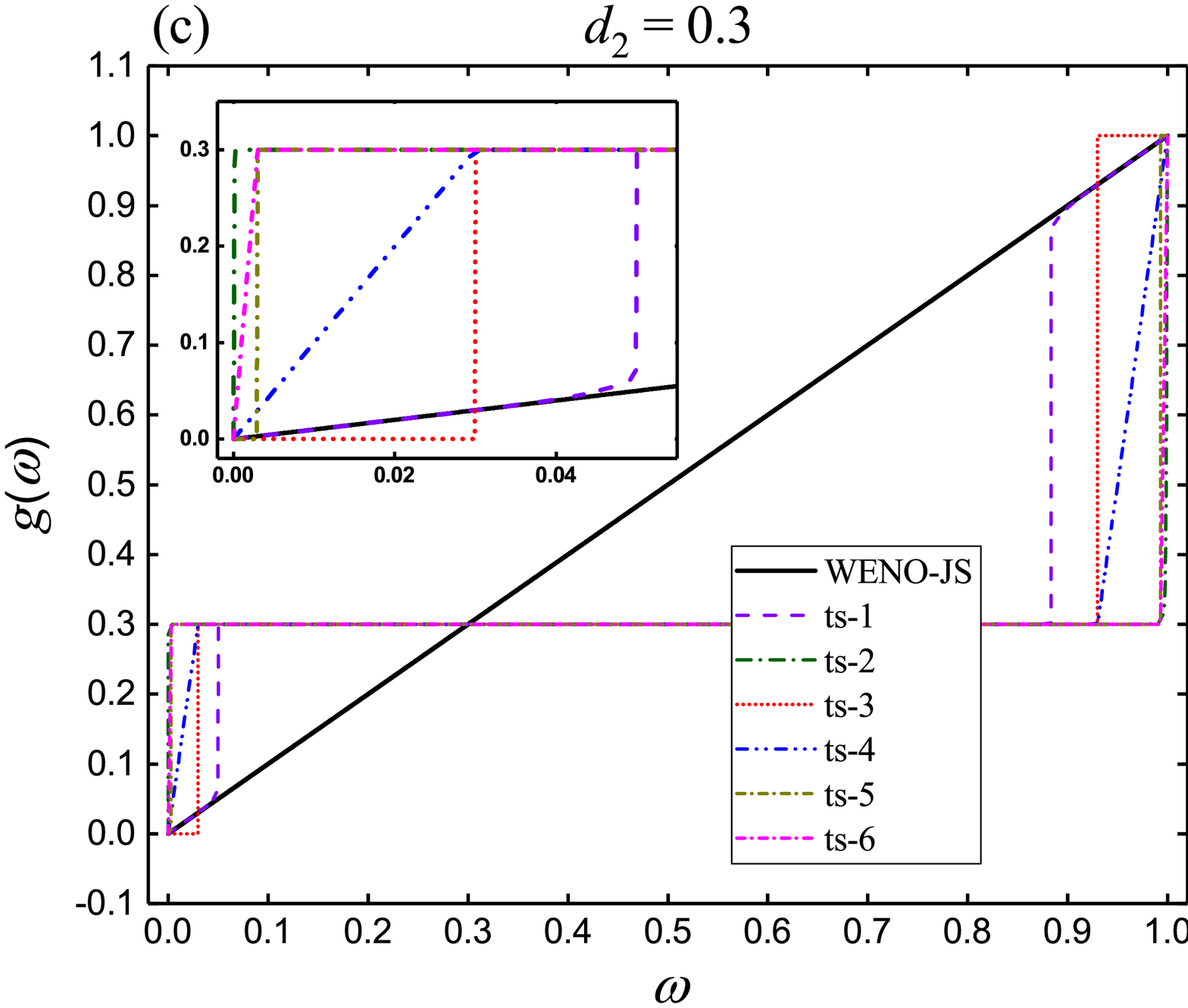}
\caption{The mapping functions of the test schemes shown in Table 
\ref{table:TestSchemes:analysis:2}, $d_{0}=0.1, d_{1} = 0.6, 
d_{2} = 0.3$.}
\label{fig:gOmega:analysis:2}
\end{figure}

As shown in Table \ref{table:TestSchemes:analysis:2}, $6$ different 
test schemes (ts-$i$,$i=1,\cdots,6$) of the WENO-MAIM2 and 
MIP-WENO-ACM$k$ schemes with specified parameters leading to various 
values of $\big(g^{\mathrm{ts-}i}\big)'(0)$ are used in the 
discussion. In Fig. \ref{fig:gOmega:analysis:2}, we plot the curves 
of $\big(g^{\mathrm{ts-}i}\big)(\omega)\sim \omega, i=1,\cdots,6$, 
and we can intuitively observe that the values of 
$\big(g^{\mathrm{ts-}i}\big)'(0)$ satisfy
\begin{equation}
\big(g^{\mathrm{ts-}3}\big)'(0) = \big(g^{\mathrm{ts-}5}\big)'(0)
< \big(g^{\mathrm{ts-}1}\big)'(0) < \big(g^{\mathrm{ts-}4}\big)'(0)
< \big(g^{\mathrm{ts-}6}\big)'(0) < \big(g^{\mathrm{ts-}2}\big)'(0).
\label{eq:gOmega_0:analysis:2}
\end{equation}

\begin{table}[ht]
\centering
\caption{$6$ different test schemes (ts-$1,\cdots,6$) of the 
WENO-MAIM2 and MIP-WENO-ACM$k$ schemes with specified parameters.}
\label{table:TestSchemes:analysis:2}
\begin{tabular*}{\hsize}
{@{}@{\extracolsep{\fill}}llll@{}}
\hline
ts-i  & WENO-X          & Parameters  
                        &  $\big(g^{\mathrm{X}}\big)'(0)$\\
\hline
ts-1  & WENO-MAIM2      & $k=10, A=1.0\mathrm{e-}6,Q=1.0,
\mathrm{CFS}_{s}=0.05$  & $1$\\
ts-2  & WENO-MAIM2      & $k=10, A=1.0\mathrm{e-}6,Q=0.25,
\mathrm{CFS}_{s}=0.05$  & $\gg 1$\\
ts-3  & MIP-WENO-ACM$k$ & $k_{s}=0,\mathrm{CFS}_{s}=0.1d_{s}$ 
                        & $0$\\
ts-4  & MIP-WENO-ACM$k$ & $k_{s}=10,\mathrm{CFS}_{s}=0.1d_{s}$ 
                        & $10$\\
ts-5  & MIP-WENO-ACM$k$ & $k_{s}=0,\mathrm{CFS}_{s}=0.01d_{s}$ 
                        & $0$\\
ts-6  & MIP-WENO-ACM$k$ & $k_{s}=100,\mathrm{CFS}_{s}=0.01d_{s}$ 
                        & $100$\\
\hline
\end{tabular*}
\end{table}

\begin{figure}[ht]
\centering
\includegraphics[height=0.5\textwidth]
{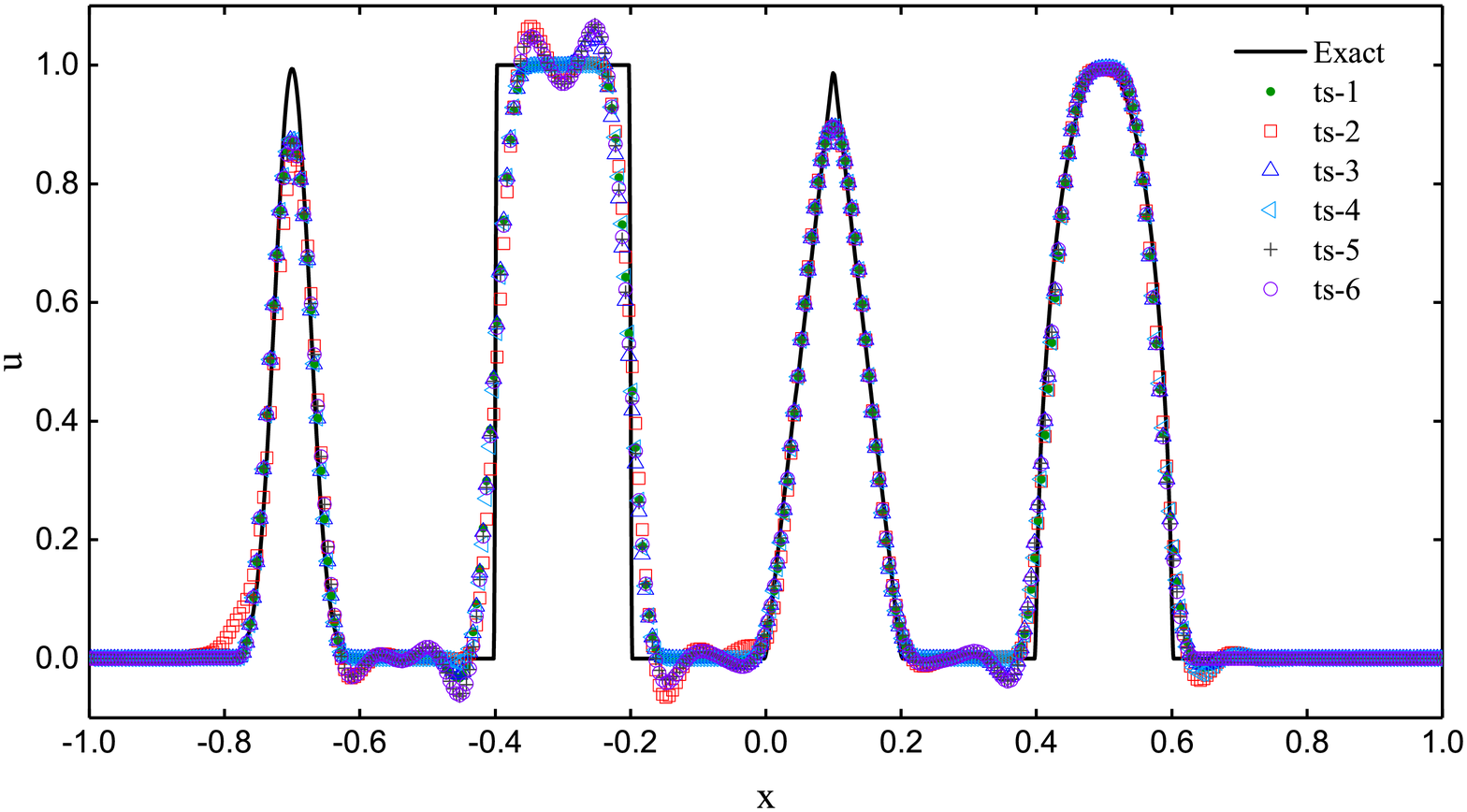}
\caption{Performance of the test schemes shown in Table 
\ref{table:TestSchemes:analysis:2} for the SLP with $N=400$ at long 
output time $t=200$.}
\label{fig:x-u:analysis:2}
\end{figure}

Fig. \ref{fig:x-u:analysis:2} shows the calculating results, and 
Table \ref{table:errors:analysis:2} shows the $L_{1}, L_{2}, L_{
\infty}$ errors and the order of these errors (in brackets in 
descending manner). From 
Fig. \ref{fig:x-u:analysis:2} and Table \ref{table:errors:analysis:2}
, we can observe that: (1) all the $6$ test schemes present spurious 
oscillations; (2) ts-$2$ shows more in number and bigger in size of 
the spurious oscillations than ts-$1$, and the $L_{1}, L_{2}, 
L_{\infty}$ errors of ts-$2$ are larger than those of ts-$1$; (3) 
however, although $\big(g^{\mathrm{ts-}4}\big)'(0) > \big(g^{\mathrm{
ts-}3}\big)'(0)$ and $\big(g^{\mathrm{ts-}3}\big)'(0) = 0$, ts-$4$ 
shows fewer in number and smaller in size of the spurious 
oscillations than ts-$3$, and the $L_{1}, L_{2}, L_{\infty}$ 
errors of ts-$4$ are smaller than those of ts-$3$; (4) in addition, 
although $\big(g^{\mathrm{ts-}5}\big)'(0) \ll \big(g^{\mathrm{ts-}6}
\big)'(0)$ and $\big(g^{\mathrm{ts-}5}\big)'(0) = 0$, ts-$5$ shows 
comparable spurious oscillations both in number and in size with 
ts-$6$, and the $L_{1}, L_{2}, L_{\infty}$ errors of ts-$5$ are very 
close to those of ts-$6$, or in other words, the $L_{1}$ error of 
ts-$5$ is slightly smaller than that of ts-$6$, while the $L_{2},
L_{\infty}$ errors of ts-$5$ are slightly larger than those of 
ts-$6$.

\begin{table}[ht]
\begin{scriptsize}
\centering
\caption{The $L_{1}, L_{2}, L_{\infty}$ errors for the SLP with 
$N=400$ at long output time $t=200$, computed by the test schemes 
shown in Table \ref{table:TestSchemes:analysis:2}.}
\label{table:errors:analysis:2}
\begin{tabular*}{\hsize}
{@{}@{\extracolsep{\fill}}llll@{}}
\hline
Schemes, ts-$i$ 
         & $L_{1}$ error  & $L_{2}$ error  & $L_{\infty}$ error  \\ 
\hline
ts-1     & 5.71367e-02(5) & 1.06259e-01(5) & 4.76278e-01(3)\\
ts-2     & 6.80647e-02(1) & 1.09276e-01(1) & 4.91997e-01(1)\\
ts-3     & 5.91473e-02(4) & 1.07220e-01(4) & 4.89545e-01(2)\\
ts-4     & 5.55635e-02(6) & 1.04530e-01(6) & 4.52208e-01(6)\\
ts-5     & 6.60760e-02(3) & 1.08378e-01(2) & 4.74908e-01(4)\\
ts-6     & 6.62020e-02(2) & 1.08289e-01(3) & 4.69467e-01(5)\\
\hline
\end{tabular*}
\end{scriptsize}
\end{table}

\subsection{Analysis of the real-time mapping relationship}
\subsubsection{Definition of order-preserving/non-order-preserving 
mapping and important numerical experiments}
From the discussion above, we can conclude that it is not essential 
to prevent the corresponding mapped WENO scheme from generating 
spurious oscillations or causing potential loss of accuracy near 
discontinuities that the first derivatives of the mapping functions 
tend to $0$ or a small value when $\omega$ is close to $0$. In this 
subsection, to discover the essential cause of the spurious 
oscillation generation and potential loss of 
accuracy, we will make a further analysis of the real-time mapping 
relationship $\big(g^{\mathrm{X}}\big)(\omega)\sim\omega$, that 
stands for the mapping relationship obtained from the calculation of 
some specific problem at specified output time but not directly 
obtained from the mapping function.

Before conducting the numerical experiments for the analysis, we
propose the definition of \textit{order-preserving} mapping and 
\textit{non-order-preserving} mapping.

\begin{definition}(order-preserving/non-order-preserving mapping) 
Suppose that $\big( g^{\mathrm{X}} \big)_{s}(\omega), s = 0,\cdots,
r-1$ is a monotone increasing piecewise mapping function of the 
$(2r-1)$th-order mapped WENO-X scheme. We say the set of mapping 
functions \Big\{$\big( g^{\mathrm{X}}\big)_{s}(\omega), s = 0,\cdots,
r-1$\Big\} is \textbf{order-preserving (OP)}, if for $\forall 
\omega_{a} \geq \omega_{b}$,
\begin{equation}
\big( g^{\mathrm{X}}\big)_{m}(\omega_{a}) \geq \big( 
g^{\mathrm{X}} \big)_{n}(\omega_{b}), \quad \forall m, n \in 
\{0, \cdots,r-1\},
\label{def:order-preserving_mappng}
\end{equation}
where the equality holds if and only if $\omega_{a} = \omega_{b}$. 
Otherwise, we say the set of mapping functions 
\Big\{$\big( g^{\mathrm{X}}\big)_{s}(\omega), s = 0,\cdots,r-1\Big\}$
is \textbf{non-order-preserving (non-OP)}.
\label{def:OPM}
\end{definition}

It is trivial to know that, even when the set of mapping functions 
\Big\{$\big( g^{\mathrm{X}}\big)_{s}(\omega), s = 0,\cdots,r-1$\Big\}
is \textit{non-OP}, Eq.(\ref{def:order-preserving_mappng}) may also 
hold at some points. Therefore, we add the following definition of 
\textit{OP point} and \textit{non-OP point}.

\begin{definition}(OP point, non-OP point) Let $S^{2r-1}$ denote the 
$(2r-1)$-point global stencil centered around $x_{j}$. Assume that 
$S^{2r-1}$ is subdivided into $r$-point substencils $\{S_{0},\cdots,
S_{r-1}\}$ and $\omega_{s}$ are the nonlinear weights corresponding 
to the substencils $S_{s}$ with $s=0,\cdots,r-1$, which are used as 
the independent variables by the mapping function. Suppose that 
$\big( g^{\mathrm{X}}\big)_{s}(\omega), s=0,\cdots,r-1$ is the 
mapping function of the mapped WENO-X scheme, then we say that a 
\textbf{non-OP} mapping process occurs at $x_{j}$, if $\exists m, n 
\in \{0,\cdots,r-1\}$, s.t.
\begin{equation}\left\{
\begin{array}{ll}
\begin{aligned}
&\big(\omega_{m} - \omega_{n}\big)\bigg(\big(g^{\mathrm{X}}\big)_{m}
(\omega_{m}) - \big(g^{\mathrm{X}}\big)_{n}(\omega_{n})\bigg)\leq 0, 
&\mathrm{if} \quad \omega_{m} \neq \omega_{n},\\
&\big(g^{\mathrm{X}}\big)_{m}(\omega_{m}) \neq \big(g^{\mathrm{X}}
\big)_{n}(\omega_{n}), &\mathrm{if} \quad \omega_{m}=\omega_{n}.
\end{aligned}
\end{array}\right.
\end{equation}
And we say $x_{j}$ is a \textbf{non-OP point}. Otherwise, we say 
$x_{j}$ is an \textbf{OP point}.
\label{def:MaRe}
\end{definition}

\begin{figure}[ht]
\centering
\includegraphics[height=0.35\textwidth]
{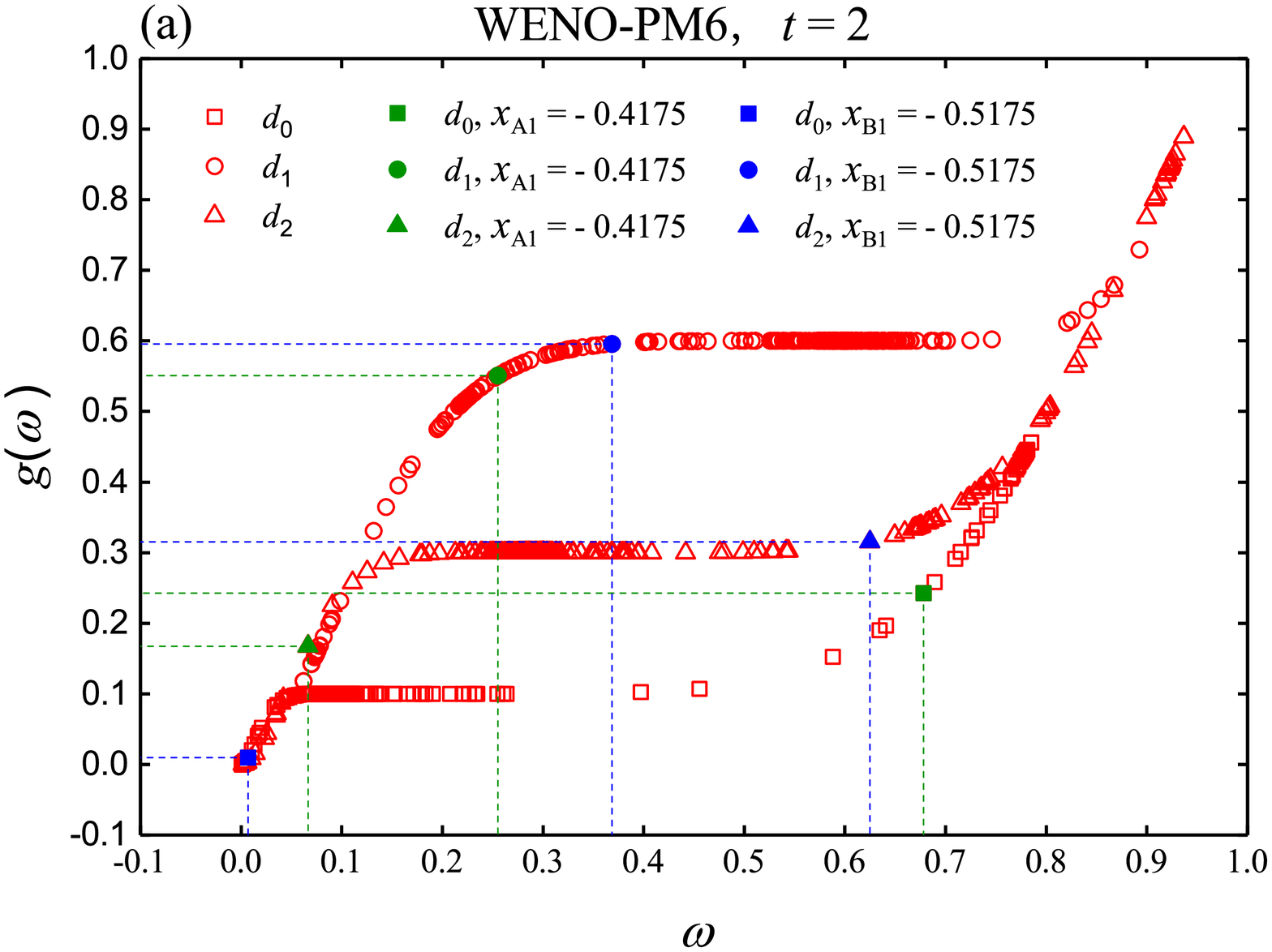}
\includegraphics[height=0.35\textwidth]
{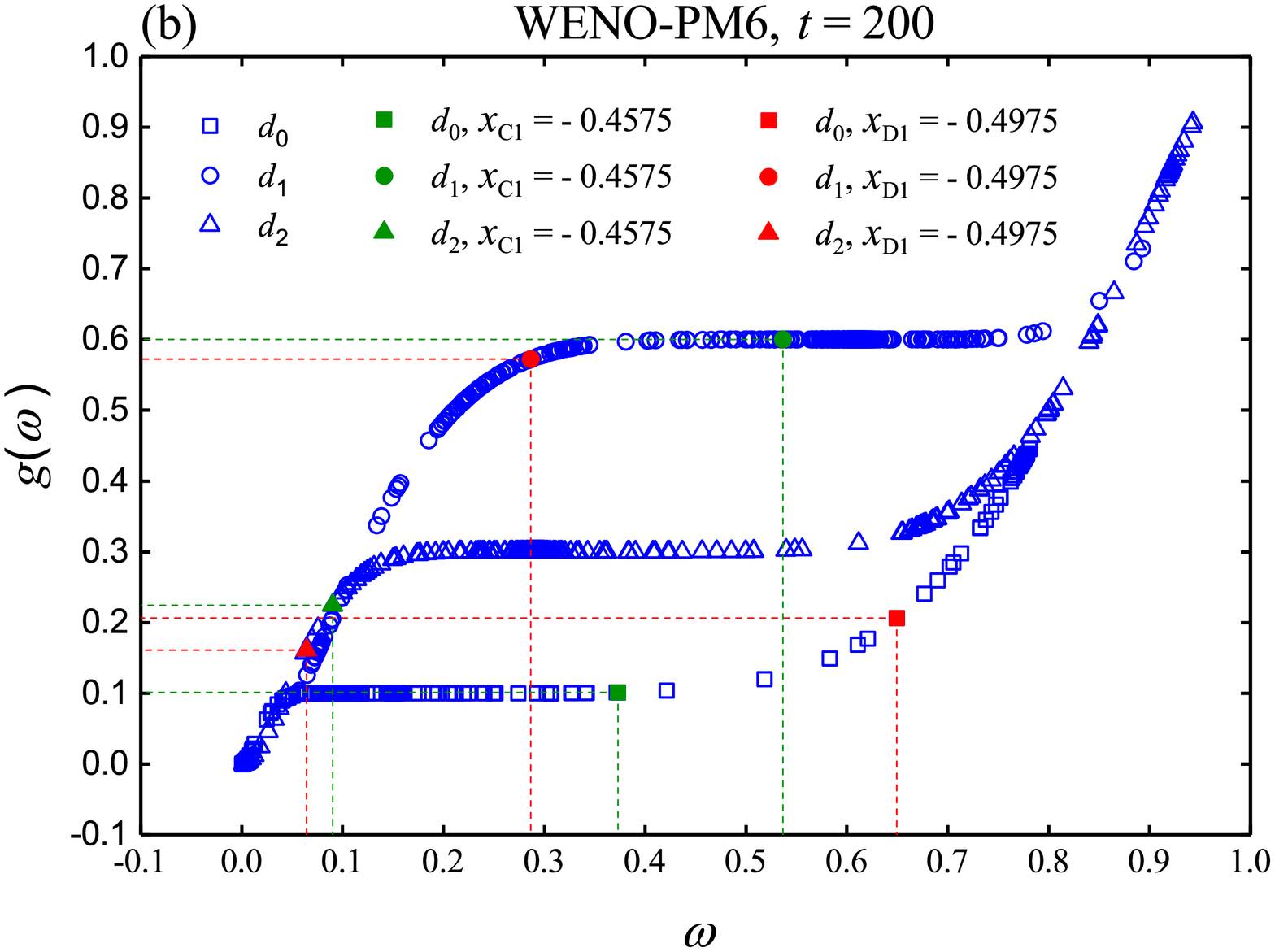}
\caption{The real-time mapping relationship $\big(g^{\mathrm{PM}6}
\big)(\omega)\sim\omega$ of the SLP. A uniform mesh size of $N = 400$
and two output times $t=2$ (left) and $t=200$ (right) are used.}
\label{fig:RT-gOmega:deepAna:PM6}
\end{figure}

\begin{figure}[ht]
\centering
\includegraphics[height=0.35\textwidth]
{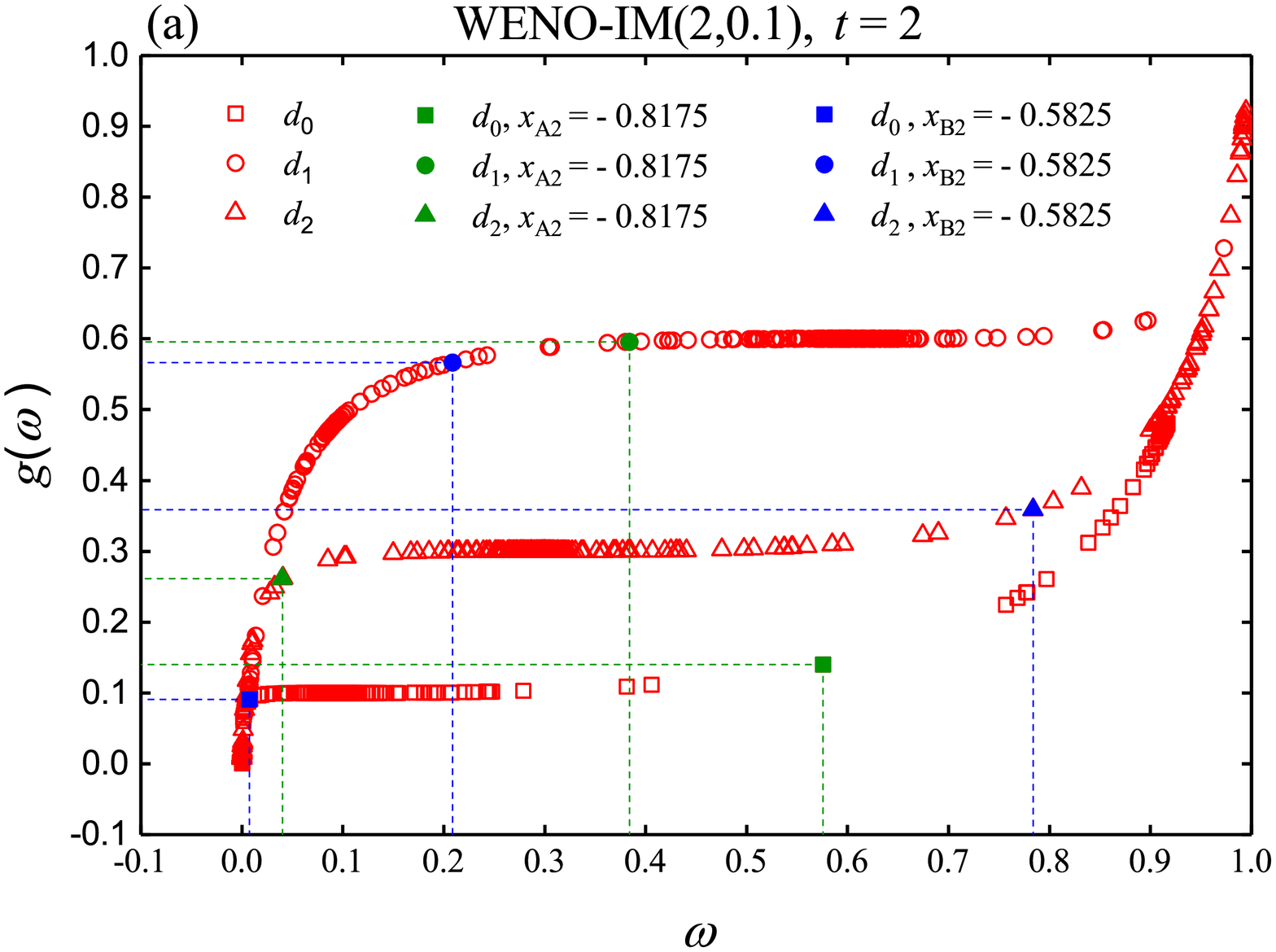}
\includegraphics[height=0.35\textwidth]
{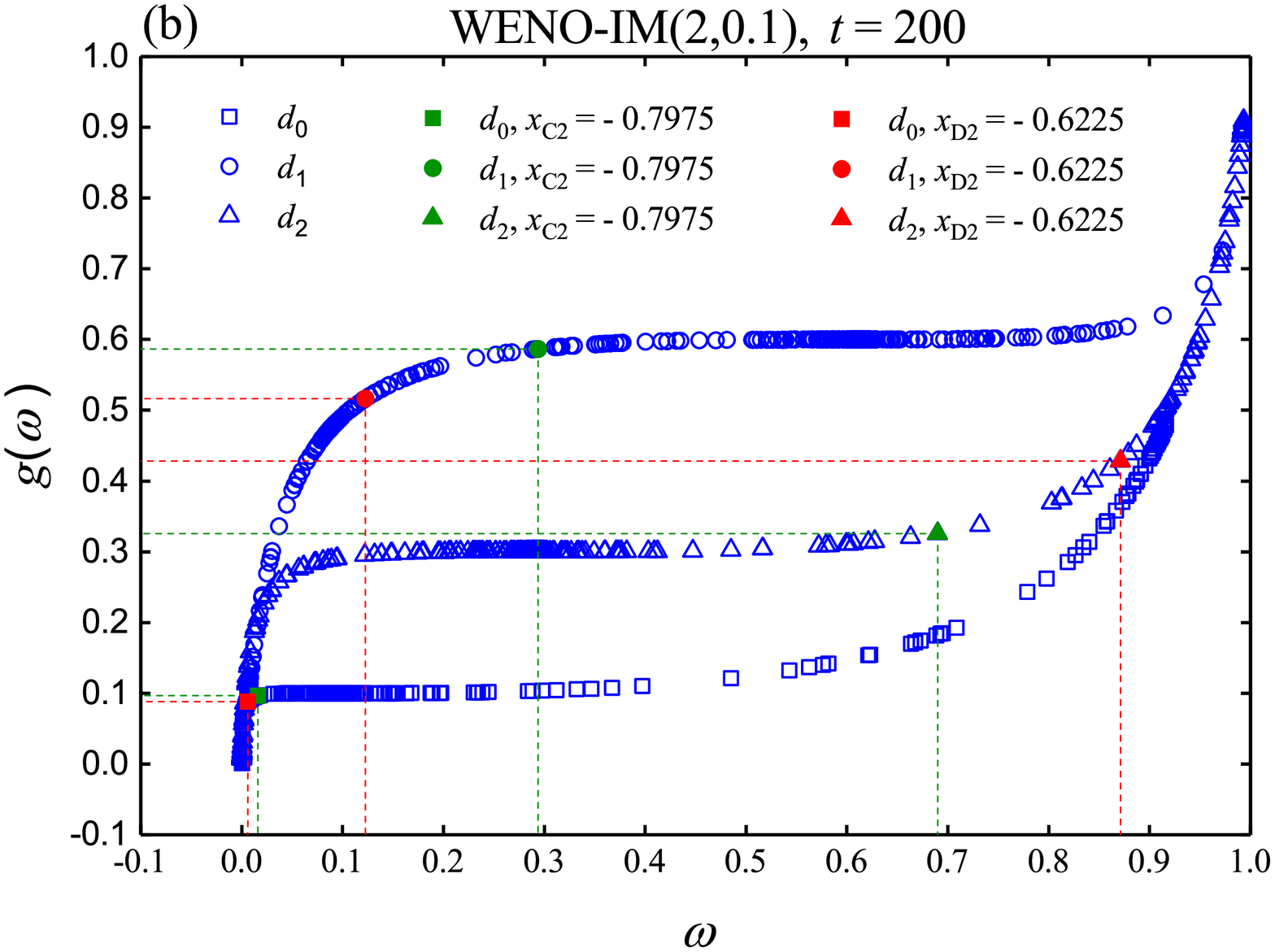}
\caption{The real-time mapping relationship $\big(g^{\mathrm{IM}(2,
0.1)}\big)(\omega)\sim\omega$ of the SLP. A uniform mesh size of
$N = 400$ and two output times $t=2$ (left) and $t=200$ (right) are 
used.}
\label{fig:RT-gOmega:deepAna:IM}
\end{figure}

\begin{figure}[ht]
\centering
\includegraphics[height=0.35\textwidth]
{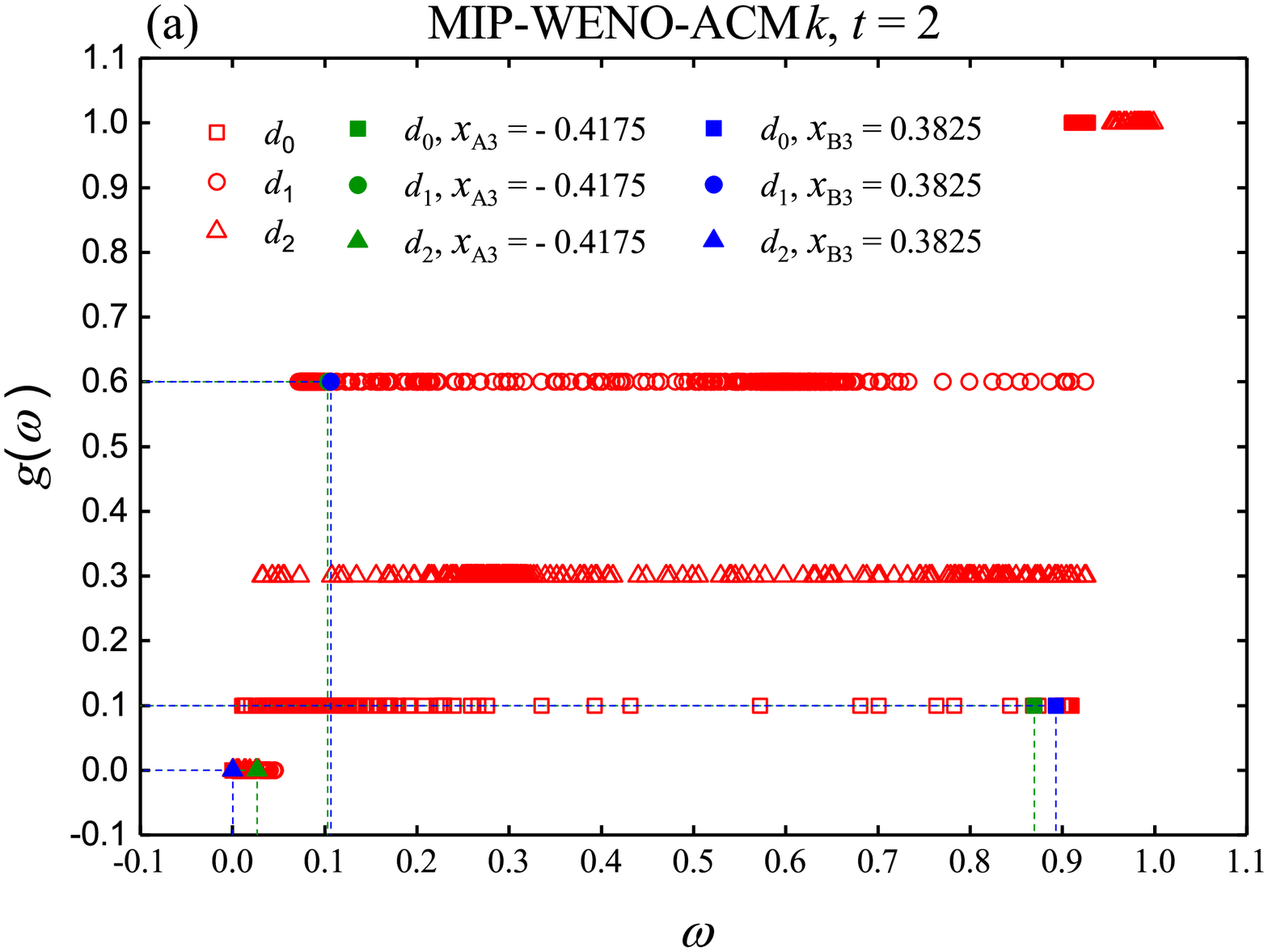}
\includegraphics[height=0.35\textwidth]
{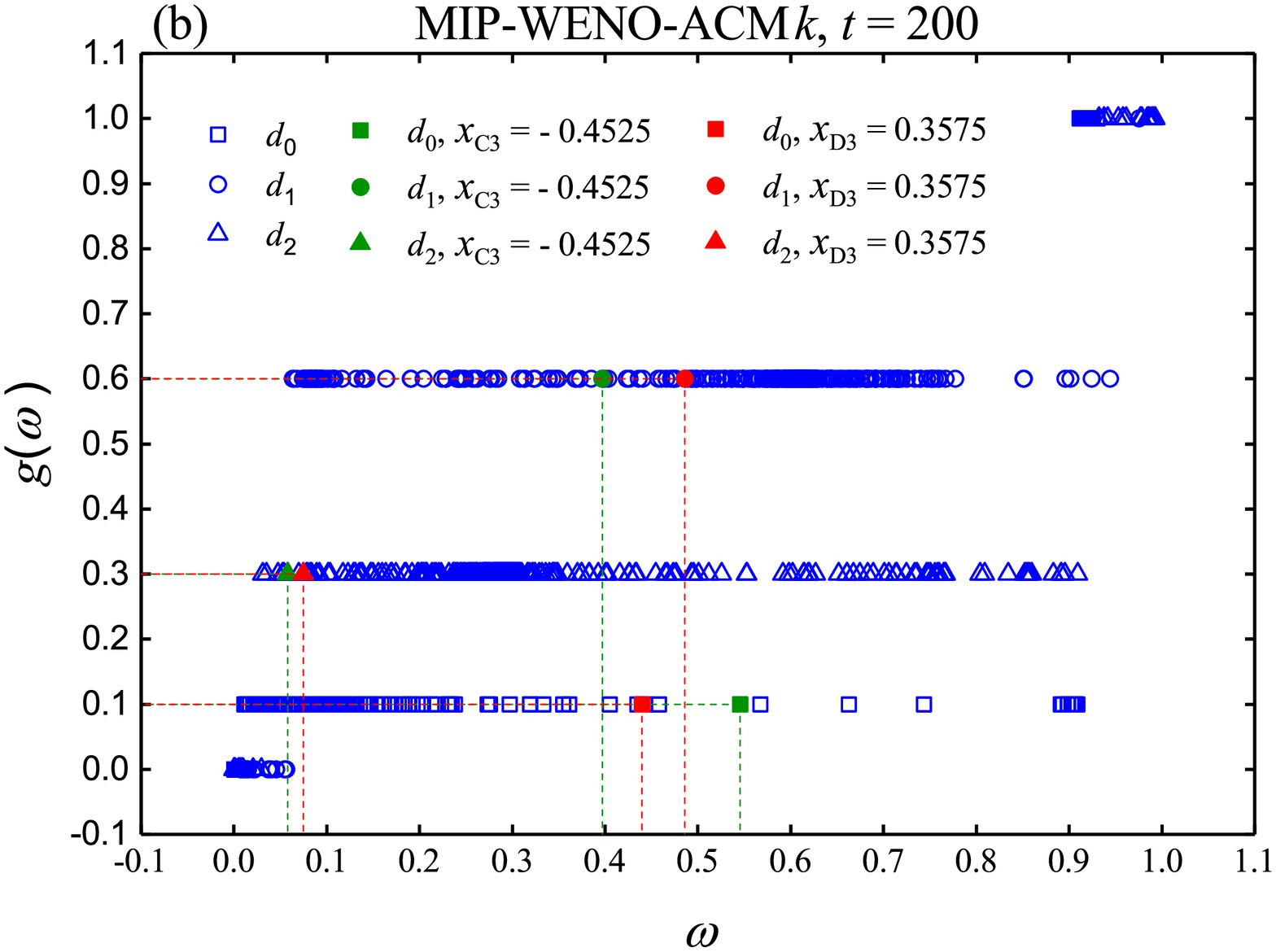}
\caption{The real-time mapping relationship $\big(g^{\mathrm{ACM}k}
\big)(\omega)\sim\omega$ of the SLP. A uniform mesh size of $N = 400$
and two output times $t=2$ (left) and $t=200$ (right) are used.}
\label{fig:RT-gOmega:deepAna:ACMk}
\end{figure}

\begin{figure}[ht]
\centering
\includegraphics[height=0.34\textwidth]
{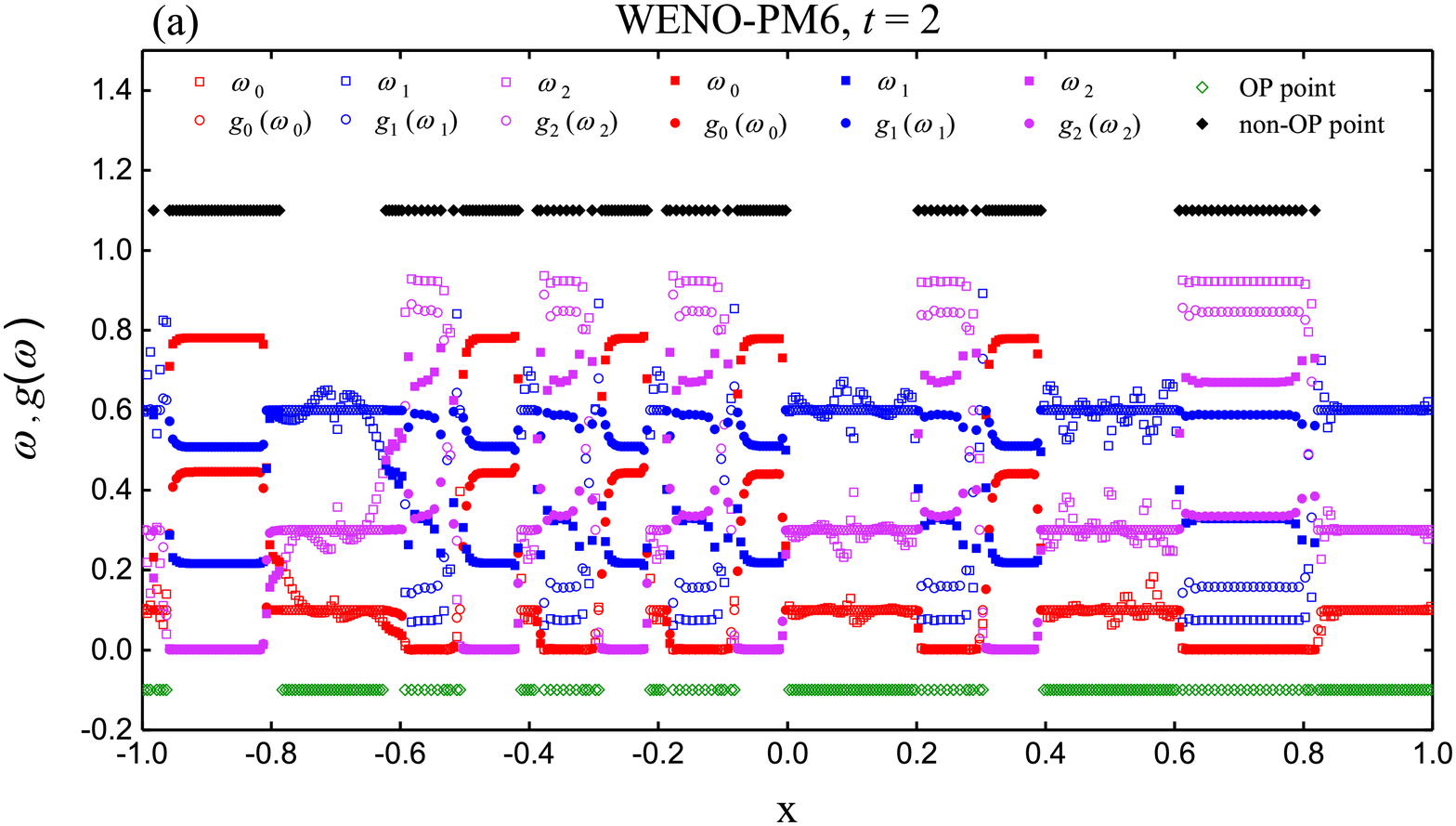}
\includegraphics[height=0.34\textwidth]
{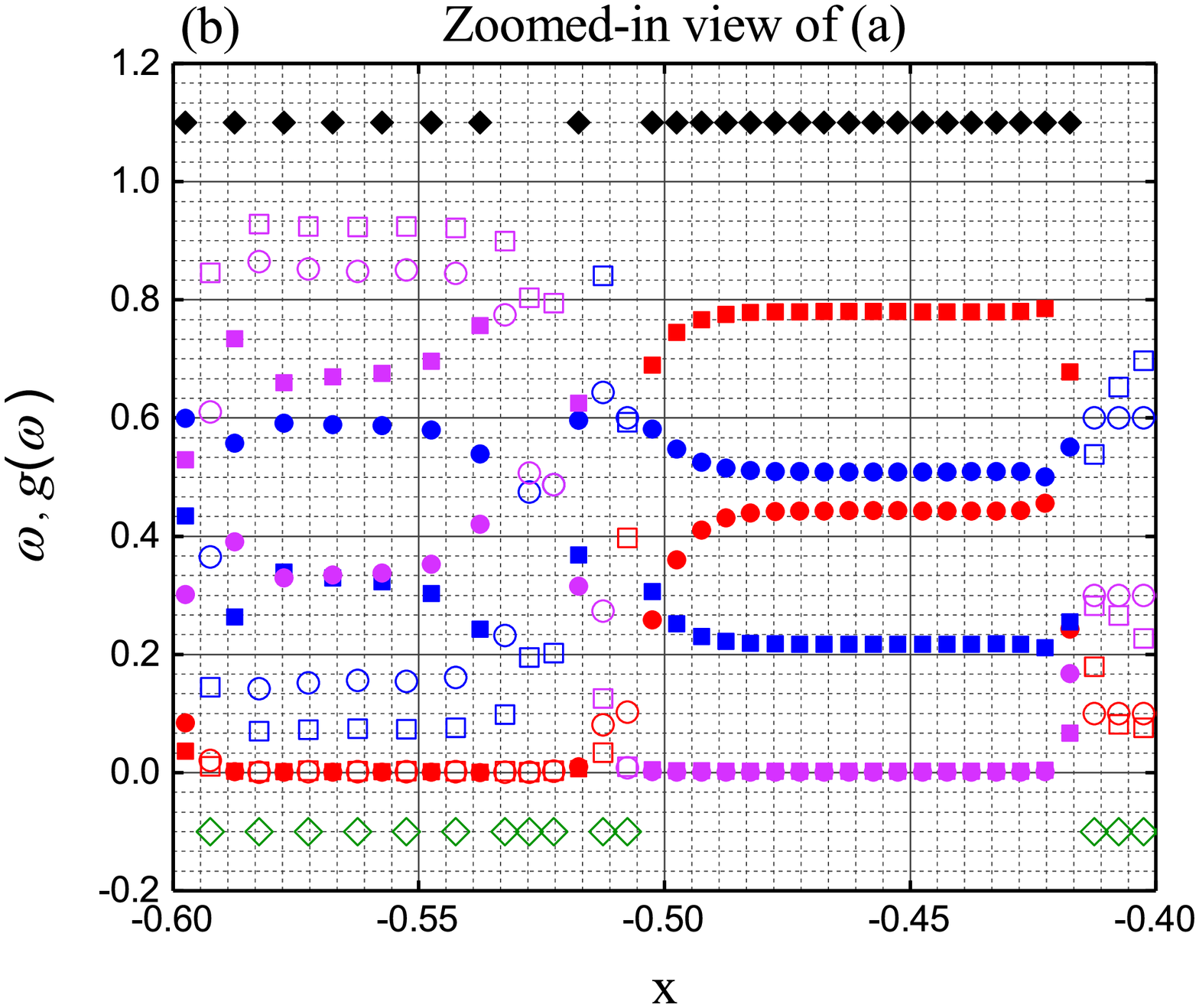}
\caption{The \textit{non-OP points} in the numerical solutions of 
the WENO-PM6 scheme. A uniform mesh size of $N = 400$ is used and 
the output time is $t=2$.}
\label{fig:RT-x-Omega:deepAna:PM6-t2}
\end{figure}

\begin{figure}[ht]
\centering
\includegraphics[height=0.34\textwidth]
{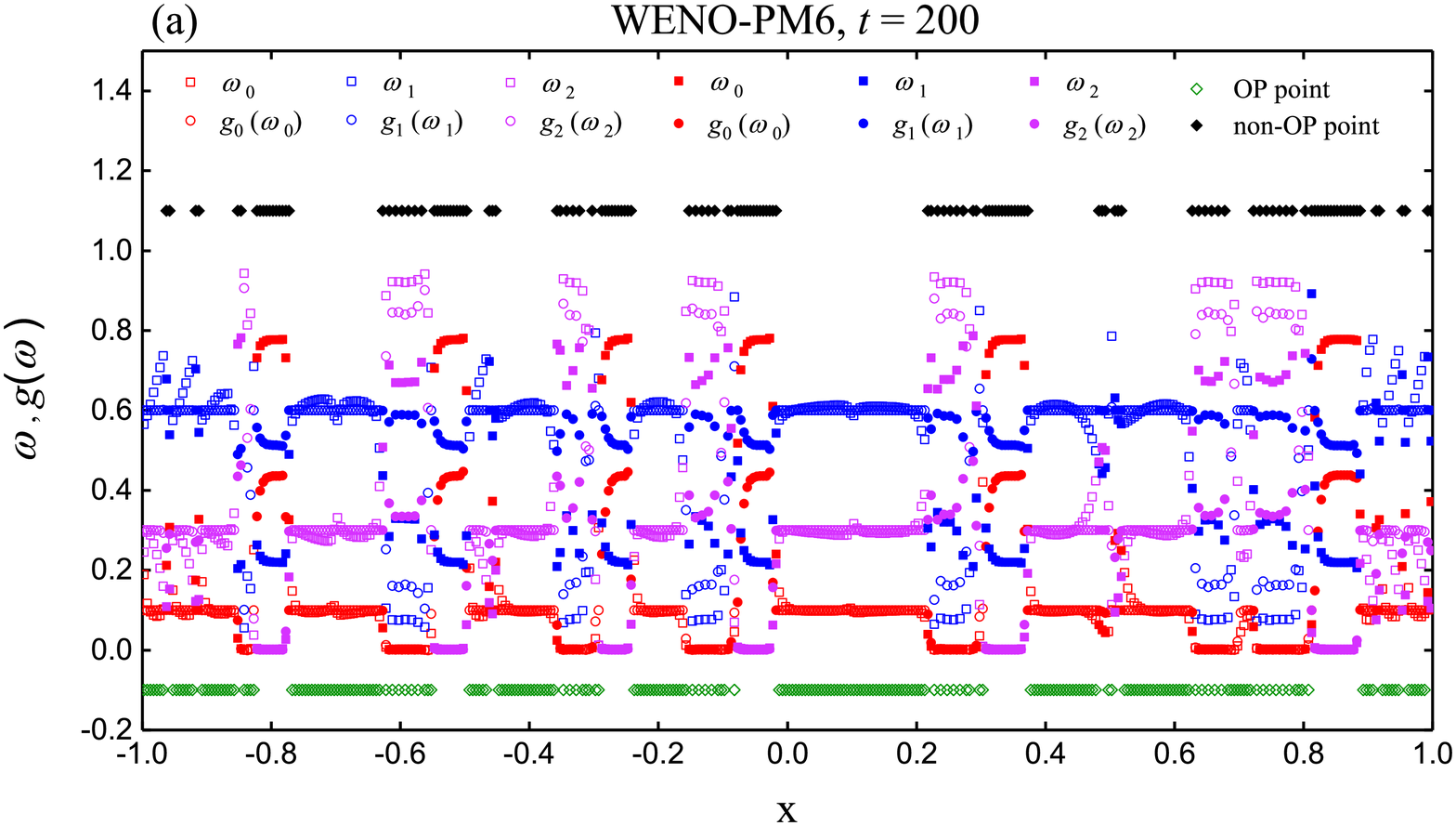}
\includegraphics[height=0.34\textwidth]
{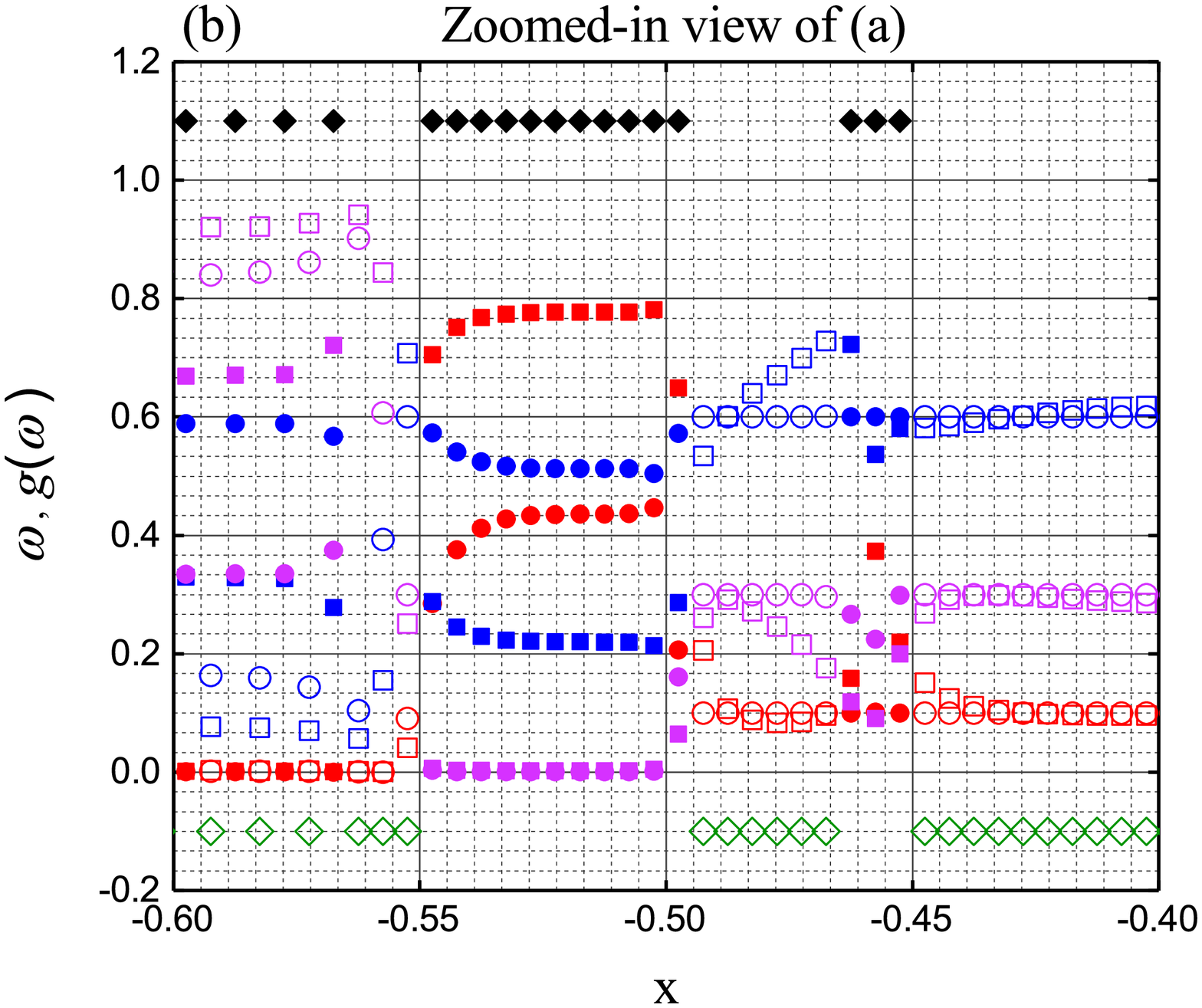}
\caption{The \textit{non-OP points} in the numerical solutions of 
the WENO-PM6 scheme. A uniform mesh size of $N = 400$ is used and 
the output time is $t=200$.}
\label{fig:RT-x-Omega:deepAna:PM6-t200}
\end{figure}

After extensive numerical experiments, we have discovered that, for 
almost all previously published mapped WENO schemes at least as far 
as we know, the non-OP mapping process will definitely occur when 
they are used for solving the problems with discontinuities. To 
demonstrate this, we still take the SLP as an example. The WENO-PM$6$
\cite{WENO-PM}, WENO-IM($2, 0.1$) \cite{WENO-IM} schemes and the 
MIP-WENO-ACM$k$ scheme with parameters $k_{s} = 0, \mathrm{CFS}_{s}=
\frac{d_{s}}{10}$ are used. A uniform mesh size of $N = 400$ and two 
output times $t=2$ (short) and $t=200$ (long) are taken in all 
calculations. In Fig. \ref{fig:RT-gOmega:deepAna:PM6} to Fig. 
\ref{fig:RT-gOmega:deepAna:ACMk}, we present the real-time mapping 
relationship $\big(g^{\mathrm{X}}\big)(\omega)\sim\omega$ of the 
considered mapped WENO schemes, where two of the 
\textit{non-OP points} are selected and highlighted in solid symbols 
for demonstration. In Table \ref{table:deepAna:2}, we present the 
computed values of the nonlinear weights, both before and after the 
mapping process, associated with these highlighted 
\textit{non-OP points}. We also give the order (in brackets in 
descending manner) of these nonlinear weights. It is evident that 
the order of the nonlinear weights has been changed when the mapping 
process is implemented at each \textit{non-OP point}. In addition, 
as shown in Fig. \ref{fig:RT-x-Omega:deepAna:PM6-t2} and Fig. 
\ref{fig:RT-x-Omega:deepAna:PM6-t200}, we see that there are many 
\textit{non-OP points} in the numerical solutions of the WENO-PM6 
scheme for both short and long output times. Similarly, we also find 
many \textit{non-OP points} in the results of the WENO-IM($2, 0.1$) 
and MIP-WENO-ACM$k$ schemes, while we do not present them here just 
for brevity.  

\begin{table}[ht]
\footnotesize
\centering
\caption{The mapping results of the SLP on highlighted 
\textit{non-OP points}, computed by WENO-PM6, WENO-IM($2,0.1$) and 
MIP-WENO-ACM$k$, with a uniform mesh size of $N=400$ and output 
times $t=2, 200$.}
\label{table:deepAna:2}
\begin{tabular*}{\hsize}
{@{}@{\extracolsep{\fill}}llllllllll@{}}
\hline
\space&\space&\space&\space &\multicolumn{3}{l}{Before the mapping 
(the order)} &\multicolumn{3}{l}{After the mapping (the order)}\\
\cline{5-7}  \cline{8-10}
Schemes, X & Time, $t$ & Point & Position, $x$  &$\omega_{0}$ & 
$\omega_{1}$ & $\omega_{2}$ &$g_{0}(\omega_{0})$ &$g_{1}(\omega_{1})$
&$g_{2}(\omega_{2})$ \\
\hline
WENO-PM6 & 2   & A1 & -0.4175 & 0.67828(1) & 0.25528(2) & 0.06644(3) 
                             & 0.24252(2) & 0.55068(1) & 0.16737(3)\\
{}       & 2   & B1 & -0.5175 & 0.00678(3) & 0.36849(2) & 0.62473(1) 
                             & 0.00980(3) & 0.59595(1) & 0.31538(2)\\
{}       & 200 & C1 & -0.4575 & 0.37291(2) & 0.53663(1) & 0.09046(3) 
                             & 0.10125(3) & 0.60000(1) & 0.22432(2)\\
{}       & 200 & D1 & -0.4975 & 0.64949(1) & 0.28636(2) & 0.06415(3) 
                             & 0.20605(2) & 0.57222(1) & 0.16098(3)\\
\hline
WENO-IM($2,0.1$) 
         & 2   & A2 & -0.8175 & 0.57568(1) & 0.38416(2) & 0.04016(3) 
                             & 0.14033(3) & 0.59583(1) & 0.26127(2)\\
{}       & 2   & B2 & -0.5825 & 0.00768(3) & 0.20854(2) & 0.78378(1) 
                             & 0.09071(3) & 0.56674(1) & 0.35871(2)\\
{}       & 200 & C2 & -0.7975 & 0.01609(3) & 0.29405(2) & 0.68986(1) 
                             & 0.09643(3) & 0.58680(1) & 0.32586(2)\\
{}       & 200 & D2 & -0.6225 & 0.00622(3) & 0.12285(2) & 0.87093(1) 
                             & 0.08832(3) & 0.51677(1) & 0.42834(2)\\
\hline
MIP-WENO-ACM$k$ 
         & 2   & A3 & -0.4175 & 0.86977(1) & 0.10334(2) & 0.02689(3) 
                             & 0.10000(2) & 0.60000(1) & 0.00000(3)\\
{}       & 2   & B3 &  0.3825 & 0.89299(1) & 0.10671(2) & 0.00030(3) 
                             & 0.10000(2) & 0.60000(1) & 0.00000(3)\\
{}       & 200 & C3 & -0.4525 & 0.54547(1) & 0.39684(2) & 0.05769(3) 
                             & 0.10000(3) & 0.60000(1) & 0.30000(2)\\
{}       & 200 & D3 &  0.3575 & 0.43952(2) & 0.48568(1) & 0.07480(3) 
                             & 0.10000(3) & 0.60000(1) & 0.30000(2)\\
\hline
\end{tabular*}
\end{table}

\subsubsection{Effects of the non-OP mapping process 
on the numerical solutions}
Without loss of generality, we assume that the weights $\omega_{s}^{
\mathrm{JS}}$, which would be substituted into some mapping function,
satisfy $\omega_{0}^{\mathrm{JS}} > \omega_{1}^{\mathrm{JS}} > 
\omega_{2}^{\mathrm{JS}}$, the mapped weights 
$\omega_{s}^{\mathrm{non-OP}}$ computed by some set of mapping 
functions which is \textit{non-OP} satisfy 
$\omega_{1}^{\mathrm{non-OP}} > \omega_{2}^{\mathrm{non-OP}} > 
\omega_{0}^{\mathrm{non-OP}}$, and the mapped weights
$\omega_{s}^{\mathrm{OP}}$ computed by some set of mapping functions 
which is \textit{OP} satisfy $\omega_{0}^{\mathrm{OP}} > 
\omega_{1}^{\mathrm{OP}} > \omega_{2}^{\mathrm{OP}}$.

According to Eq.(\ref{eq:approx_WENO}), by adding and subtracting 
$\sum_{s = 0}^{2}d_{s}^{\mathrm{non-OP}}u_{j + 1/2}^{s}$ to 
$u_{j + 1/2}^{\mathrm{non-OP}} = \sum_{s = 0}^{2}\omega_{s}^{
\mathrm{non-OP}}u_{j + 1/2}^{s}$, we obtain
\begin{equation}
u_{j+1/2}^{\mathrm{non-OP}}=\sum\limits_{s=0}^{2}d_{s}^{
\mathrm{non-OP}}u_{j+1/2}^{s} + \sum\limits_{s=0}^{2}(\omega_{s}^{
\mathrm{non-OP}} - d_{s}^{\mathrm{non-OP}})u_{j+1/2}^{s}.
\label{eq:approxWENO:non-OP}
\end{equation}
In smooth regions, the first term in the right-hand side of 
Eq.(\ref{eq:approxWENO:non-OP}) satisfies 
\cite{WENO-M,WENO-Z,WENO-IM}
\begin{equation}
\sum\limits_{s=0}^{2}d_{s}^{\mathrm{non-OP}}u_{j+1/2}^{s}=
u(x_{j+1/2}) + O(\Delta x^{5}).
\label{eq:approxWENO:non-OP:ideal}
\end{equation}
Similarly, we have
\begin{equation}
u_{j+1/2}^{\mathrm{OP}}=\sum\limits_{s=0}^{2}d_{s}^{\mathrm{OP}}
u_{j+1/2}^{s} + \sum\limits_{s=0}^{2}(\omega_{s}^{\mathrm{OP}} - 
d_{s}^{\mathrm{OP}})u_{j+1/2}^{s}, 
\label{eq:approxWENO:OP}
\end{equation}
and
\begin{equation}
\sum\limits_{s=0}^{2}d_{s}^{\mathrm{OP}}u_{j+1/2}^{s}=u(x_{j+1/2})+
O(\Delta x^{5}).
\label{eq:approxWENO:OP:ideal}
\end{equation}
So the second term in the right-hand side of 
Eq.(\ref{eq:approxWENO:non-OP}) or Eq.(\ref{eq:approxWENO:OP}) must 
be at least an $O(\Delta x^{6})$ quantity to ensure the convergence 
rate to be approximated at $5$th-order, and this is the key point 
that the mapped WENO methods have focused on.

However, in the parts of solutions with discontinuities, 
Eq.(\ref{eq:approxWENO:non-OP:ideal}) and 
Eq.(\ref{eq:approxWENO:OP:ideal}) usually do not hold. Furthermore, 
it is easy to know that the non-OP mapping process will amplify the 
effect from the relatively non-smooth stencils and decrease the 
effect from the relatively smooth stencils, so that we can probably 
get the following inequality
\begin{equation}
\Bigg\lvert \sum\limits_{s = 0}^{2}d_{s}^{\mathrm{non-OP}}
u_{j + 1/2}^{s} - u(x_{j + 1/2}) \Bigg\rvert > \Bigg\lvert \sum
\limits_{s = 0}^{2}d_{s}^{\mathrm{OP}}u_{j + 1/2}^{s} - 
u(x_{j + 1/2}) \Bigg\rvert.
\label{eq:theoreticalAna4non-OP:firstTerm}
\end{equation}

Now, we analyze the effect of the second term in the right-hand side 
of Eq.(\ref{eq:approxWENO:non-OP}) and Eq.(\ref{eq:approxWENO:OP}). 
Suppose that $\big( g^{\mathrm{non-OP}}\big)_{s}(d_{s}^{
\mathrm{non-OP}})=d_{s}^{\mathrm{non-OP}}, \big( g^{\mathrm{non-OP}} 
\big)_{s}^{'}(d_{s}^{\mathrm{non-OP}}) = \cdots = \big( 
g^{\mathrm{non-OP}}\big)_{s}^{(n - 1)}(d_{s}^{\mathrm{non-OP}}) = 0, 
\big(g^{\mathrm{non-OP}} \big)_{s}^{(n)}(d_{s}^{\mathrm{non-OP}})\neq
0$. Then, evaluation at $\omega_{s}^{\mathrm{JS}}$ of the Taylor 
series approximations of $\big(g^{\mathrm{non-OP}}\big)_{s}
(\omega)$ about $d_{s}^{\mathrm{non-OP}}$ yields
\begin{equation}
\begin{array}{ll}
\begin{aligned}
\alpha_{s}^{\mathrm{non-OP}} &= \big( g^{\mathrm{non-OP}} \big)_{s}
(d_{s}^{\mathrm{non-OP}}) + \big( g^{\mathrm{non-OP}} \big)_{s}^{'}
(d_{s}^{\mathrm{non-OP}})\big(\omega_{s}^{\mathrm{JS}} - d_{s}^{
\mathrm{non-OP}}\big)+\cdots\\
\space &+\dfrac{\big( g^{\mathrm{non-OP}} 
\big)_{s}^{(n)}(d_{s}^{\mathrm{non-OP}})}{n!}\big(\omega_{s}^{
\mathrm{JS}} - d_{s}^{\mathrm{non-OP}} \big)^{n} + O\bigg(
\big(\omega_{s}^{\mathrm{JS}} - d_{s}^{\mathrm{non-OP}} \big)^{n+1}
\bigg)\\
\space &=d_{s}^{\mathrm{non-OP}} + \dfrac{\big( g^{
\mathrm{non-OP}} \big)_{s}^{(n)}(d_{s}^{\mathrm{non-OP}})}{n!}\big(
\omega_{s}^{\mathrm{JS}} - d_{s}^{\mathrm{non-OP}} \big)^{n} + 
O\bigg(\big(\omega_{s}^{\mathrm{JS}} - d_{s}^{\mathrm{non-OP}} 
\big)^{n+1}\bigg).
\end{aligned}
\end{array}
\label{eq:proof:alpha_non-OP:01}
\end{equation}
Similarly, when $\big(g^{\mathrm{OP}}\big)_{s}(d_{s}^{\mathrm{OP}})
=d_{s}^{\mathrm{OP}}, \big( g^{\mathrm{OP}} \big)_{s}^{'}(d_{s}^{
\mathrm{OP}}) = \cdots = \big( g^{\mathrm{OP}} \big)_{s}^{(n - 1)}
(d_{s}^{\mathrm{OP}}) = 0, \big( g^{\mathrm{OP}} \big)_{s}^{(n)}(
d_{s}^{\mathrm{OP}}) \neq 0$, we have
\begin{equation}
\begin{array}{ll}
\begin{aligned}
\alpha_{s}^{\mathrm{OP}} &= \big( g^{\mathrm{OP}} \big)_{s}(d_{s}^{
\mathrm{OP}}) + \big( g^{\mathrm{OP}} \big)_{s}^{'}(d_{s}^{
\mathrm{OP}})\big(\omega_{s}^{\mathrm{JS}} - d_{s}^{\mathrm{OP}}
\big)+\cdots+\dfrac{\big( g^{\mathrm{OP}} \big)_{s}^{(n)}(d_{s}^{
\mathrm{OP}})}{n!}\big(\omega_{s}^{\mathrm{JS}}-d_{s}^{\mathrm{OP}}
\big)^{n} + O\bigg(\big(\omega_{s}^{\mathrm{JS}}-d_{s}^{\mathrm{OP}}
\big)^{n+1}\bigg)\\
\space &=d_{s}^{\mathrm{OP}} + \dfrac{\big( g^{
\mathrm{OP}} \big)_{s}^{(n)}(d_{s}^{\mathrm{OP}})}{n!}\big(
\omega_{s}^{\mathrm{JS}} - d_{s}^{\mathrm{OP}} \big)^{n} + O\bigg(
\big(\omega_{s}^{\mathrm{JS}}-d_{s}^{\mathrm{OP}} \big)^{n+1}\bigg).
\end{aligned}
\end{array}
\label{eq:proof:alpha_OP:01}
\end{equation}
Then, from Eq.(\ref{eq:proof:alpha_non-OP:01}) and Eq.(
\ref{eq:proof:alpha_OP:01}), we obtain
\begin{equation}
\dfrac{\alpha_{s}^{\mathrm{non-OP}} - d_{s}^{\mathrm{non-OP}}}{
\alpha_{s}^{\mathrm{OP}} - d_{s}^{\mathrm{OP}}} \approx \dfrac{\big( 
g^{\mathrm{non-OP}} \big)_{s}^{(n)}(d_{s}^{\mathrm{non-OP}})}{\big( 
g^{\mathrm{OP}} \big)_{s}^{(n)}(d_{s}^{\mathrm{OP}})} \times \Bigg(
\dfrac{\omega_{s}^{\mathrm{JS}} - d_{s}^{\mathrm{non-OP}}}{
\omega_{s}^{\mathrm{JS}} - d_{s}^{\mathrm{OP}}} \Bigg)^{n}.
\label{eq:proof:alpha_non-OP-OP:01}
\end{equation}

One may probably get $\dfrac{\omega_{s}^{\mathrm{JS}} - 
d_{s}^{\mathrm{non-OP}}}{\omega_{s}^{\mathrm{JS}} - d_{s}^{
\mathrm{OP}}}\gg 1$ at a \textit{non-OP point} where a non-OP 
mapping process occurs. Then, as $n$ is a possibly large positive 
integer, Eq.(\ref{eq:proof:alpha_non-OP-OP:01}) yields
\begin{equation}
\alpha_{s}^{\mathrm{non-OP}} - d_{s}^{\mathrm{non-OP}} \gg 
\alpha_{s}^{\mathrm{OP}} - d_{s}^{\mathrm{OP}}. 
\label{eq:proof:alpha_non-OP-alpha_OP}
\end{equation}
Thus, according to Eq.(46) in \cite{WENO-M}, it is trivial to know 
that
\begin{equation}
\omega_{s}^{\mathrm{non-OP}} - d_{s}^{\mathrm{non-OP}} \gg 
\omega_{s}^{\mathrm{OP}} - d_{s}^{\mathrm{OP}}. 
\label{eq:omega_non-OP-omega_OP}
\end{equation}

Now, from Eq.(\ref{eq:approxWENO:non-OP})(\ref{eq:approxWENO:OP})(\ref{eq:theoreticalAna4non-OP:firstTerm})(\ref{eq:omega_non-OP-omega_OP}) and considering the accumulation of 
the errors when $t$ gets larger, we can conclude that the non-OP 
mapping process might probably be the essential cause of the 
spurious oscillation generation and potential loss of 
accuracy when the mapped WENO schemes are used to simulate the 
problems with discontinuities for long output times. 

\begin{figure}[ht]
\centering
  \includegraphics[height=0.38\textwidth]
  {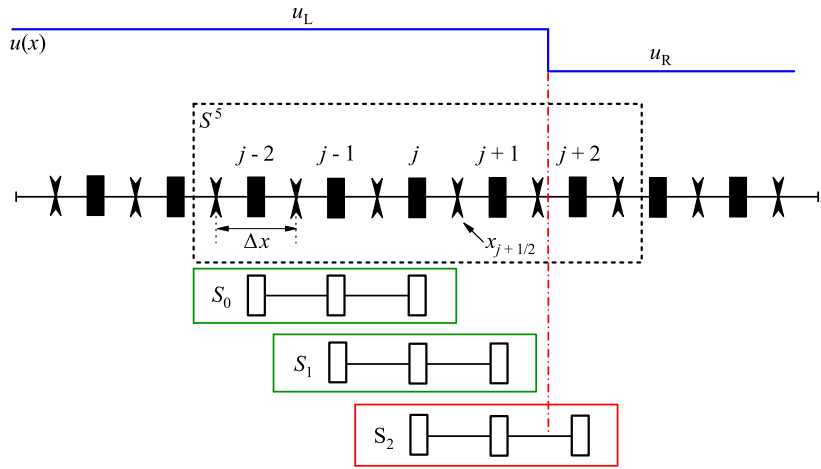}
\caption{Schematic of fifth-order WENO stencils with an isolated 
discontinuity on one of the substencils.}
\label{fig:illustration}
\end{figure}

To illustrate this, and for brevity in the discussion but without 
loss of generality, we assume that there is a global stencil $S^{5}$ 
which is divided into $3$ substencils $S_{0}, S_{1}, S_{2}$, and 
there is an isolated discontinuity on $S_{2}$, as shown in Fig. 
\ref{fig:illustration}. We suppose $u_{\mathrm{L}}=1, u_{\mathrm{R}}=
-1$ and $u_{j+1/2}^{0}=u_{j+1/2}^{1}=1, u_{j+1/2}^{2}=-1$. According 
to Eq.(\ref{eq:approx_WENO}), we calculate the approximation of 
$u(x_{j+1/2})$ on the stencil $S^{5}$, by applying the mapped 
weights of the \textit{non-OP points} C1, A2, C3 in Table 
\ref{table:deepAna:2} to the corresponding substencils in Fig.
\ref{fig:illustration}, respectively. For comparison, we also 
calculate the results by applying the weights of the WENO-JS scheme 
and the mapped weights of some \textit{OP points}. Here, for 
simplicity but without loss of generality, we directly use the same 
values of the mapped weights of the \textit{non-OP points} C1, A2, 
C3 in Table \ref{table:deepAna:2}, but change their order, to set 
the values of  the mapped weights of the \textit{OP points}. We 
present the computing conditions and results in detail in Table 
\ref{table:illustration}. 
From Table \ref{table:illustration}, we find that the errors 
computed by using the mapped weights of the \textit{non-OP points} 
are much larger than the solutions computed by using the mapped 
weights of the \textit{OP points}. Although we only present a 
pseudo-test example here, it is highly conducive to describe and 
understand the way that the non-OP mapping process causes the 
potential loss of accuracy. In practice tests with short output 
times, this phenomenon, as well as the spurious oscillation, may not 
be easy to be noticed. However, the errors are accumulated and will 
be demonstrated when the output time gets larger, and then the 
spurious oscillation and potential loss of accuracy can be observed. 
We will show these through numerical experiments in Section
\ref{NumericalExperiments}. 

\begin{table}[ht]
\footnotesize
\centering
\caption{The computing conditions and results for comparison.}
\label{table:illustration}
\renewcommand\arraystretch{1.3}
\begin{tabular*}{\hsize}
{@{}@{\extracolsep{\fill}}llll@{}}
\hline
\space    &WENO-PM6, Point C1 &WENO-IM(2,0.1), Point A2
          &MIP-WENO-ACM$k$, Point C3\\
\hline
($\omega_{0}^{\mathrm{JS}},
  \omega_{1}^{\mathrm{JS}},
  \omega_{2}^{\mathrm{JS}}$) 
& (0.37291, 0.53663, 0.09046) & (0.57568, 0.38416, 0.04016)
& (0.54547, 0.39684, 0.05769) \\
($\omega_{0}^{\mathrm{non-OP}},
  \omega_{1}^{\mathrm{non-OP}},
  \omega_{2}^{\mathrm{non-OP}}$) 
& (0.10939, 0.64825, 0.24236) & (0.14069, 0.59737, 0.26194)
& (0.10000, 0.60000, 0.30000) \\
($\omega_{0}^{\mathrm{OP}},
  \omega_{1}^{\mathrm{OP}},
  \omega_{2}^{\mathrm{OP}}$) 
& (0.24236, 0.64825, 0.10939) & (0.59737, 0.26194, 0.14069)
& (0.60000, 0.30000, 0.10000) \\
($u_{j+1/2}^{\mathrm{JS}},
  u_{j+1/2}^{\mathrm{non-OP}},
  u_{j+1/2}^{\mathrm{OP}}$) 
& (0.81908, 0.51528, 0.78122) & (0.91968, 0.47611, 0.71862)
& (0.88462, 0.40000, 0.80000) \\
$u_{j+1/2}^{\mathrm{JS}}-u(x_{j+1/2})$
& 0.18092(18.09\%) & 0.08031(8.03\%)  & 0.11538(11.54\%) \\
$u_{j+1/2}^{\mathrm{non-OP}}-u(x_{j+1/2})$
& 0.48472(48.47\%) & 0.52389(52.39\%) & 0.60000(60.00\%) \\
$u_{j+1/2}^{\mathrm{OP}}-u(x_{j+1/2})$
& 0.21878(21.88\%) & 0.28138(28.14\%) & 0.20000(20.00\%) \\
\hline
\end{tabular*}
\end{table}


\section{Design and properties of the order-preserving mapping 
functions}
\label{MOP-WENO-ACMk}
\subsection{Design of the new mapping function}
To design a new mapped WENO scheme that can prevent the 
\textit{non-OP mapping process}, we devise a set of mapping 
functions that is OP.

Let $\mathcal{D} = \{d_{0},d_{1},\cdots,d_{r-1}\}$ be an array of 
all the ideal weights of the $(2r-1)$th-order WENO schemes. We build 
a new array by sorting the elements of $\mathcal{D}$ in ascending 
order, that is, $\mathcal{\widetilde{D}} = \{\widetilde{d}_{0},
\widetilde{d}_{1},\cdots,\widetilde{d}_{r-1}\}$. In other words, the 
arrays $\mathcal{D}$ and $\mathcal{\widetilde{D}}$ have the same 
elements with different arrangements, and the elements of 
$\mathcal{\widetilde{D}}$ satisfy
\begin{equation}
\widetilde{d}_{0} < \widetilde{d}_{1} < \cdots < \widetilde{d}_{r-1}.
\label{eq:tilde_ds}
\end{equation}

The following notations are introduced to simplify the expressions
\begin{equation}
\hat{\Omega}_{1}=[0,\mathrm{CFS}_{0}),\hat{\Omega}_{2}=\Big[
\mathrm{CFS}_{0},
\dfrac{\widetilde{d}_{0} + \widetilde{d}_{1}}{2}\Big), \cdots,
\hat{\Omega}_{r+1}=\Big[\dfrac{\widetilde{d}_{r-2}+
\widetilde{d}_{r-1}}{2}
,\mathrm{CFS}_{1}\Big],\hat{\Omega}_{r+2}=\Big(
\mathrm{CFS}_{1},1\Big],
\label{eq:Omega_i}
\end{equation}
where $0<\mathrm{CFS}_{0} \leq \widetilde{d}_{0}$ and 
$\widetilde{d}_{r - 1} \leq \mathrm{CFS}_{1} \leq 1$. 
It is easy to verify that: (1) $\hat{\Omega} = \hat{\Omega}_{1}\cup
\hat{\Omega}_{2}\cup\cdots\cup\hat{\Omega}_{r+2}$; (2) if $i \neq j$ 
, then $\hat{\Omega}_{i}\cap\hat{\Omega}_{j} = \varnothing$, $\forall
i,j=1,2,\cdots,r+2$. 

Now, we give a new mapping function as follows
\begin{equation}
\big( g^{\mathrm{MOP-ACM}k} \big)_{s}(\omega) = \left\{
\begin{array}{lll}
k_{0} \omega,& {} & \omega \in \hat{\Omega}_{1}, \\
\widetilde{d}_{0},   & {} & \omega \in \hat{\Omega}_{2}, \\
\widetilde{d}_{1},   & {} & \omega \in \hat{\Omega}_{3}, \\
{}       & \vdots & {}                 \\
\widetilde{d}_{r-1}, & {} & \omega \in \hat{\Omega}_{r+1}, \\
1 - k_{1}(1 - \omega), & {} & \omega \in \hat{\Omega}_{r+2},
\end{array}
\right.
\label{eq:mapping_MOP-ACMk}
\end{equation}
where
\begin{equation}
k_{0} \in \Bigg[0, \dfrac{\widetilde{d}_{0}}{\mathrm{CFS}_{0}}\Bigg],
k_{1} \in \Bigg[0, \dfrac{1 - \widetilde{d}_{r-1}}{1 - 
\mathrm{CFS}_{1}}\Bigg].
\label{eq:rangeOf_k}
\end{equation}

Actually, we can verify that $\big( g^{\mathrm{MOP-ACM}k} \big)_{s}
(\omega)$ is independent of the parameter $s$, that is, $\big( 
g^{\mathrm{MOP-ACM}k} \big)_{0}(\omega) = \cdots = \big( 
g^{\mathrm{MOP-ACM}k}\big)_{r-1}(\omega)$, which is significantly 
different from the previously published mapping functions. Thus, for 
simplicity, we will drop the subscript $s$ of
$\big(g^{\mathrm{MOP-ACM}k}\big)_{s}(\omega)$ in the following.

\subsection{Properties of the new mapping function}
\begin{lemma}
The mapping function $g^{\mathrm{MOP-ACM}k}(\omega)$ defined by 
Eq.(\ref{eq:mapping_MOP-ACMk}) is \textit{order-preserving (OP)}.
\end{lemma}
\textbf{Proof.}
According to 
Eq.(\ref{eq:tilde_ds})(\ref{eq:Omega_i})(\ref{eq:mapping_MOP-ACMk}), 
we can easily verify that $g^{\mathrm{MOP-ACM}k}(\omega)$ satisfies 
Definition \ref{def:OPM}.
$\hfill\square$ \\

As mentioned earlier, $g^{\mathrm{MOP-ACM}k}(\omega)$ is a monotone 
increasing piecewise mapping function. It satisfies the following 
properties.

\begin{lemma}
Let $\overline{\Omega}_{i}=\{\omega \in \hat{\Omega}_{i}$ and
$\omega \neq \partial \hat{\Omega}_{i}\}$, then the mapping function
$g^{\mathrm{MOP-ACM}k}(\omega)$ defined by 
Eq.(\ref{eq:mapping_MOP-ACMk}) satisfies the following properties: \\

C1. for $\forall \omega \in \overline{\Omega}_{i}, i= 1,\cdots,r+2$,
$\big( g^{\mathrm{MOP-ACM}k} \big)'(\omega) \geq 0$;

C2. for $\forall \omega \in \Omega, 0 \leq g^{\mathrm{MOP-ACM}k}(
\omega) \leq 1$;

C3. for $\forall s \in \{ 0,1,\cdots,r-1\}, d_{s} \in \Omega_{s+2},
g^{\mathrm{MOP-ACM}k}(d_{s}) = d_{s},
\big(g^{\mathrm{MOP-ACM}k}\big)'(d_{s})=\big(g^{\mathrm{MOP-ACM}k}
\big)''(d_{s})=\cdots=0$;

C4. $g^{\mathrm{MOP-ACM}k}(0) = 0, g^{\mathrm{MOP-ACM}k}(1) = 1, 
\big( g^{\mathrm{MOP-ACM}k}\big)'(0) = k_{0}, 
\big( g^{\mathrm{MOP-ACM}k}\big)'(1) = k_{1}$.  
\label{lemma:mappingFunction:MOP-ACMk}
\end{lemma}

\begin{figure}[ht]
\centering
\includegraphics[height=0.38\textwidth]
{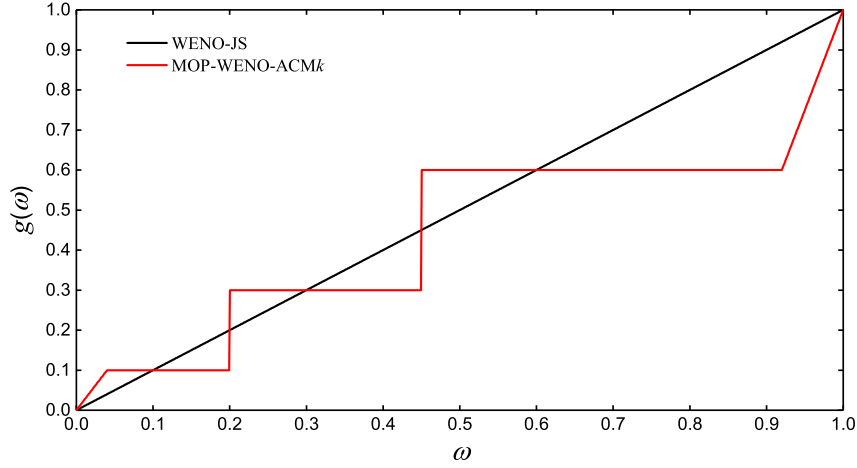}\quad
\caption{The mapping function of the MOP-WENO-ACM$k$ scheme with
$\mathrm{CFS}_{0} = 0.04, \mathrm{CFS}_{1} = 0.92, k_{0} = 2.5, k_{1}
=5$.}
\label{fig:gOmega:MOP-ACMk}
\end{figure}

As the proof of Lemma \ref{lemma:mappingFunction:MOP-ACMk} is very 
easy, we do not state them here and we can observe these properties 
intuitively in Fig. \ref{fig:gOmega:MOP-ACMk}.

Now, we give the new mapped WENO scheme with 
\textit{monotone increasing piecewise} and 
\textit{order-preserving mapping}, denoted as 
MOP-WENO-ACM$k$. The mapped weights are given by
\begin{equation}
\omega_{s}^{\mathrm{MOP-ACM}k} = \dfrac{\alpha _{s}^{\mathrm{MOP-ACM}
k}}{\sum_{l = 0}^{r-1} \alpha _{l}^{\mathrm{MOP-ACM}k}}, \alpha_{s}^{
\mathrm{MOP-ACM}k}=g^{\mathrm{MOP-ACM}k}(\omega^{\mathrm{JS}}_{s}).
\label{eq:nonlinearWeightsX}
\end{equation}

We present Theorem \ref{theorem:convergenceRates}, which will 
show that the MOP-WENO-ACM$k$ scheme can recover the optimal 
convergence rates for different values of $n_{\mathrm{cp}}$ in smooth
regions.

\begin{theorem}
When $\mathrm{CFS}_{0} \ll \widetilde{d}_{0}$ and $\mathrm{CFS}_{1}
\gg \widetilde{d}_{r-1}$, for $\forall n_{\mathrm{cp}} < r - 1$, the
$(2r-1)$th-order MOP-WENO-ACM$k$ scheme can achieve the optimal 
convergence rate of accuracy if the mapping function 
$g^{\mathrm{MOP-ACM}k}(\omega)$ is applied to the weights of the 
$(2r-1)$th-order WENO-JS scheme.
\label{theorem:convergenceRates}
\end{theorem}

We can prove Theorem \ref{theorem:convergenceRates} by employing the
Taylor series analysis and using Lemma 
\ref{lemma:mappingFunction:MOP-ACMk} of this paper and Lemma 1 and 
Lemma 2 in the statement of page 456 to 457 in \cite{WENO-IM}. 
The detailed proof process is almost identical to the one in 
\cite{WENO-M}.

We can get the fifth-order MOP-WENO-ACM$k$ scheme by choosing $r=3$ 
in Eq.(\ref{eq:mapping_MOP-ACMk})(\ref{eq:nonlinearWeightsX}). In 
this case, it is trivial to obtain the arrays $\mathcal{D},
\widetilde{\mathcal{D}}$, that is, $\mathcal{D}=\{0.1,0.6,0.3\},
\widetilde{\mathcal{D}}=\{0.1,0.3,0.6\}$. In other words, we have
$\widetilde{d}_{0}=0.1,\widetilde{d}_{1}=0.3,\widetilde{d}_{2}=0.6$. 
Then, we can write the mapping function explicitly as follows
\begin{equation}
g^{\mathrm{MOP-ACM}k}(\omega) = \left\{
\begin{array}{lll}
k_{0} \omega,& {} & \omega \in \hat{\Omega}_{1}, \\
0.1,         & {} & \omega \in \hat{\Omega}_{2}, \\
0.3,         & {} & \omega \in \hat{\Omega}_{3}, \\
0.6,         & {} & \omega \in \hat{\Omega}_{4}, \\
1 - k_{1}(1 - \omega), & {} & \omega \in \hat{\Omega}_{5},
\end{array}
\right.
\label{eq:mapping_MOP-ACMk:5th}
\end{equation}
where $\hat{\Omega}_{1} = [0, \mathrm{CFS}_{0}), \hat{\Omega}_{2} = 
[\mathrm{CFS}_{0},0.2), \hat{\Omega}_{3} = [0.2, 0.45),
\hat{\Omega}_{4} = [0.45, \mathrm{CFS}_{1}], \hat{\Omega}_{5} = \Big(
\mathrm{CFS}_{1},1\Big]$, $0 < \mathrm{CFS}_{0} \leq 0.1, 0.6 \leq 
\mathrm{CFS}_{1}<1$, $k_{0}\in \Bigg[0,\dfrac{0.1}{\mathrm{CFS}_{0}}
\Bigg],k_{1} \in \Bigg[0, \dfrac{0.4}{1 - \mathrm{CFS}_{1}}\Bigg]$.

In Fig. \ref{fig:gOmega:MOP-ACMk}, we plot the curve of 
$g^{\mathrm{MOP-ACM}k}(\omega)$ varying with $\omega$ by taking 
$\mathrm{CFS}_{0} = 0.04, \mathrm{CFS}_{1} = 0.92$ and $k_{0} = 2.5, 
k_{1} = 5$.

\section{Numerical experiments}
\label{NumericalExperiments}
In this section, we compare the numerical performance of the 
MOP-WENO-ACM$k$ scheme with the WENO-JS scheme \cite{WENO-JS} and 
its various versions with mapping, WENO-M \cite{WENO-M}, WENO-PM$6$ 
\cite{WENO-PM}, WENO-IM(2, 0.1) \cite{WENO-IM} and MIP-WENO-ACM$k$ 
which is proposed in subsection \ref{subsec:MIP-WENO-ACMk}. In all 
the numerical experiments below, MOP-WENO-ACM$k$ refers to the 
definition in Eq.(\ref{eq:nonlinearWeightsX}) and 
Eq.(\ref{eq:mapping_MOP-ACMk:5th}) with $k_{0}=k_{1}=0, 
\mathrm{CFS}_{0} = 0.01, \mathrm{CFS}_{1} = 0.94$, and the 
parameters in the MIP-WENO-ACM$k$ scheme are chosen to be $k_{s} = 0,
\mathrm{CFS}_{s} = d_{s}/10$.

The numerical presentation of this section starts with the accuracy 
test of the one-dimensional linear advection equation with four 
kinds of initial conditions, followed by the solutions at long 
output times of one-dimensional linear advection equation with two 
kinds of initial condition with discontinuities, and finishes with 
2D calculations on the two-dimensional Riemann problem and the 
shock-vortex interaction. 

\subsection{Accuracy test}
\begin{example}
\bf{(Accuracy test without critical points \cite{WENO-IM})} 
\rm{We consider the one-dimensional linear advection equation 
Eq.(\ref{eq:LAE}) with the periodic boundary condition and the 
following initial condition} 
\label{LAE1}
\end{example} 
\begin{equation}
u(x, 0) = \sin ( \pi x ). 
\label{eq:LAE:IC1}
\end{equation}
It is easy to know that the initial condition in 
Eq.(\ref{eq:LAE:IC1}) has no critical points. The CFL number is set 
to be $(\Delta x)^{2/3}$ to prevent the convergence rates of error 
from being influenced by time advancement. The $L_{1}, L_{2}, 
L_{\infty}$ norms of the error are given as
\begin{equation*}
\displaystyle
\begin{array}{l}
L_{1} = h \cdot \displaystyle\sum\limits_{j} \big\lvert u_{j}^{
\mathrm{exact}} - (u_{h})_{j} \big\rvert, \\
L_{2} = \sqrt{h \cdot \displaystyle\sum\limits_{j} (u_{j}^{
\mathrm{exact}} - (u_{h})_{j})^{2}}, \\
L_{\infty} = \displaystyle\max_{j} \big\lvert u_{j}^{
\mathrm{exact}} - (u_{h})_{j} \big\rvert,
\end{array}
\end{equation*}
where $h = \Delta x$ is the uniform spatial step size, $(u_{h})_{j}$ 
is the numerical solution and $u_{j}^{\mathrm{exact}}$ is the exact 
solution.

The $L_{1}, L_{2}, L_{\infty}$ errors and corresponding convergence 
orders of various considered WENO schemes for Example 1 at output 
time $t = 2.0$ are shown in Table \ref{table_LAE1}. The results of 
the three rows are $L_{1}$-, $L_{2}$- and $L_{\infty}$- norm errors 
and convergence orders in turn (similarly hereinafter). 
Unsurprisingly, the MOP-WENO-ACM$k$ scheme has gained the 
fifth-order convergence like the other considered schemes. It can be 
found that the MOP-WENO-ACM$k$, MIP-WENO-ACM$k$, WENO-M, WENO-PM6 
and WENO-IM(2, 0.1) schemes give more accurate numerical solutions 
than the WENO-JS scheme in general.

\begin{table}[ht]
\begin{scriptsize}
\centering
\caption{Convergence properties of considered schemes solving $u_{t} 
+ u_{x} = 0$ with initial condition $u(x, 0) = \sin(\pi x)$.}
\label{table_LAE1}
\begin{tabular*}{\hsize}
{@{}@{\extracolsep{\fill}}lllllll@{}}
\hline
$N$       
& 10                  & 20                  & 40
& 80                  & 160                 & 320    \\
\hline
WENO-JS 
& 6.18628e-02(-)      & 2.96529e-03(4.3821) & 9.27609e-05(4.9985) 
& 2.89265e-06(5.0031) & 9.03392e-08(5.0009) & 2.82330e-09(4.9999)  \\
\space
& 4.72306e-02(-)      & 2.42673e-03(4.2826) & 7.64332e-05(4.9887) 
& 2.33581e-06(5.0322) & 7.19259e-08(5.0213) & 2.23105e-09(5.0107)  \\
\space
& 4.87580e-02(-)      & 2.57899e-03(4.2408) & 9.05453e-05(4.8320) 
& 2.90709e-06(4.9610) & 8.85753e-08(5.0365) & 2.72458e-09(5.0228)  \\
WENO-M 
& 2.01781e-02(-)      & 5.18291e-04(5.2829) & 1.59422e-05(5.0228) 
& 4.98914e-07(4.9979) & 1.56021e-08(4.9990) & 4.99356e-10(4.9977)  \\
\space
& 1.55809e-02(-)      & 4.06148e-04(5.2616) & 1.25236e-05(5.0193) 
& 3.91875e-07(4.9981) & 1.22541e-08(4.9991) & 3.83568e-10(4.9976)  \\
\space
& 1.47767e-02(-)      & 3.94913e-04(5.2256) & 1.24993e-05(4.9816) 
& 3.91808e-07(4.9956) & 1.22538e-08(4.9988) & 3.83541e-10(4.9977)  \\
WENO-PM6 
& 1.74869e-02(-)      & 5.02923e-04(5.1198) & 1.59130e-05(4.9821) 
& 4.98858e-07(4.9954) & 1.56020e-08(4.9988) & 4.88355e-10(4.9977)  \\
\space
& 1.35606e-02(-)      & 3.95215e-04(5.1006) & 1.25010e-05(4.9825) 
& 3.91831e-07(4.9957) & 1.22541e-08(4.9989) & 3.83568e-10(4.9976)  \\
\space
& 1.27577e-02(-)      & 3.94515e-04(5.0151) & 1.24960e-05(4.9805) 
& 3.91795e-07(4.9952) & 1.22538e-08(4.9988) & 3.83543e-10(4.9977)  \\
WENO-IM(2,0.1) 
& 1.58051e-02(-)      & 5.04401e-04(4.9697) & 1.59160e-05(4.9860) 
& 4.98863e-07(4.9957) & 1.56020e-08(4.9988) & 4.88355e-10(4.9977)  \\
\space
& 1.23553e-02(-)      & 3.96236e-04(4.9626) & 1.25033e-05(4.9860) 
& 3.91836e-07(4.9959) & 1.22541e-08(4.9989) & 3.83568e-10(4.9976)  \\
\space
& 1.19178e-02(-)      & 3.94458e-04(4.9171) & 1.24963e-05(4.9803) 
& 3.91797e-07(4.9953) & 1.22538e-08(4.9988) & 3.83547e-10(4.9977)  \\
MIP-WENO-ACM$k$ 
& 1.52184e-02(-)      & 5.02844e-04(4.9196) & 1.59130e-05(4.9818) 
& 4.98858e-07(4.9954) & 1.56020e-08(4.9988) & 4.88355e-10(4.9977)  \\
\space
& 1.19442e-02(-)      & 3.95138e-04(4.9178) & 1.25010e-05(4.9822) 
& 3.91831e-07(4.9957) & 1.22541e-08(4.9989) & 3.83568e-10(4.9976)  \\
\space
& 1.17569e-02(-)      & 3.94406e-04(4.8977) & 1.24960e-05(4.9801) 
& 3.91795e-07(4.9952) & 1.22538e-08(4.9988) & 3.83543e-10(4.9977)  \\
MOP-WENO-ACM$k$ 
& 3.29609e-02(-)      & 5.02844e-04(6.0345) & 1.59130e-05(4.9818) 
& 4.98858e-07(4.9954) & 1.56020e-08(4.9988) & 4.88355e-10(4.9977)  \\
\space
& 2.72363e-02(-)      & 3.95138e-04(6.1070) & 1.25010e-05(4.9822) 
& 3.91831e-07(4.9957) & 1.22541e-08(4.9989) & 3.83568e-10(4.9976)  \\
\space
& 2.70295e-02(-)      & 3.94406e-04(6.0987) & 1.24960e-05(4.9801) 
& 3.91795e-07(4.9952) & 1.22538e-08(4.9988) & 3.83543e-10(4.9977)  \\
\hline
\end{tabular*}
\end{scriptsize}
\end{table}

\begin{example}
\bf{(Accuracy test with first-order critical points \cite{WENO-M})} 
\rm{We consider the one-dimensional linear advection equation 
Eq.(\ref{eq:LAE}) with the periodic boundary condition and the 
following initial condition}
\label{LAE2}
\end{example}
\begin{equation}
u(x, 0) = \sin \bigg( \pi x - \dfrac{\sin(\pi x)}{\pi} \bigg). 
\label{eq:LAE:IC2}
\end{equation}
The particular initial condition Eq.(\ref{eq:LAE:IC2}) has two 
first-order critical points, which both have a non-vanishing third 
derivative. As mentioned earlier, the CFL number is set to be 
$(\Delta x)^{2/3}$.

Table \ref{table_LAE2} shows the $L_{1}, L_{2}, L_{\infty}$ errors 
and corresponding convergence orders of the various considered WENO 
schemes at output time $t = 2.0$. We can see that the $L_{\infty}$ 
convergence order of the WENO-JS scheme drops by almost 2 orders 
leading to an overall accuracy loss shown with $L_{1}$ and $L_{2}$
convergence orders. It is evident that the MOP-WENO-ACM$k$, 
MIP-WENO-ACM$k$, WENO-M, WENO-PM6 and WENO-IM(2,0.1) schemes can 
retain the optimal orders even in the presence of critical points. 
It is noteworthy that when the grid number is too small, like 
$N \leq 40$, in terms of accuracy, the MOP-WENO-ACM$k$ scheme 
provides less accurate results than those of the MIP-WENO-ACM$k$ 
scheme. The cause of this kind of accuracy loss is that the mapping 
function of the MOP-WENO-ACM$k$ scheme has narrower optimal weight 
intervals (standing for the intervals about $\omega = d_{s}$ over 
which the mapping process attempts to use the corresponding optimal 
weights, see \cite{WENO-MAIMi,WENO-ACM}) than the MIP-WENO-ACM$k$ 
scheme. However, this issue can surely be overcome by increasing the 
grid number. Therefore, we can find that, as expected, the 
MOP-WENO-ACM$k$ scheme gives equally accurate numerical solutions 
like those of the MIP-WENO-ACM$k$ scheme when the grid number
$N \geq 80$.

\begin{table}[ht]
\begin{scriptsize}
\centering
\caption{Convergence properties of considered schemes solving $u_{t} 
+ u_{x} = 0$ with initial condition $u(x, 0) = \sin(\pi x - \sin(\pi 
x)/\pi )$.}
\label{table_LAE2}
\begin{tabular*}{\hsize}
{@{}@{\extracolsep{\fill}}lllllll@{}}
\hline
$N$       
& 10                  & 20                  & 40
& 80                  & 160                 & 320    \\
\hline
WENO-JS 
& 1.24488e-01(-)      & 1.01260e-02(3.6199) & 7.22169e-04(3.8096) 
& 3.42286e-05(4.3991) & 1.58510e-06(4.4326) & 7.95517e-08(4.3165)  \\
\space
& 1.09463e-01(-)      & 8.72198e-03(3.6496) & 6.76133e-04(3.6893) 
& 3.63761e-05(4.2162) & 2.29598e-06(3.9858) & 1.68304e-07(3.7700)  \\
\space
& 1.24471e-01(-)      & 1.43499e-02(3.1167) & 1.09663e-03(3.7099) 
& 9.02485e-05(3.6030) & 8.24022e-06(3.4531) & 8.31702e-07(3.3085)  \\
WENO-M 
& 7.53259e-02(-)      & 3.70838e-03(4.3443) & 1.45082e-04(4.6758) 
& 4.80253e-06(4.9169) & 1.52120e-07(4.9805) & 4.77083e-09(4.9948)  \\
\space
& 6.39017e-02(-)      & 3.36224e-03(4.2484) & 1.39007e-04(4.5962) 
& 4.52646e-06(4.9406) & 1.42463e-07(4.9897) & 4.45822e-09(4.9980)  \\
\space
& 7.49250e-02(-)      & 5.43666e-03(3.7847) & 2.18799e-04(4.6350) 
& 6.81451e-06(5.0049) & 2.14545e-07(4.9893) & 6.71080e-09(4.9987)  \\
WENO-PM6 
& 9.51313e-02(-)      & 4.82173e-03(4.3023) & 1.55428e-04(4.9552) 
& 4.87327e-06(4.9952) & 1.52750e-07(4.9956) & 4.77729e-09(4.9988)  \\
\space
& 7.83600e-02(-)      & 4.29510e-03(4.1894) & 1.43841e-04(4.9001) 
& 4.54036e-06(4.9855) & 1.42488e-07(4.9939) & 4.45807e-09(4.9983)  \\
\space
& 9.32356e-02(-)      & 5.91037e-03(3.9796) & 2.09540e-04(4.8180) 
& 6.83270e-06(4.9386) & 2.14532e-07(4.9932) & 6.71079e-09(4.9986)  \\
WENO-IM(2,0.1) 
& 8.38131e-02(-)      & 4.30725e-03(4.2823) & 1.51327e-04(4.8310) 
& 4.85592e-06(4.9618) & 1.52659e-07(4.9914) & 4.77654e-09(4.9982)  \\
\space
& 6.71285e-02(-)      & 3.93700e-03(4.0918) & 1.41737e-04(4.7958) 
& 4.53602e-06(4.9656) & 1.42479e-07(4.9926) & 4.45805e-09(4.9982)  \\
\space
& 7.62798e-02(-)      & 5.84039e-03(3.7072) & 2.10531e-04(4.7940) 
& 6.82606e-06(4.9468) & 2.14534e-07(4.9918) & 6.71079e-09(4.9986)  \\
MIP-WENO-ACM$k$ 
& 8.75629e-02(-)      & 4.39527e-03(4.3163) & 1.52219e-04(4.8517) 
& 4.86436e-06(4.9678) & 1.52735e-07(4.9931) & 4.77728e-09(4.9987)  \\
\space
& 6.98131e-02(-)      & 4.02909e-03(4.1150) & 1.42172e-04(4.8247) 
& 4.53770e-06(4.9695) & 1.42486e-07(4.9931) & 4.45807e-09(4.9983)  \\
\space
& 7.91292e-02(-)      & 5.89045e-03(3.7478) & 2.09893e-04(4.8107) 
& 6.83017e-06(4.9416) & 2.14533e-07(4.9926) & 6.71079e-09(4.9986)  \\
MOP-WENO-ACM$k$ 
& 9.08634e-02(-)      & 7.09246e-03(3.6793) & 2.59429e-04(4.7729) 
& 4.86436e-06(5.7369) & 1.52735e-07(4.9931) & 4.77728e-09(4.9987)  \\
\space
& 7.58160e-02(-)      & 6.88532e-03(3.4609) & 2.51208e-04(4.7766) 
& 4.53770e-06(5.7908) & 1.42486e-07(4.9931) & 4.45807e-09(4.9983)  \\
\space
& 9.29135e-02(-)      & 1.01479e-02(3.1947) & 4.03069e-04(4.6540) 
& 6.83017e-06(5.8830) & 2.14533e-07(4.9926) & 6.71079e-09(4.9986)  \\
\hline
\end{tabular*}
\end{scriptsize}
\end{table}

\begin{example}
\bf{(High resolution performance test with high-order critical
points)} 
\rm{We consider the one-dimensional linear advection equation 
Eq.(\ref{eq:LAE}) with the periodic boundary condition and the 
following initial condition \cite{WENO-PM}} 
\label{LAE4}
\end{example} 
\begin{equation}
u(x, 0) = \sin^{9} ( \pi x ) 
\label{eq:LAE:IC4}
\end{equation}
It is trivial to verify that the initial condition in 
Eq.(\ref{eq:LAE:IC4}) has high-order critical points. Again, the CFL 
number is set to be $(\Delta x)^{2/3}$. 

We calculate this problem using the MOP-WENO-ACM$k$, MIP-WENO-ACM$k$
, WENO-JS and WENO-M schemes. Table \ref{table_LAE4} shows the $L_{1}
, L_{2}, L_{\infty}$ errors of these WENO schemes at several output 
times with a uniform mesh size of $\Delta x = 1/200$. We also 
present the corresponding increased errors (in percentage) compared 
to the errors of the MIP-WENO-ACM$k$ scheme which gives the most 
accurate results. Taking the $L_{1}$-norm error as an example, its 
increased error at output time $t$ is calculated by
$\frac{L_{1}^{\mathrm{X}}(t)-L_{1}^{*}(t)}{L_{1}^{*}(t)}\times100\%$
, where $L_{1}^{*}(t)$ and $L_{1}^{\mathrm{X}}(t)$ are the
$L_{1}$-norm error of the MIP-WENO-ACM$k$ scheme and the WENO-X 
scheme (X = WENO-JS, WENO-M, or MOP-WENO-ACM$k$) at output time $t$. 
Clearly, the WENO-JS scheme gives the largest increased errors for 
both short and long output times. At short output times, like $t 
\leq 100$, the solutions computed by the WENO-M scheme are most 
close to those of the MIP-WENO-ACM$k$ scheme, leading to the 
smallest increased errors. However, when the output time increases 
to $t \geq 200$, the solutions computed by the MOP-WENO-ACM$k$ 
scheme are most close to those of the MIP-WENO-ACM$k$ scheme. 
Furthermore, when the output time is large, the errors of the WENO-M 
scheme increase significantly leading to evidently larger 
increased errors, while the increased errors of the MOP-WENO-ACM$k$ 
scheme do not get larger as it provides comparably small errors. The 
errors of the MOP-WENO-ACM$k$ are not as small as those of the 
MIP-WENO-ACM$k$ scheme. The cause of this kind of accuracy loss is 
the same as stated in Example \ref{LAE2} that the mapping function 
of the MOP-WENO-ACM$k$ scheme has narrower optimal weight intervals 
than the MIP-WENO-ACM$k$ scheme. As mentioned earlier, one can 
surely address this issue by increasing the grid number. In order to 
verify this, we calculate this problem using the same schemes at the 
same output times but with a larger grid number of $N = 800$, and 
the results are shown in Table \ref{table_LAE4-N800}. From Table 
\ref{table_LAE4-N800}, we can see that the errors of the MOP-WENO-ACM
$k$ scheme are closer to those of the MIP-WENO-ACM$k$ scheme when 
the grid number increases from $N=200$ to $N=800$, resulting in the 
significantly decreasing of the increased errors. Moreover, for long 
output times, the increased errors of the MOP-WENO-ACM$k$ scheme are 
much smaller than those of the WENO-JS and WENO-M schemes. Actually, 
it is an important advantage of the MOP-WENO-ACM$k$ scheme that can 
maintain comparably high resolution for long output times. In the 
next subsection we will demonstrate this again.

\begin{figure}[ht]
\centering
  \includegraphics[height=0.39\textwidth]
  {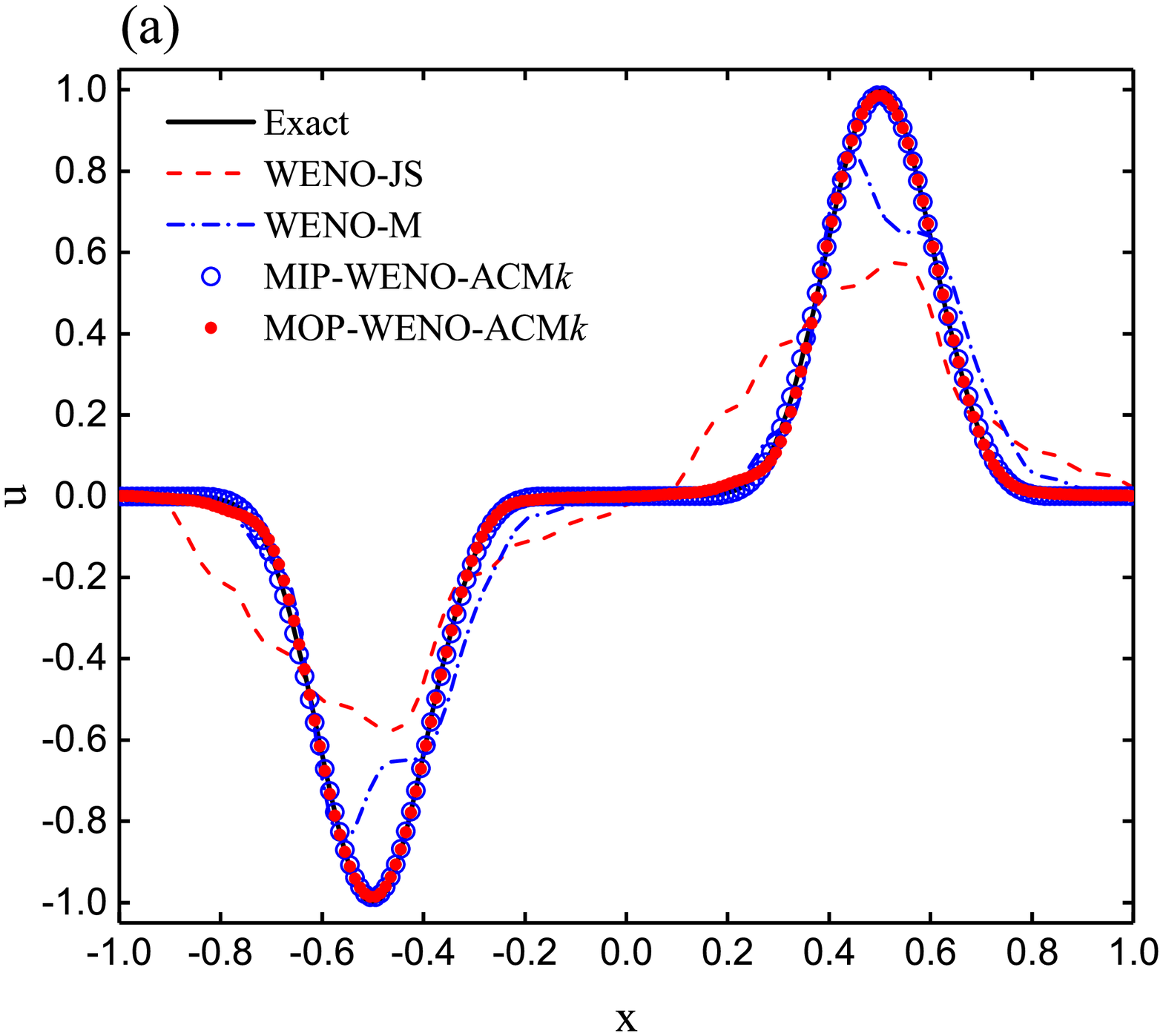}
  \includegraphics[height=0.39\textwidth]
  {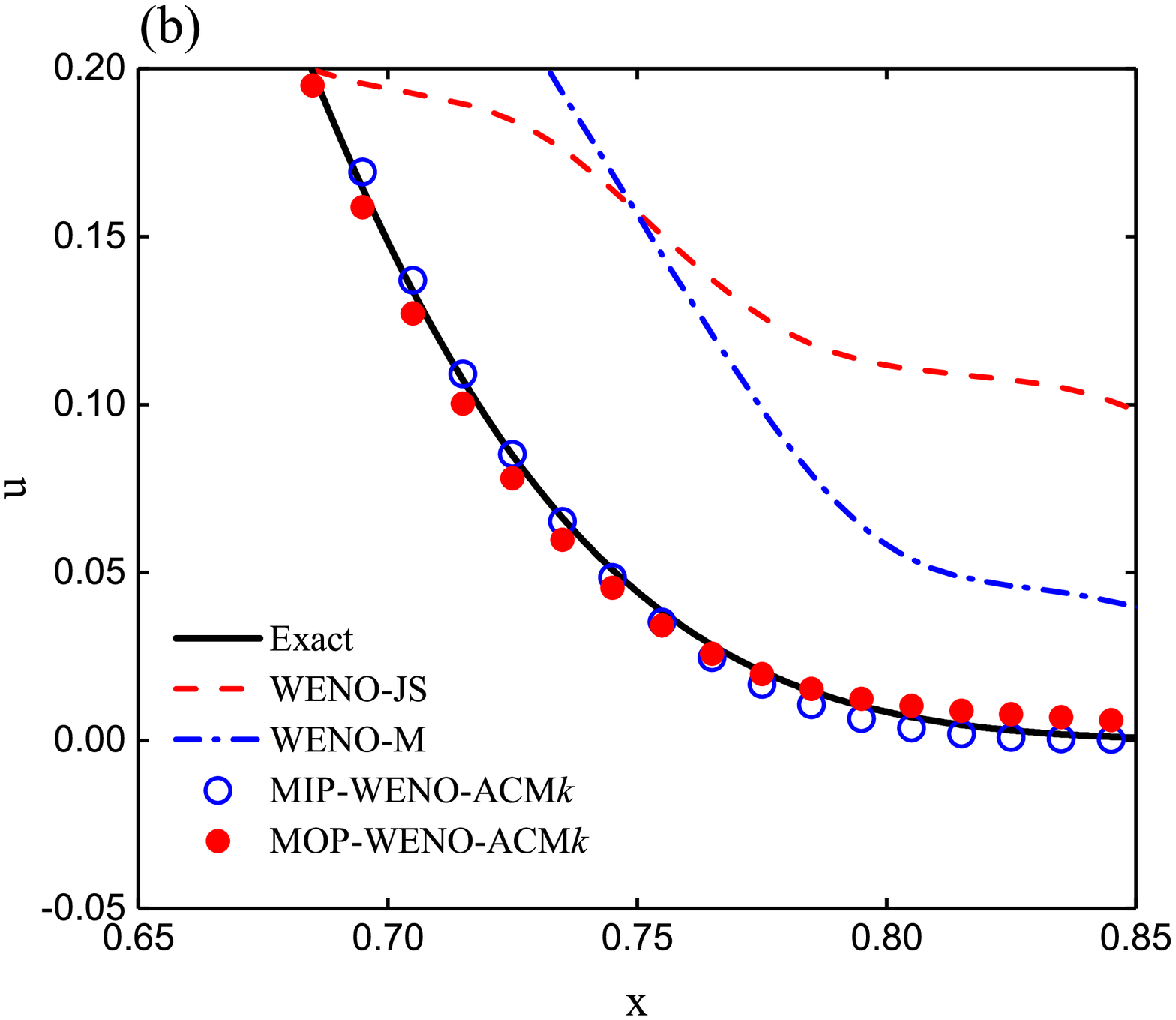}
  \caption{Performance of the MOP-WENO-ACM$k$, MIP-WNEO-ACM$k$, 
  WENO-JS and WENO-M schemes for Example \ref{LAE4} at output time 
  $t = 1000$ with a uniform mesh size of $\Delta x= 1/200$.}
\label{fig:ex:LAE4}
\end{figure}

Fig. \ref{fig:ex:LAE4} shows the performance of the WENO-JS, 
WENO-M, MIP-WENO-ACM$k$ and MOP-WENO-ACM$k$ schemes at output time
$t = 1000$ with a uniform mesh size of $\Delta x= 1/200$. Clearly, 
the MIP-WENO-ACM$k$ and MOP-WENO-ACM$k$ schemes give the highest 
resolution, followed by the WENO-M scheme, whose resolution 
decreases significantly. The WENO-JS scheme shows the lowest 
resolution.

\begin{table}[ht]
\begin{scriptsize}
\centering
\caption{Performance of various considered schemes solving $u_{t} + 
u_{x} = 0$ with $u(x, 0) = \sin^{9} (\pi x), \Delta x= 1/200$.}
\label{table_LAE4}
\begin{tabular*}{\hsize}
{@{}@{\extracolsep{\fill}}llllllll@{}}
\hline
\space      & MIP-WENO-ACM$k$             
            & \multicolumn{2}{l}{MOP-WENO-ACM$k$}
            & \multicolumn{2}{l}{WENO-JS}
            & \multicolumn{2}{l}{WENO-M}  \\
\cline{2-2}  \cline{3-4}  \cline{5-6}  \cline{7-8}            
Scheme      & Errors 
            & Errors & Increased errors
            & Errors & Increased errors 
            & Errors & Increased errors    \\         
\hline
$t = 1$     & 8.43356e-06            
            & 1.73735e-05 & 106.00\% 
            & 3.87826e-05 & 359.86\% 
            & 8.84565e-06 & 4.89\%  \\
\space      & 8.20366e-06            
            & 1.85069e-05 & 125.59\%
            & 3.62689e-05 & 342.11\% 
            & 8.31248e-06 & 1.33\%  \\
\space      & 1.38389e-05            
            & 6.40784e-05 & 363.03\%
            & 6.69118e-05 & 383.51\% 
            & 1.38461e-05 & 0.05\%  \\
$t = 10$    & 8.42873e-05            
            & 1.55900e-04 & 84.96\%
            & 3.86931e-04 & 359.06\% 
            & 8.90890e-05 & 5.70\%  \\
\space      & 8.19107e-05            
            & 1.63558e-04 & 99.68\%
            & 3.52611e-04 & 330.48\%    
            & 8.32089e-05 & 1.58\%  \\
\space      & 1.38205e-04            
            & 5.22964e-04 & 278.40\%
            & 5.36940e-04 & 288.51\%    
            & 1.38348e-04 & 0.10\%  \\
$t = 30$    & 2.52378e-04            
            & 5.66584e-04 & 124.50\%
            & 1.17988e-03 & 367.51\%    
            & 2.73430e-04 & 8.34\%  \\
\space      & 2.45090e-04            
            & 6.35696e-04 & 159.37\%
            & 1.06511e-03 & 334.58\%    
            & 2.51737e-04 & 2.71\%  \\
\space      & 4.13398e-04            
            & 1.99465e-03 & 382.50\%
            & 1.58134e-03 & 282.52\%    
            & 4.13887e-04 & 0.12\%  \\
$t = 50$    & 4.19825e-04            
            & 1.09581e-03 & 161.02\%
            & 2.05488e-03 & 389.46\%    
            & 4.81901e-04 & 14.79\%  \\
\space      & 4.07429e-04            
            & 1.26383e-03 & 210.20\%
            & 1.84782e-03 & 353.53\%    
            & 4.39983e-04 & 7.99\%\\
\space      & 6.86983e-04            
            & 3.84112e-03 & 459.13\%
            & 2.69500e-03 & 292.30\%    
            & 6.87879e-04 & 0.13\%  \\
$t = 100$   & 8.35747e-04            
            & 2.72470e-03 & 226.02\%
            & 5.42288e-03 & 548.87\%     
            & 1.29154e-03 & 54.54\%  \\
\space      & 8.09679e-04            
            & 3.23726e-03 & 299.82\%
            & 5.17716e-03 & 539.41\%
            & 1.28740e-03 & 59.00\%  \\
\space      & 1.36404e-03            
            & 9.83147e-03 & 620.76\%
            & 1.20056e-02 & 780.15\%
            & 3.32665e-03 & 143.88\% \\
$t = 200$   & 1.65557e-03            
            & 4.11740e-03 & \textbf{148.70\%}
            & 2.35657e-02 & \textbf{1323.42\%}
            & 5.74021e-03 & \textbf{246.72\%} \\
\space      & 1.59929e-03            
            & 3.71649e-03 & \textbf{132.38\%}
            & 2.68753e-02 & \textbf{1580.45\%}
            & 7.66721e-03 & \textbf{379.41\%} \\
\space      & 2.68955e-03            
            & 6.66166e-03 & \textbf{147.69\%}
            & 6.47820e-02 & \textbf{2308.66\%}
            & 2.37125e-02 & \textbf{781.65\%} \\
$t = 500$   & 3.95849e-03            
            & 8.34435e-03 & \textbf{110.80\%}
            & 1.55650e-01 & \textbf{3832.05\%}
            & 4.89290e-02 & \textbf{1136.05\%} \\
\space      & 3.84802e-03            
            & 7.96980e-03 & \textbf{107.11\%}
            & 1.46859e-01 & \textbf{3716.48\%}
            & 6.23842e-02 & \textbf{1521.20\%} \\
\space      & 6.45564e-03            
            & 1.83215e-02 & \textbf{183.81\%}
            & 2.57663e-01 & \textbf{3891.29\%}
            & 1.78294e-01 & \textbf{2661.83\%} \\
$t = 1000$  & 7.24723e-03            
            & 1.54830e-02 & \textbf{113.64\%}
            & 2.91359e-01 & \textbf{3920.28\%}
            & 1.34933e-01 & \textbf{1761.86\%} \\
\space      & 7.21626e-03            
            & 1.50017e-02 & \textbf{107.89\%}
            & 2.66692e-01 & \textbf{3595.71\%}
            & 1.46524e-01 & \textbf{1930.47\%} \\
\space      & 1.21593e-02            
            & 3.16523e-02 & \textbf{160.31\%}
            & 4.44664e-01 & \textbf{3556.96\%}
            & 3.17199e-01 & \textbf{2508.69\%} \\
\hline
\end{tabular*}
\end{scriptsize}
\end{table}

\begin{table}[ht]
\begin{scriptsize}
\centering
\caption{Performance of various considered schemes solving $u_{t} + 
u_{x} = 0$ with $u(x, 0) = \sin^{9} (\pi x), \Delta x= 1/800$.}
\label{table_LAE4-N800}
\begin{tabular*}{\hsize}
{@{}@{\extracolsep{\fill}}llllllll@{}}
\hline
\space      & MIP-WENO-ACM$k$             
            & \multicolumn{2}{l}{MOP-WENO-ACM$k$}
            & \multicolumn{2}{l}{WENO-JS}
            & \multicolumn{2}{l}{WENO-M}  \\
\cline{2-2}  \cline{3-4}  \cline{5-6}  \cline{7-8}            
Scheme      & Errors 
            & Errors & Increased errors
            & Errors & Increased errors 
            & Errors & Increased errors    \\         
\hline
$t = 1$     & 8.28796e-09            
            & 8.33502e-09 & 0.57\% 
            & 4.23539e-08 & 411.03\% 
            & 8.28916e-09 & 0.01\%  \\
\space      & 8.06772e-09            
            & 8.07086e-09 & 0.04\% 
            & 3.78411e-08 & 369.04\% 
            & 8.06776e-09 & 0.00\%  \\
\space      & 1.36173e-08            
            & 1.36173e-08 & 0.00\% 
            & 7.82744e-08 & 474.82\% 
            & 1.36173e-08 & 0.00\%  \\
$t = 10$    & 8.28794e-08            
            & 8.47930e-08 & 2.31\% 
            & 4.23531e-07 & 411.02\% 
            & 8.28912e-08 & 0.01\%  \\
\space      & 8.06769e-08            
            & 8.11193e-08 & 0.55\% 
            & 3.76804e-07 & 367.05\% 
            & 8.06774e-08 & 0.00\%  \\
\space      & 1.36172e-07            
            & 1.36172e-07 & 0.00\% 
            & 6.95290e-07 & 410.60\% 
            & 1.36173e-07 & 0.00\%  \\
$t = 30$    & 2.48654e-07            
            & 2.73030e-07 & 9.80\% 
            & 1.27070e-06 & 411.03\% 
            & 2.48687e-07 & 0.01\%  \\
\space      & 2.42057e-07            
            & 2.58215e-07 & 6.68\% 
            & 1.12915e-06 & 366.48\% 
            & 2.42059e-07 & 0.00\%  \\
\space      & 4.08575e-07            
            & 5.62109e-07 & 37.58\% 
            & 1.98978e-06 & 387.00\% 
            & 4.08577e-07 & 0.00\%  \\
$t = 50$    & 4.14434e-07            
            & 5.04422e-07 & 21.71\% 
            & 2.11778e-06 & 411.01\% 
            & 4.14491e-07 & 0.01\%  \\
\space      & 4.03460e-07            
            & 5.02965e-07 & 24.66\% 
            & 1.88080e-06 & 366.17\% 
            & 4.03463e-07 & 0.00\%  \\
\space      & 6.81002e-07            
            & 1.82747e-06 & 168.35\% 
            & 3.25711e-06 & 378.28\% 
            & 6.81005e-07 & 0.00\%  \\
$t = 100$   & 8.28891e-07            
            & 9.73202e-07 & 17.41\% 
            & 4.74028e-06 & 471.88\% 
            & 8.29015e-07 & 0.01\%  \\
\space      & 8.06973e-07            
            & 9.07161e-07 & 12.42\% 
            & 4.32403e-06 & 435.83\% 
            & 8.06977e-07 & 0.00\%  \\
\space      & 1.36206e-06            
            & 1.79160e-06 & 31.54\% 
            & 1.09481e-05 & 703.79\% 
            & 1.36207e-06 & 0.00\%  \\
$t = 200$   & 1.65782e-06            
            & 1.78369e-06 & \textbf{7.59\%} 
            & 7.29285e-05 & \textbf{4299.06\%} 
            & 2.27991e-06 & \textbf{37.52\%}  \\
\space      & 1.61399e-06            
            & 1.64768e-06 & \textbf{2.09\%} 
            & 1.60499e-04 & \textbf{9844.24\%} 
            & 2.59031e-06 & \textbf{60.49\%}  \\
\space      & 2.72415e-06            
            & 2.72415e-06 & \textbf{0.00\%} 
            & 9.51604e-04 & \textbf{34832.14\%} 
            & 1.22731e-05 & \textbf{350.53\%}  \\
$t = 500$   & 4.14451e-06            
            & 4.84739e-06 & \textbf{16.96\%} 
            & 3.11698e-02 & \textbf{751974.43\%} 
            & 1.41413e-03 & \textbf{34020.56\%}  \\
\space      & 4.03492e-06            
            & 4.47345e-06 & \textbf{10.87\%} 
            & 4.08456e-02 & \textbf{1012202.60\%} 
            & 3.28891e-03 & \textbf{81411.16\%}  \\
\space      & 6.81018e-06            
            & 8.79296e-06 & \textbf{29.11\%} 
            & 8.63989e-02 & \textbf{1268572.78\%} 
            & 1.90785e-02 & \textbf{280046.78\%}  \\
$t = 1000$  & 8.28868e-06            
            & 8.61232e-06 & \textbf{3.90\%} 
            & 1.01278e-01 & \textbf{1221783.34\%} 
            & 1.83325e-02 & \textbf{221075.14\%}  \\
\space      & 8.06938e-06            
            & 8.14436e-06 & \textbf{0.93\%} 
            & 1.13316e-01 & \textbf{1404171.46\%} 
            & 3.30753e-02 & \textbf{409786.51\%}  \\
\space      & 1.36194e-05            
            & 1.36194e-05 & \textbf{0.00\%} 
            & 2.13485e-01 & \textbf{1567406.64\%}
            & 1.38215e-01 & \textbf{1014739.13\%}  \\
\hline
\end{tabular*}
\end{scriptsize}
\end{table}

\begin{example}
\bf{(Accuracy test with discontinuous initial condition)} 
\rm{We consider the SLP modeled by the one-dimensional linear 
advection equation Eq.(\ref{eq:LAE}) with the initial condition 
Eq.(\ref{eq:LAE:SLP}). In this problem, the CFL number is taken to be
$0.1$.}
\label{LAE3}
\end{example} 

Table \ref{table_LAE3} shows the $L_{1}, L_{2}, L_{\infty}$ errors 
and corresponding convergence orders of various considered WENO 
schemes for this example at output times $t = 2$ and $t = 2000$. At 
the short output time $t = 2$, we find that: (1) for all considered 
schemes, the $L_{1}$ and $L_{2}$ orders are approximately $1.0$ and
$0.4$ to $0.5$, respectively, and the $L_{\infty}$ orders are all 
negative; (2) the MIP-WENO-ACM$k$, WENO-M, WENO-PM6 and WENO-IM(2, 
0.1) schemes present more accurate results than the MOP-WENO-ACM$k$ 
and WENO-JS schemes. At the long output time $t = 2000$, we find 
that: (1) for the WENO-JS and WENO-M schemes, the $L_{1}$, $L_{2}$ 
orders decrease to very small values and even become negative; (2) 
however, for the MOP-WENO-ACM$k$, MIP-WENO-ACM$k$, WENO-PM6 and 
WENO-IM(2, 0.1) schemes, their $L_{1}$ orders are clearly larger 
than $1.0$, and their $L_{2}$ orders increase to approximately $0.6$ 
to $0.9$; (3) for all considered schemes, the $L_{\infty}$ orders 
are very small and even become negative; (4) in terms of accuracy, 
on the whole, the MOP-WENO-ACM$k$ scheme produces accurate and 
comparable results as the other considered mapped WENO schemes 
except the WENO-M scheme. However, if we take a closer look, we can 
find that: (1) the resolution of the result computed by the WENO-M 
scheme is significantly lower than that of the MOP-WENO-ACM$k$ 
scheme; (2) the WENO-PM6, WENO-IM(2,0.1) and MIP-WENO-ACM$k$ schemes 
generate spurious oscillations but the MOP-WENO-ACM$k$ scheme does 
not. More results will be presented carefully to demonstrate this in 
the following subsection.

\begin{table}[ht]
\begin{scriptsize}
\centering
\caption{Convergence properties of various considered schemes 
solving $u_{t} + u_{x} = 0$ with initial condition 
Eq.(\ref{eq:LAE:SLP}).}
\label{table_LAE3}
\begin{tabular*}{\hsize}
{@{}@{\extracolsep{\fill}}lllllll@{}}
\hline
\space  &\multicolumn{3}{l}{t = 2}  &\multicolumn{3}{l}{t = 2000}  \\
\cline{2-4}  \cline{5-7}
$N$
& 200            & 400                  & 800 
& 200            & 400                  & 800   \\
\hline
WENO-JS     
& 6.30497e-02(-) & 2.81654e-02(1.2103)  & 1.41364e-02(0.9945)  
& 6.12899e-01(-) & 5.99215e-01(0.0326)  & 5.50158e-01(0.1232)   \\
\space
& 1.08621e-01(-) & 7.71111e-02(0.4943)  & 5.69922e-02(0.4362)  
& 5.08726e-01(-) & 5.01160e-01(0.0216)  & 4.67585e-01(0.1000)   \\
\space
& 4.09733e-01(-) & 4.19594e-01(-0.0343) & 4.28463e-01(-0.0302)  
& 7.99265e-01(-) & 8.20493e-01(-0.0378) & 8.14650e-01(0.0103)   \\
WENO-M     
& 4.77201e-02(-) & 2.23407e-02(1.0949)  & 1.11758e-02(0.9993)  
& 3.81597e-01(-) & 3.25323e-01(0.2302)  & 3.48528e-01(-0.0994)  \\
\space
& 9.53073e-02(-) & 6.91333e-02(0.4632)  & 5.09232e-02(0.4411)  
& 3.59205e-01(-) & 3.12970e-01(0.1988)  & 3.24373e-01(-0.0516)  \\
\space
& 3.94243e-01(-) & 4.05856e-01(-0.0419) & 4.16937e-01(-0.0389)  
& 6.89414e-01(-) & 6.75473e-01(0.0295)  & 6.25645e-01(0.1106)   \\
WENO-PM6     
& 4.66681e-02(-) & 2.13883e-02(1.1256)  & 1.06477e-02(1.0063)  
& 2.17323e-01(-) & 1.05197e-01(1.0467)  & 4.47030e-02(1.2347)   \\
\space
& 9.45566e-02(-) & 6.82948e-02(0.4694)  & 5.03724e-02(0.4391)  
& 2.28655e-01(-) & 1.47518e-01(0.6323)  & 9.34250e-02(0.6590)   \\
\space
& 3.96866e-01(-) & 4.06118e-01(-0.0332) & 4.15277e-01(-0.0322)  
& 5.63042e-01(-) & 5.04977e-01(0.1570)  & 4.71368e-01(0.0994)   \\
WENO-IM(2,0.1)     
& 4.40293e-02(-) & 2.02331e-02(1.1217)  & 1.01805e-02(0.9909)  
& 2.17411e-01(-) & 1.12590e-01(0.9493)  & 5.18367e-02(1.1190)   \\
\space
& 9.19118e-02(-) & 6.68479e-02(0.4594)  & 4.95333e-02(0.4325)  
& 2.30000e-01(-) & 1.64458e-01(0.4839)  & 9.98968e-02(0.7192)   \\
\space
& 3.86789e-01(-) & 3.98769e-01(-0.0441) & 4.09515e-01(-0.0383)  
& 5.69864e-01(-) & 4.82180e-01(0.2410)  & 4.73102e-01(0.02784)   \\
MIP-WENO-ACM$k$     
& 4.45059e-02(-) & 2.03667e-02(1.1278)  & 1.02183e-02(0.9951)  
& 2.21312e-01(-) & 1.10365e-01(1.0038)  & 4.76589e-02(1.2115)   \\
\space
& 9.24356e-02(-) & 6.70230e-02(0.4638)  & 4.96081e-02(0.4341)  
& 2.28433e-01(-) & 1.48498e-01(0.6213)  & 9.40843e-02(0.6584)   \\
\space
& 3.92505e-01(-) & 4.04024e-01(-0.0417) & 4.13511e-01(-0.0335)  
& 5.36242e-01(-) & 5.13503e-01(0.0625)  & 5.15898e-01(-0.0067)   \\
MOP-WENO-ACM$k$     
& 5.56533e-02(-) & 2.79028e-02(0.9961)  & 1.43891e-02(0.9554)  
& 3.83033e-01(-) & 1.77114e-01(1.1128)  & 6.70535e-02(1.4013)   \\
\space
& 9.94223e-02(-) & 7.33101e-02(0.4396)  & 5.51602e-02(0.4104)  
& 3.46814e-01(-) & 1.87369e-01(0.8883)  & 1.09368e-01(0.7767)   \\
\space
& 4.03765e-01(-) & 4.48412e-01(-0.1513) & 4.67036e-01(-0.0587)  
& 7.18464e-01(-) & 5.05980e-01(0.5058)  & 4.80890e-01(0.0734)   \\
\hline
\end{tabular*}
\end{scriptsize}
\end{table}

\subsection{Linear advection examples with discontinuities at long 
output times for comparison}
In this subsection, we will make a further study on calculating 
linear advection examples with discontinuities at long output times 
by various considered WENO schemes. The objective is to demonstrate 
the advantage of the MOP-WENO-ACM$k$ scheme that can obtain high 
resolution and do not generate spurious oscillations, especially for 
long output time simulations.

The one-dimensional linear advection problem Eq.(\ref{eq:LAE}) was 
used in this study. It was solved with the following two initial 
conditions. 

Case 1. (SLP) The initial condition is given by 
Eq.(\ref{eq:LAE:SLP}) in subsection \ref{subsec:importantDiscussion}.

Case 2. (BiCWP) The initial condition is given by
\begin{equation}
\begin{array}{l}
u(x, 0) = \left\{
\begin{array}{ll}
0,   & x \in [-1.0, -0.8] \bigcup (-0.2, 0.2] \bigcup (0.8, 1.0], \\
0.5, & x \in (-0.6, -0.4] \bigcup (0.2, 0.4]  \bigcup (0.6, 0.8], \\
1,   & x \in (-0.8, -0.6] \bigcup (-0.4, -0.2] \bigcup (0.4, 0.6],
\end{array}\right. 
\end{array}
\label{eq:LAE:BiCWP}
\end{equation}
and the periodic boundary condition is used in the two directions.

Case 1 is the SLP used earlier in this paper. Case 2 consists of 
several constant states separated by sharp discontinuities at $x = 
\pm 0.8, \pm 0.6, \pm 0.4, \pm 0.2$. We call Case 2 BiCWP for 
brevity in the presentation as the profile of the exact solution for 
this \textit{\textbf{P}roblem} looks like the 
\textit{\textbf{B}reach \textbf{i}n \textbf{C}ity \textbf{W}all}.

\begin{figure}[ht]
\centering
\includegraphics[height=0.32\textwidth]
{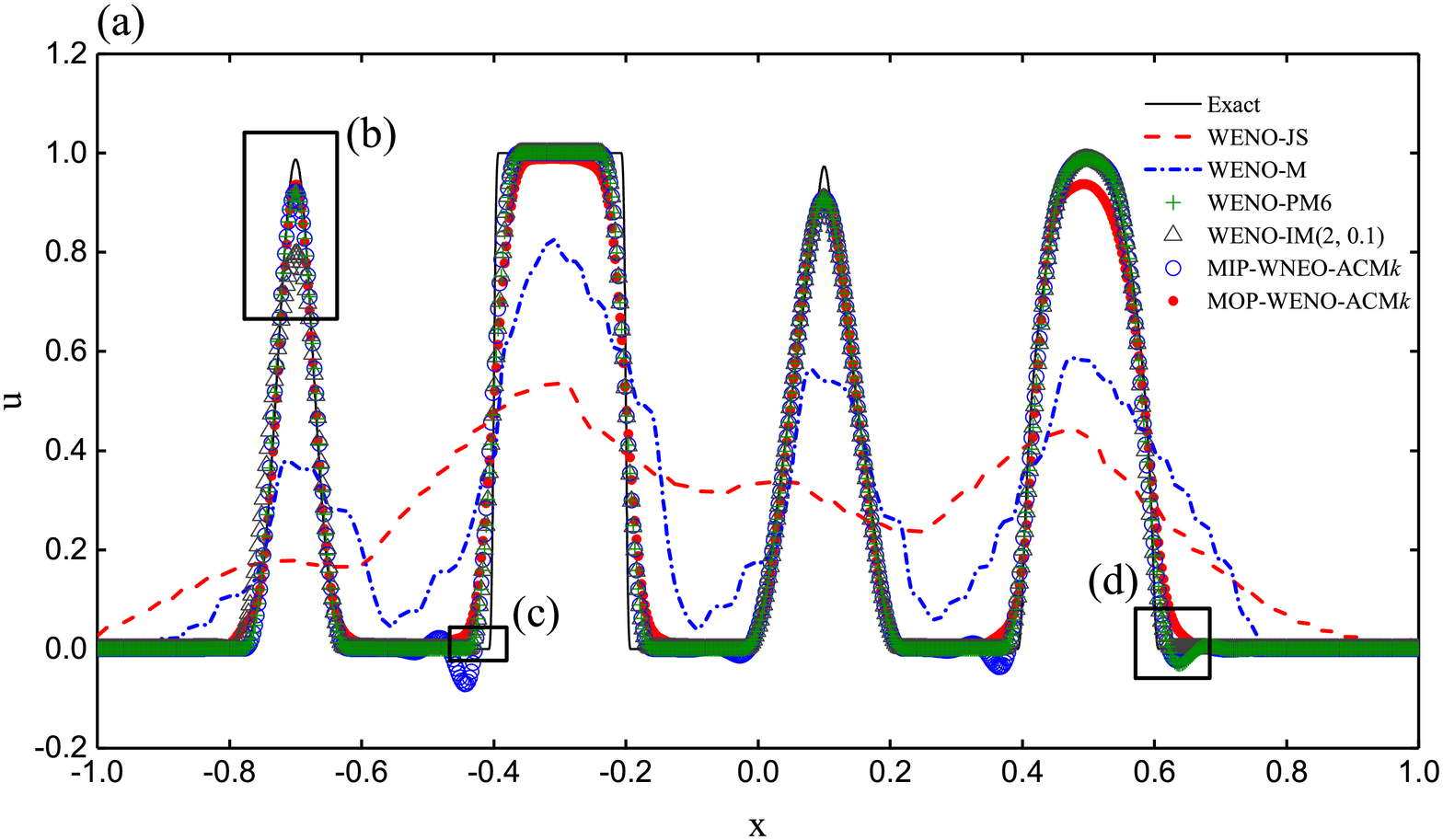}
\includegraphics[height=0.32\textwidth]
{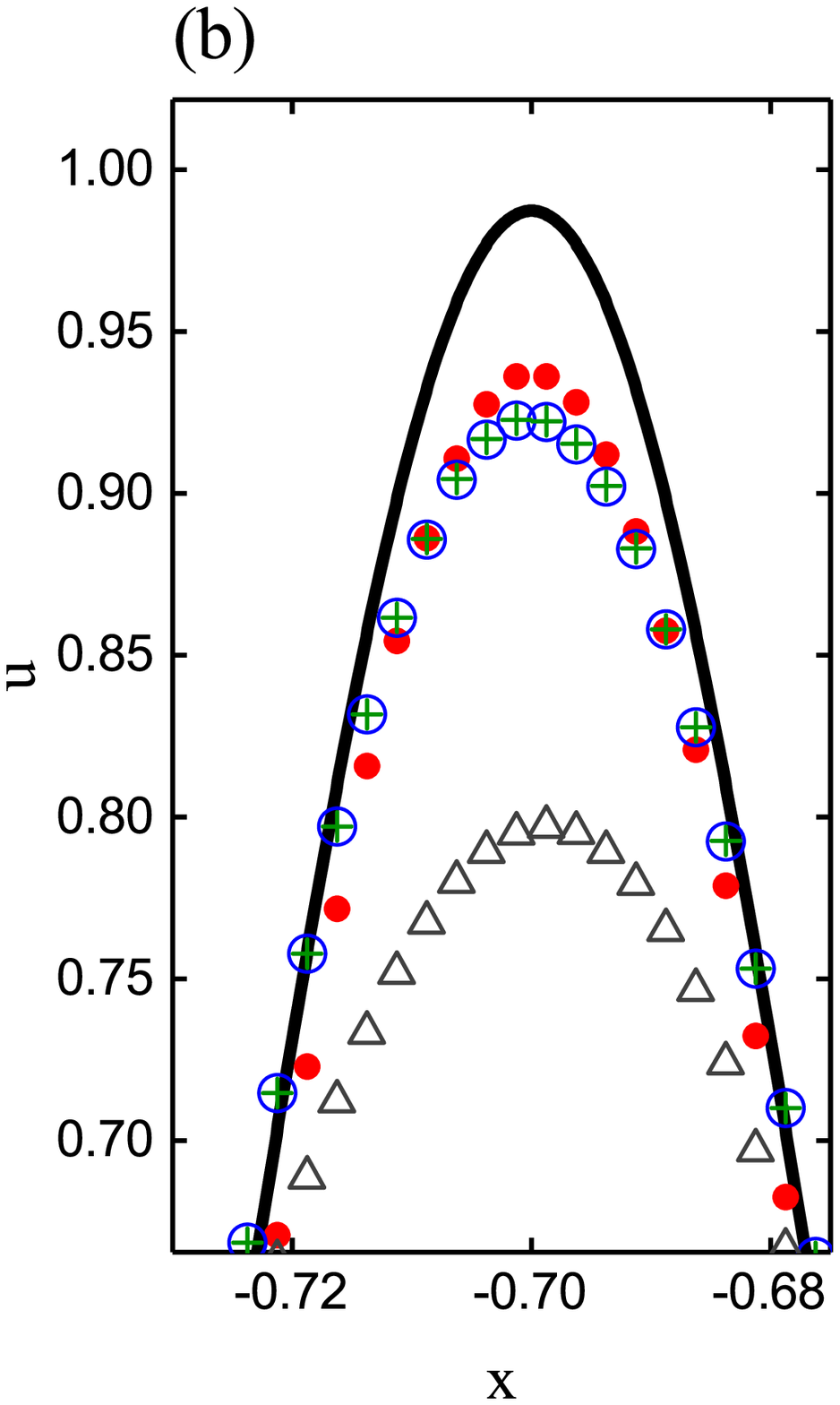}
\includegraphics[height=0.32\textwidth]
{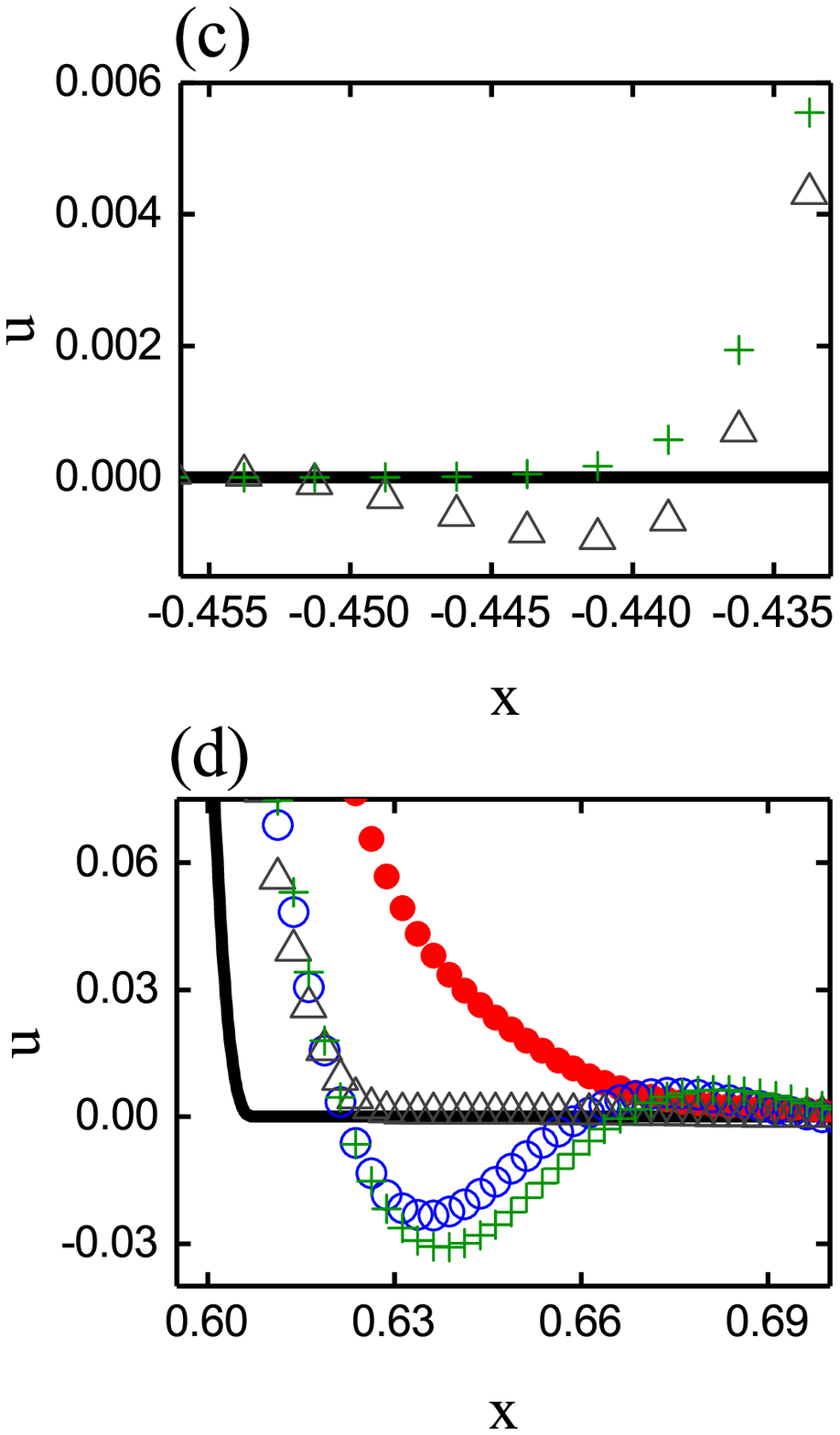}
\caption{Performance of the fifth-order MOP-WENO-ACM$k$, 
MIP-WNEO-ACM$k$, WENO-JS, WENO-M, WENO-PM6 and WENO-IM($2,0.1$) 
schemes for the SLP with $N=800$ at output time $t = 2000$.}
\label{fig:SLP:N800}
\end{figure}

\begin{figure}[ht]
\centering
\includegraphics[height=0.32\textwidth]
{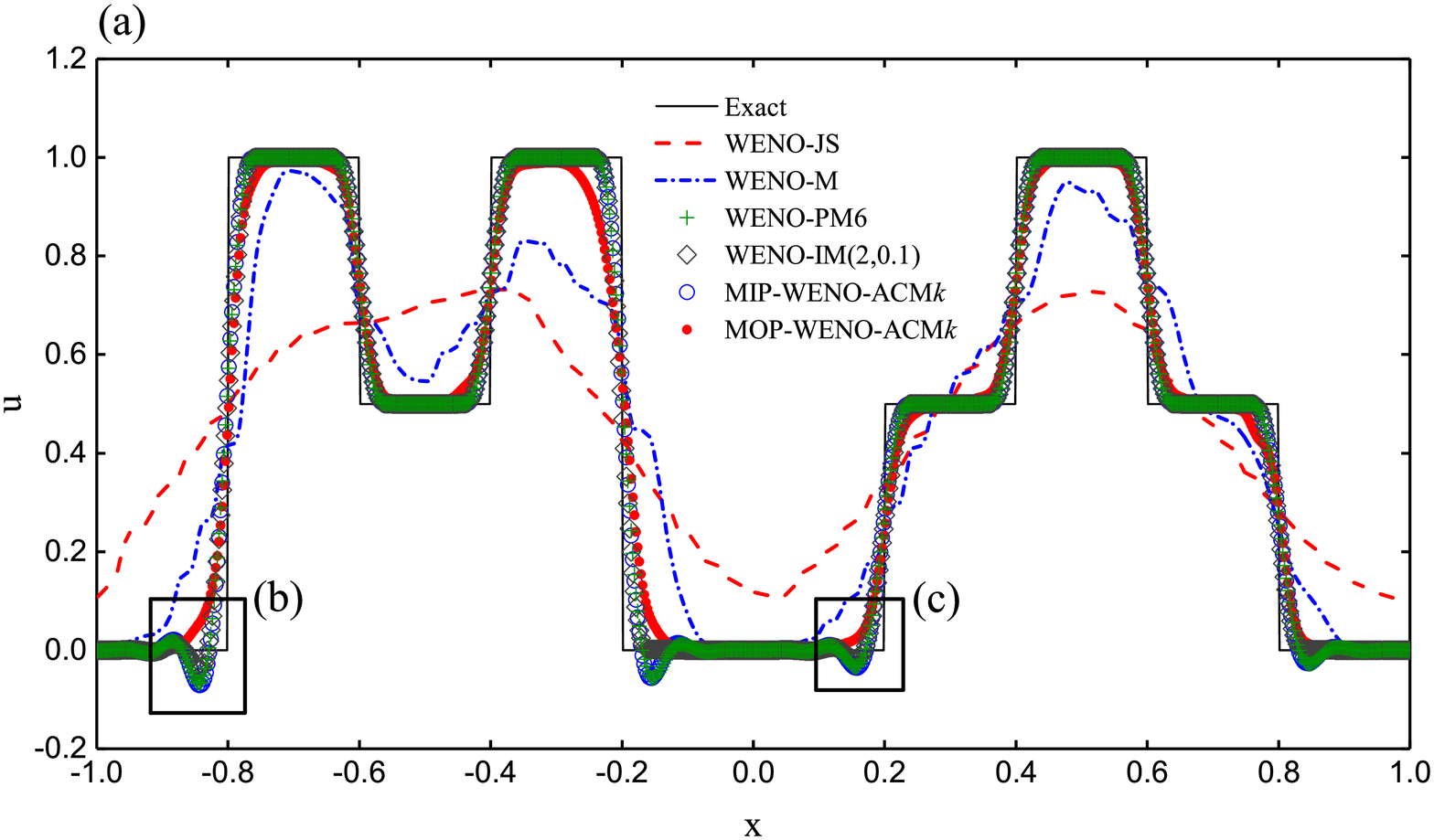}
\includegraphics[height=0.32\textwidth]
{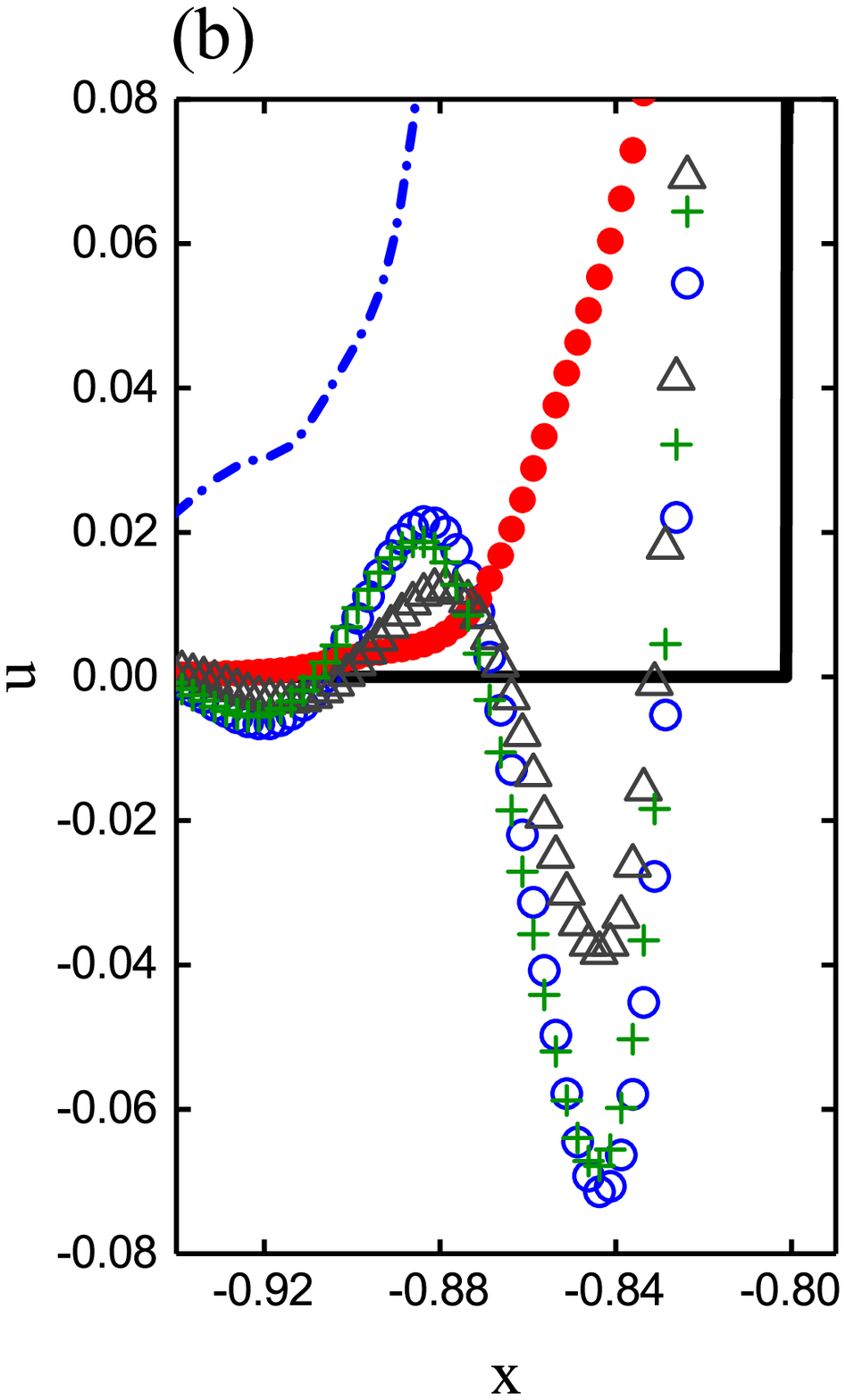}
\includegraphics[height=0.32\textwidth]
{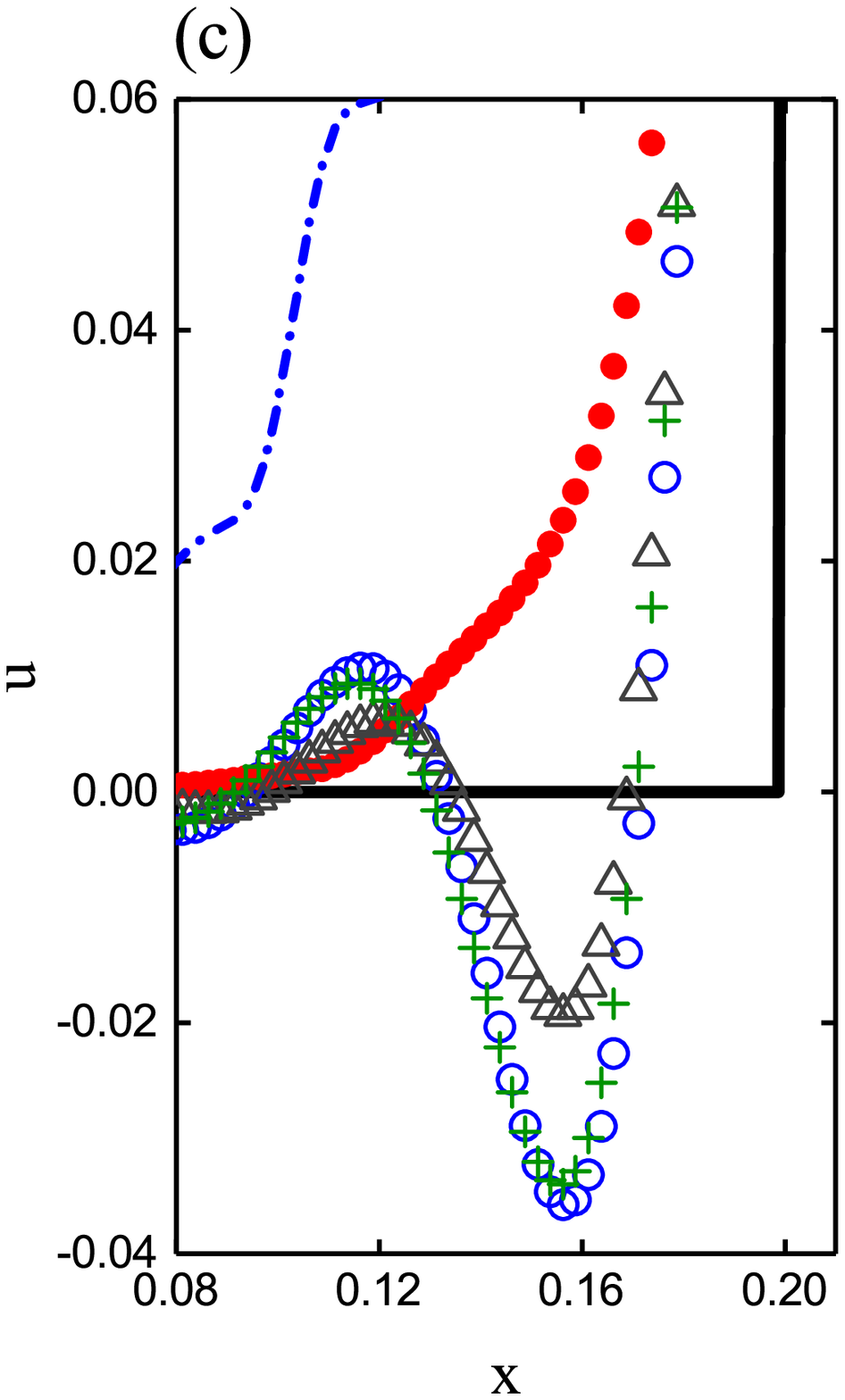}
\caption{Performance of the fifth-order MOP-WENO-ACM$k$, 
MIP-WENO-ACM$k$, WENO-JS, WENO-M, WENO-PM6 and WENO-IM($2,0.1$) 
schemes for the BiCWP with $N=800$ at long output time $t=2000$.}
\label{fig:BiCWP:800}
\end{figure}

We use the uniform mesh of $N = 800$ with the output time
$t=2000$ and use $N = 1600, 3200, 6400$ with the output time $t=200$ 
respectively to solve both SLP and BiCWP by all considered WENO 
schemes. Fig. \ref{fig:SLP:N800} and Fig. \ref{fig:BiCWP:800} show 
the comparison of various schemes when $t = 2000$ and $N=800$. We 
can observe that: (1) the MOP-WENO-ACM$k$ scheme provides the 
numerical results with a significantly higher resolution than those 
of the WENO-JS and WENO-M schemes, and it does not generate spurious 
oscillations while the WENO-PM6 and MIP-WENO-ACM$k$ schemes do, when 
solving both SLP and BiCWP; (2) when solving SLP on present 
computing condition, the WENO-IM(2, 0.1) scheme does not seem to 
generate spurious oscillations and it gives better resolution than 
the MOP-WENO-ACM$k$ scheme in most of the region; (3) however, from 
Fig. \ref{fig:SLP:N800}(b), we observe that the MOP-WENO-ACM$k$ 
scheme gives a better resolution of the Gaussian than the WENO-IM(2, 
0.1) scheme, and if we take a closer look, we can see that the 
WENO-IM(2, 0.1) scheme generate a very slight spurious oscillation 
near $x = - 0.435$ as shown in Fig. \ref{fig:SLP:N800}(c); (4) it is 
very evident as shown in Fig. \ref{fig:BiCWP:800} that, when solving 
BiCWP, the WENO-IM(2, 0.1) scheme generates the spurious 
oscillations. 

\begin{figure}[ht]
\centering
\includegraphics[height=0.32\textwidth]
{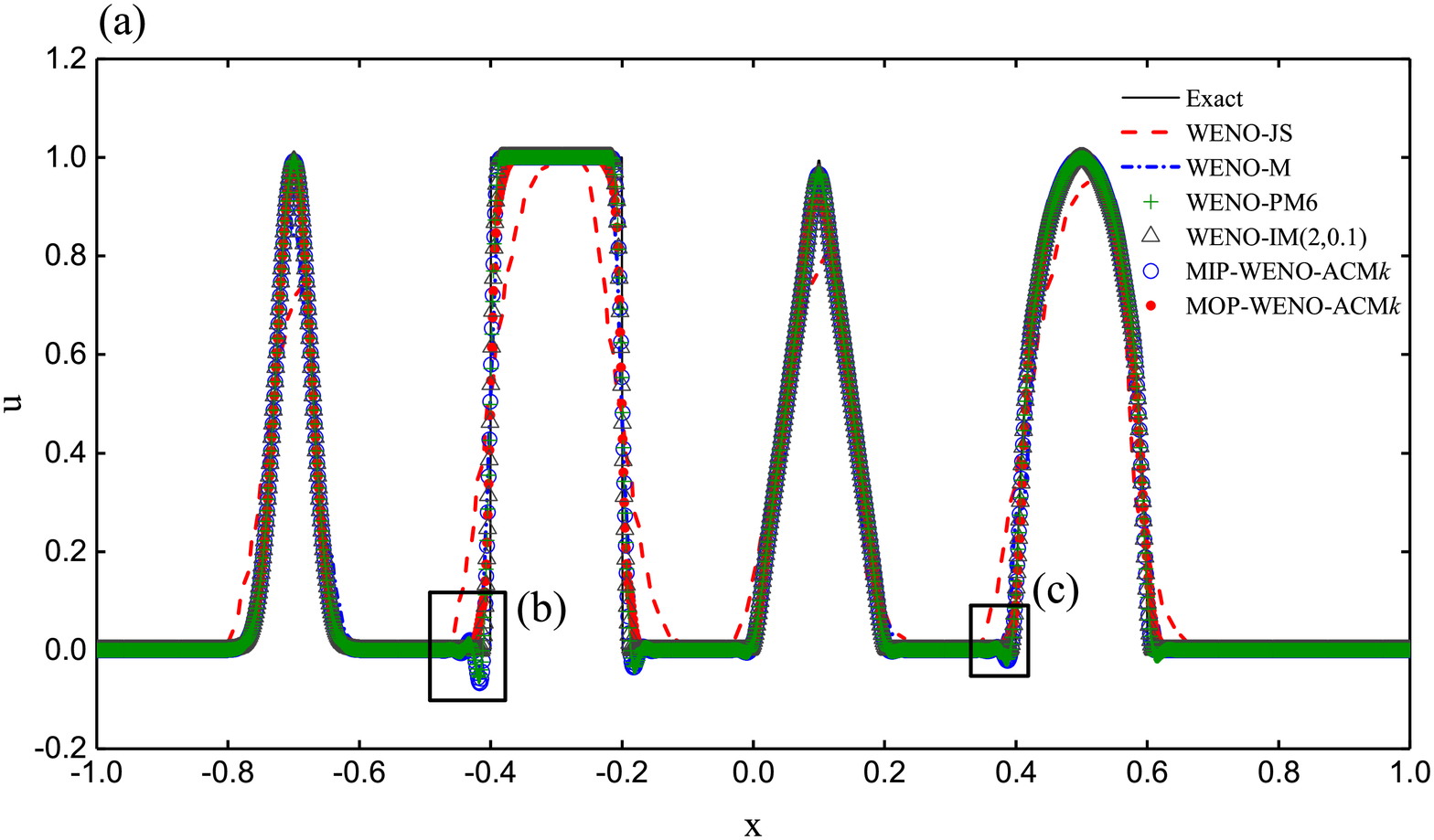}
\includegraphics[height=0.32\textwidth]
{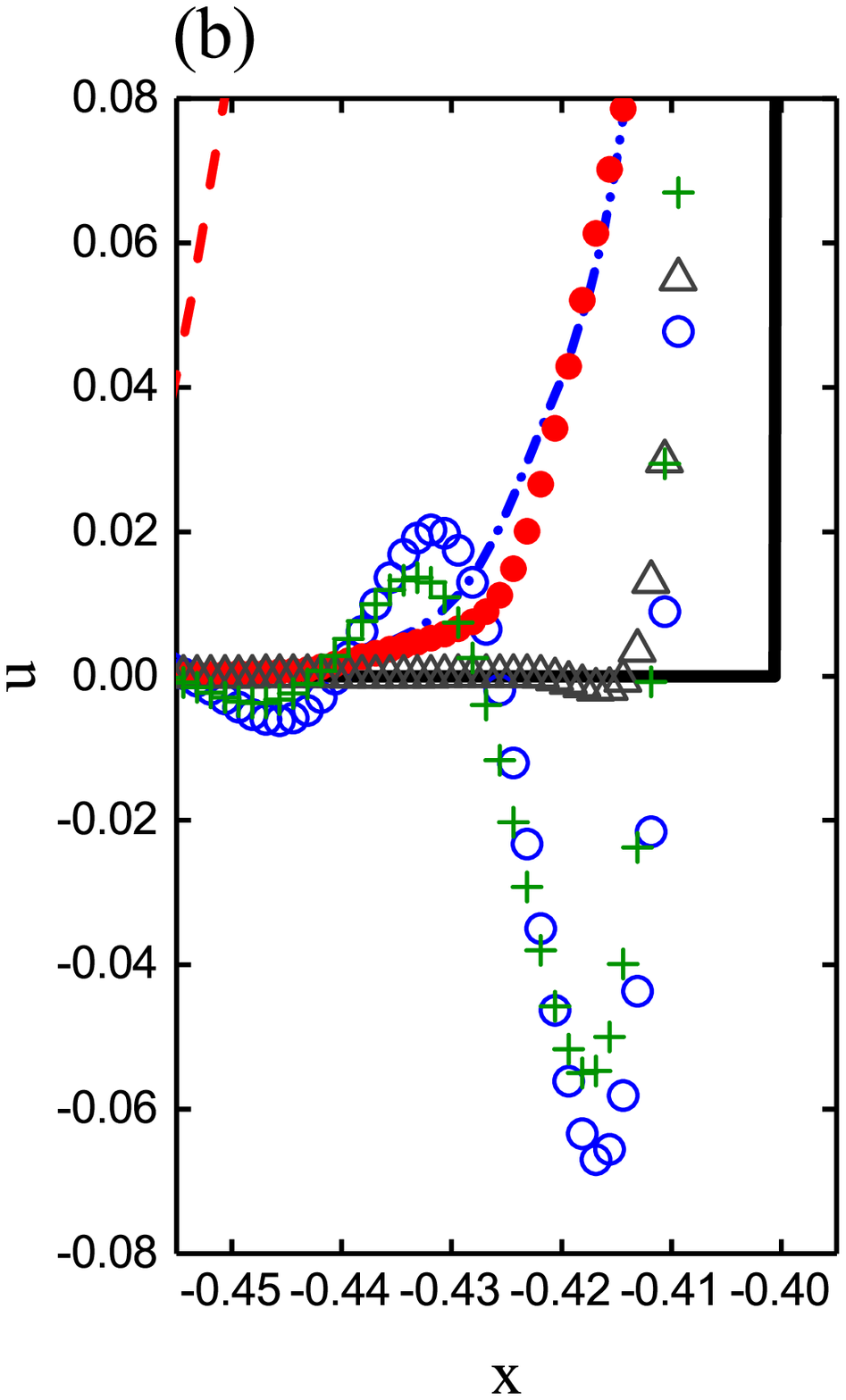}
\includegraphics[height=0.32\textwidth]
{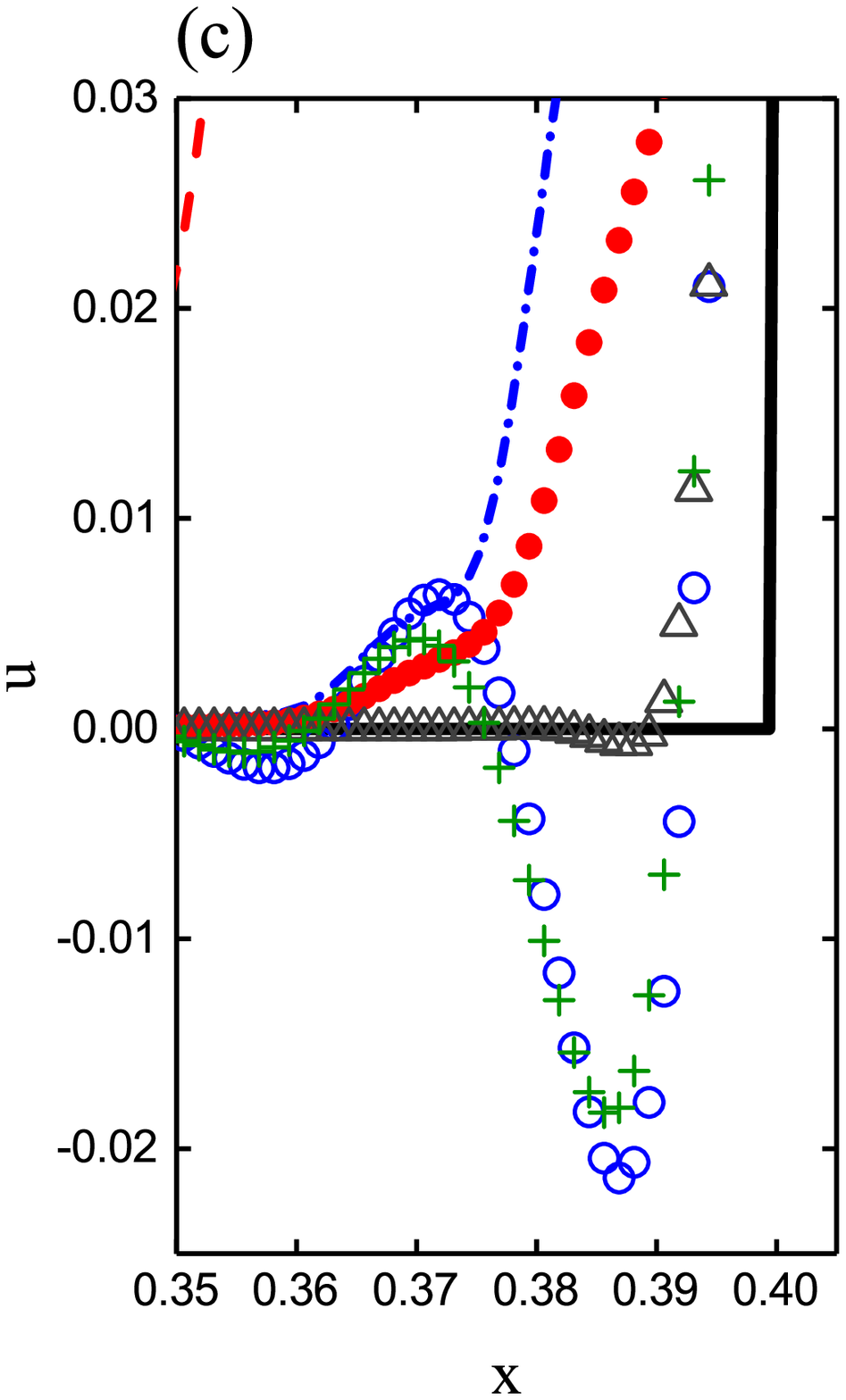}
\caption{Performance of the fifth-order MOP-WENO-ACM$k$, 
MIP-WENO-ACM$k$, WENO-JS, WENO-M, WENO-PM6 and WENO-IM($2,0.1$) 
schemes for the SLP with $N=1600$ at long output time $t=200$.}
\label{fig:SLP:1600}
\end{figure}

\begin{figure}[ht]
\centering
\includegraphics[height=0.32\textwidth]
{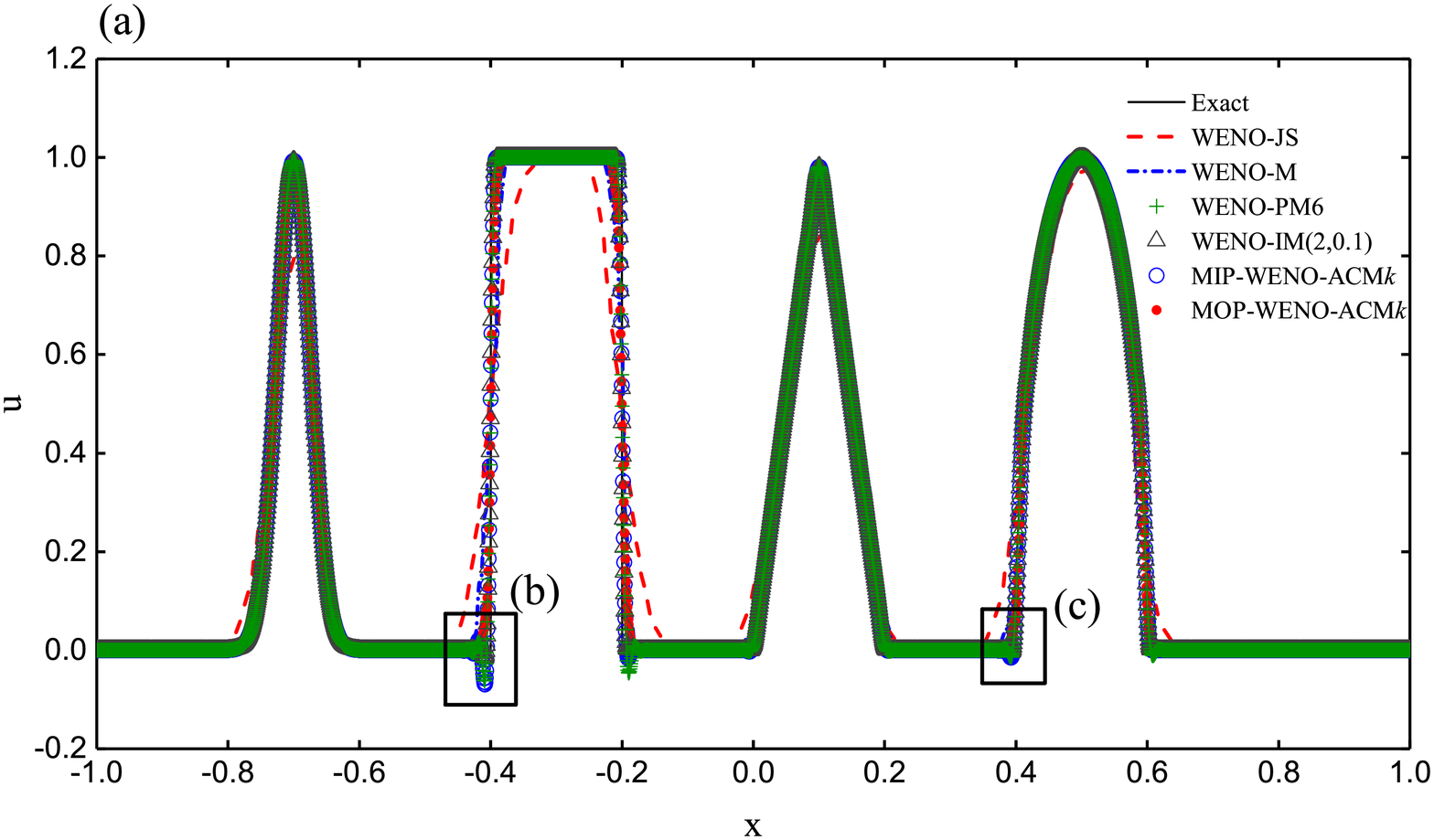}
\includegraphics[height=0.32\textwidth]
{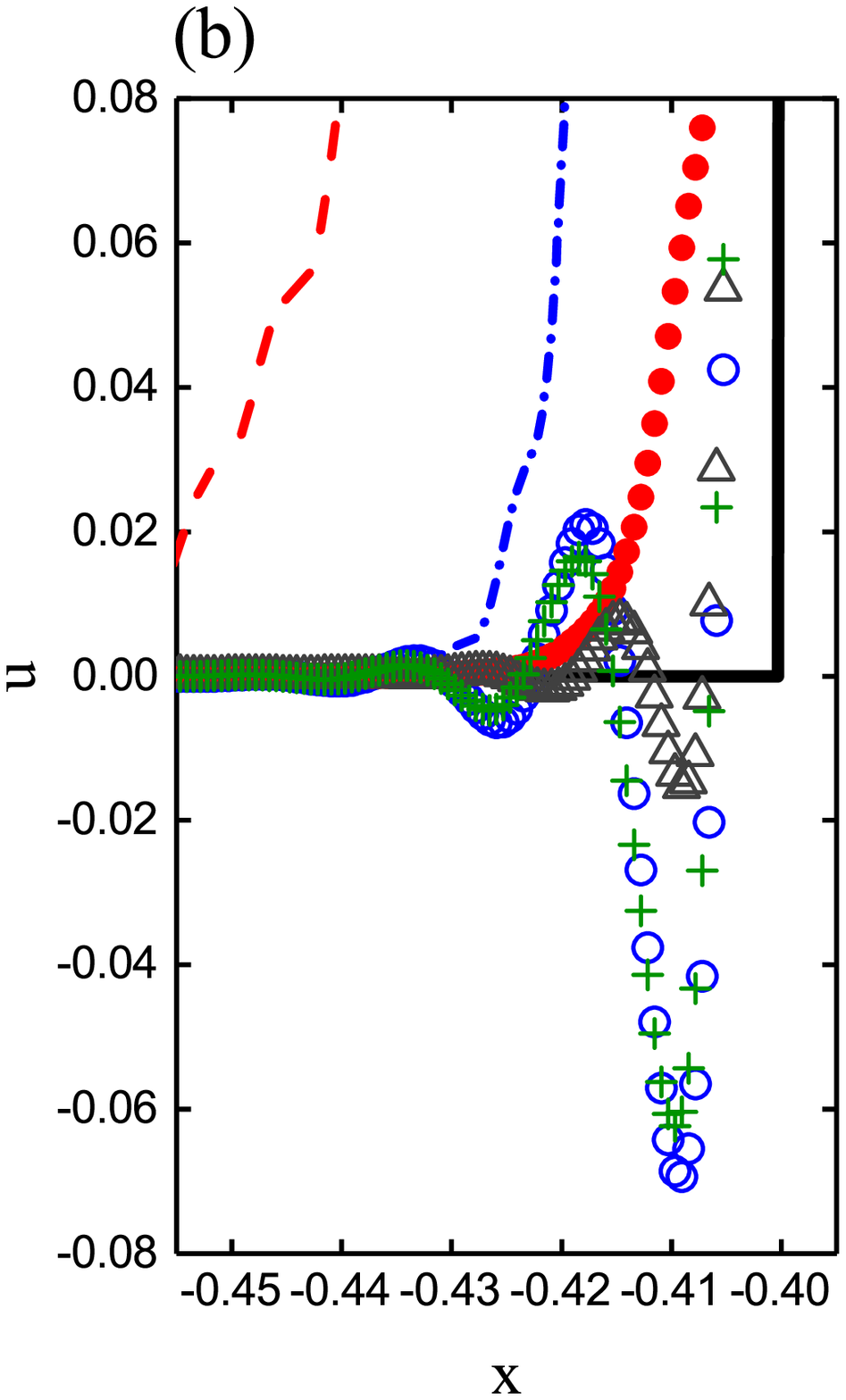}
\includegraphics[height=0.32\textwidth]
{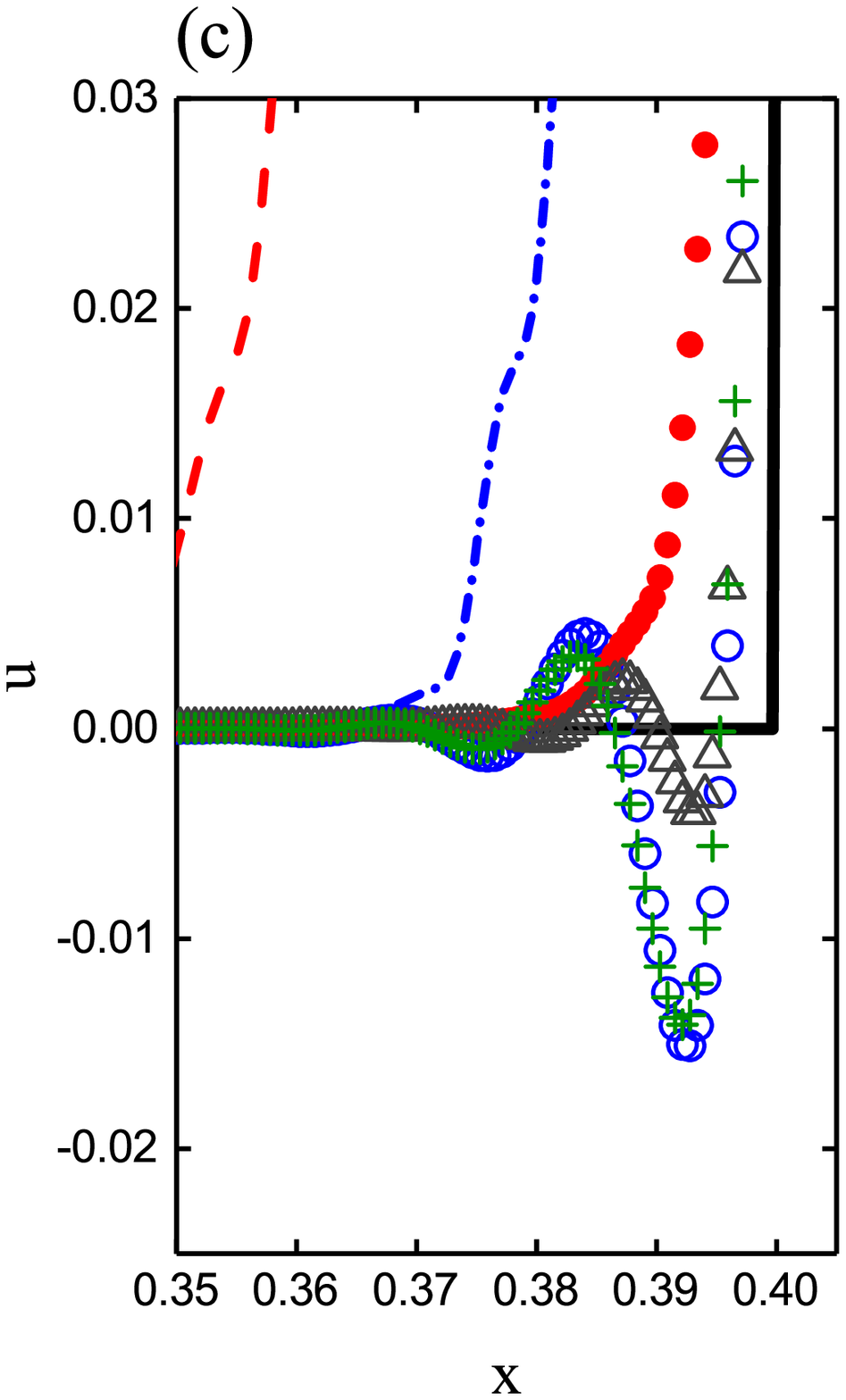}
\caption{Performance of the fifth-order MOP-WENO-ACM$k$, 
MIP-WENO-ACM$k$, WENO-JS, WENO-M, WENO-PM6 and WENO-IM($2,0.1$) 
schemes for the SLP with $N=3200$ at long output time $t=200$.}
\label{fig:SLP:3200}
\end{figure}

\begin{figure}[ht]
\centering
\includegraphics[height=0.32\textwidth]
{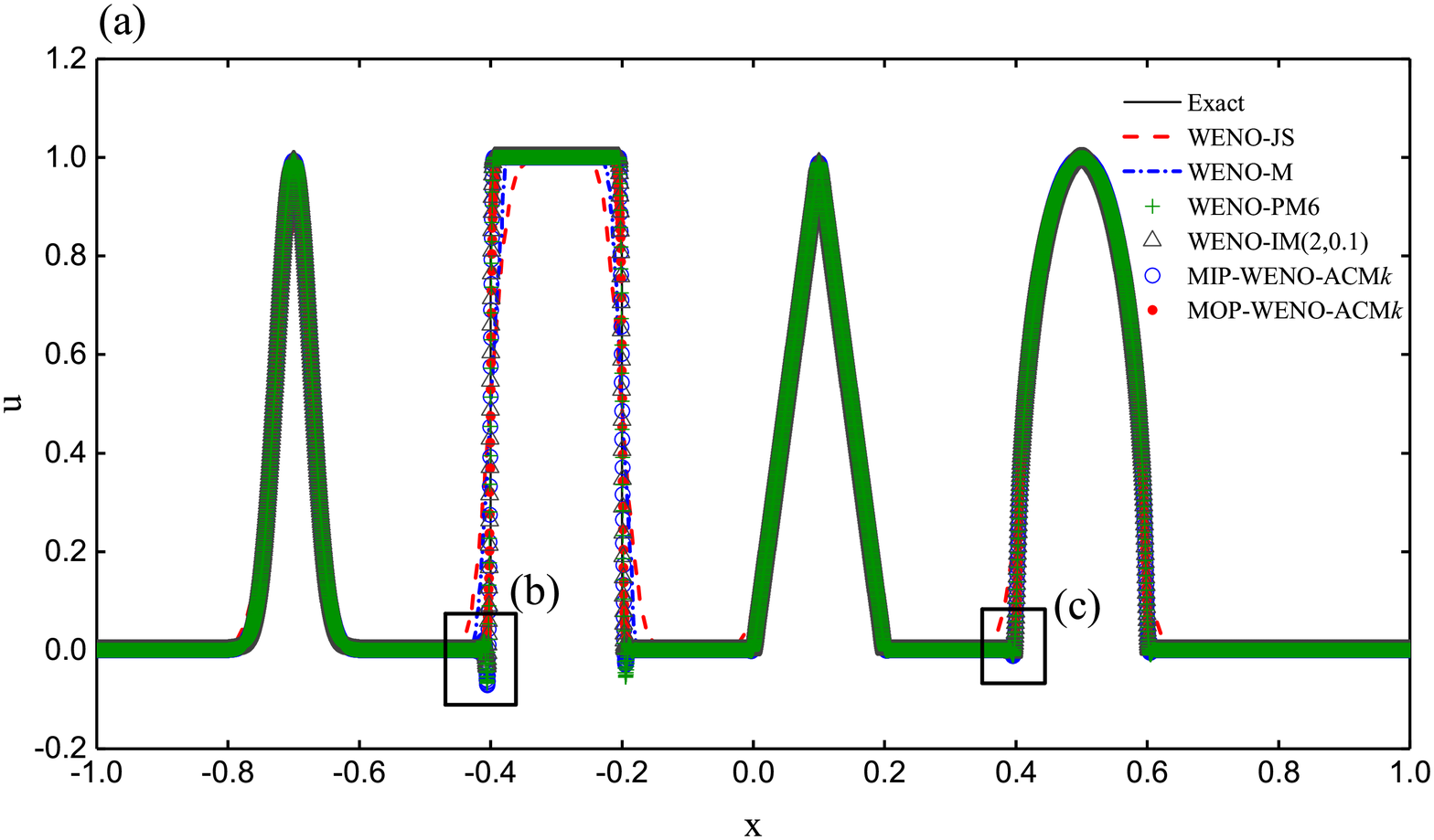}
\includegraphics[height=0.32\textwidth]
{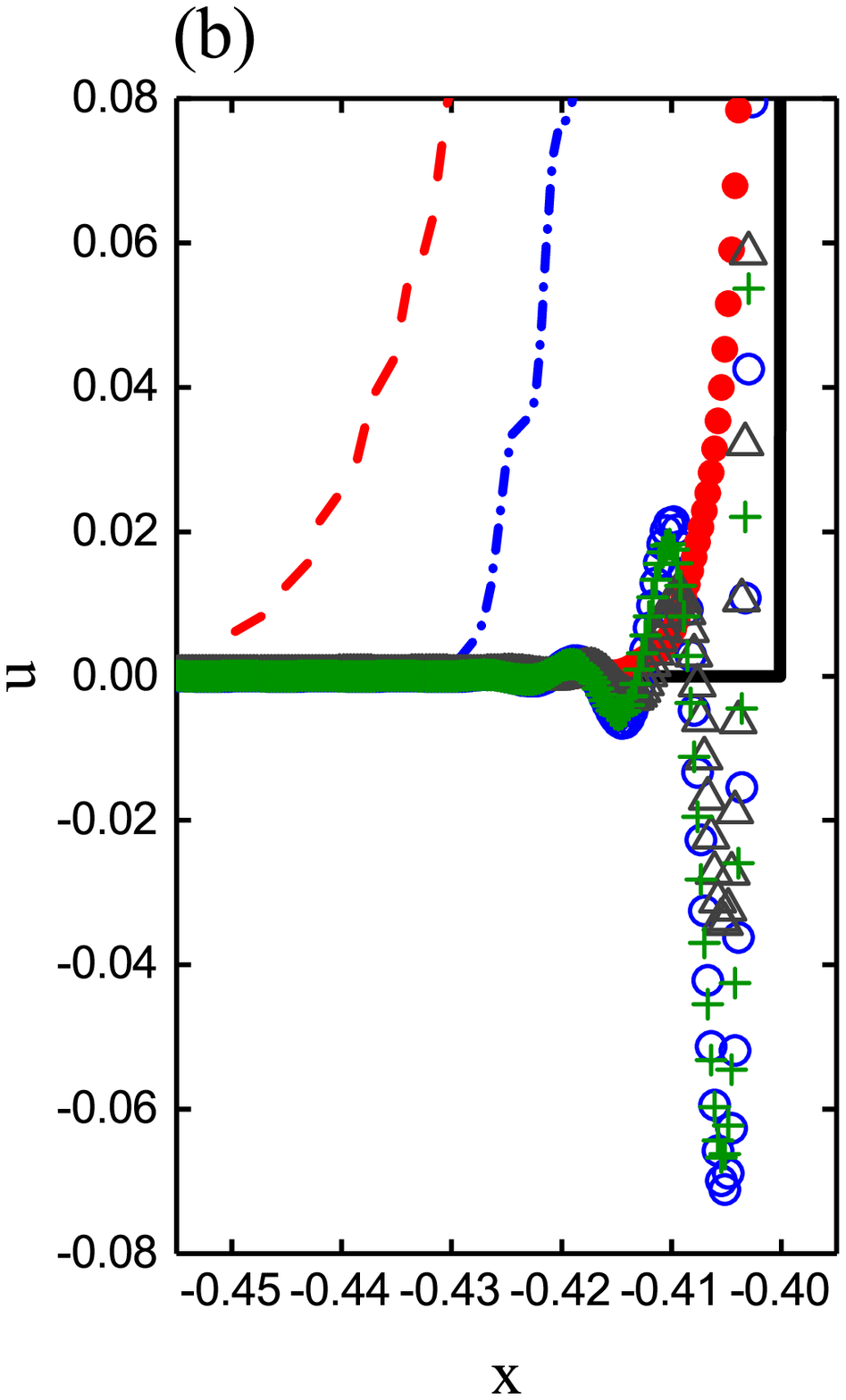}
\includegraphics[height=0.32\textwidth]
{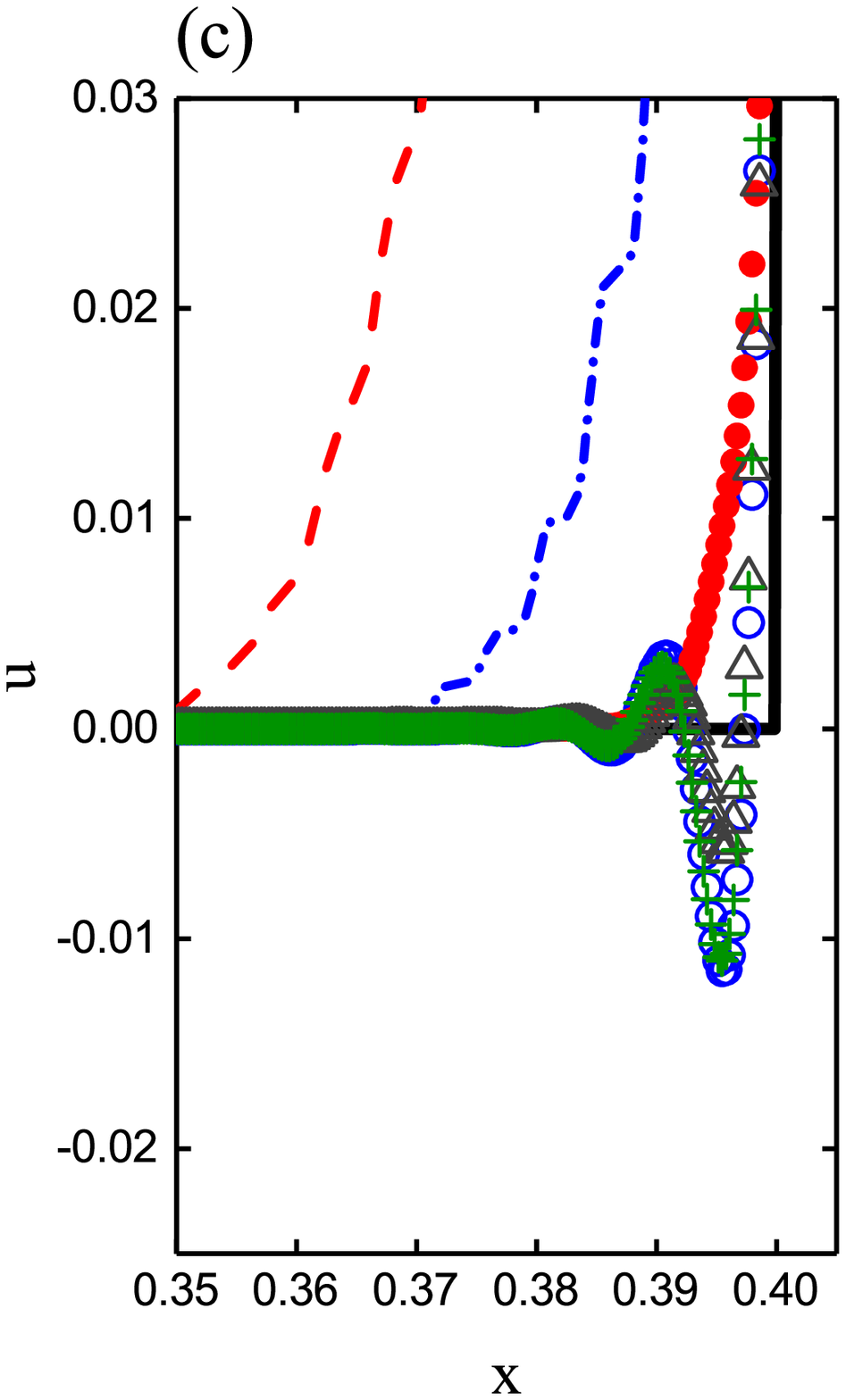}
\caption{Performance of the fifth-order MOP-WENO-ACM$k$, 
MIP-WENO-ACM$k$, WENO-JS, WENO-M, WENO-PM6 and WENO-IM($2,0.1$) 
schemes for the SLP with $N=6400$ at long output time $t=200$.}
\label{fig:SLP:6400}
\end{figure}

\begin{figure}[ht]
\centering
\includegraphics[height=0.32\textwidth]
{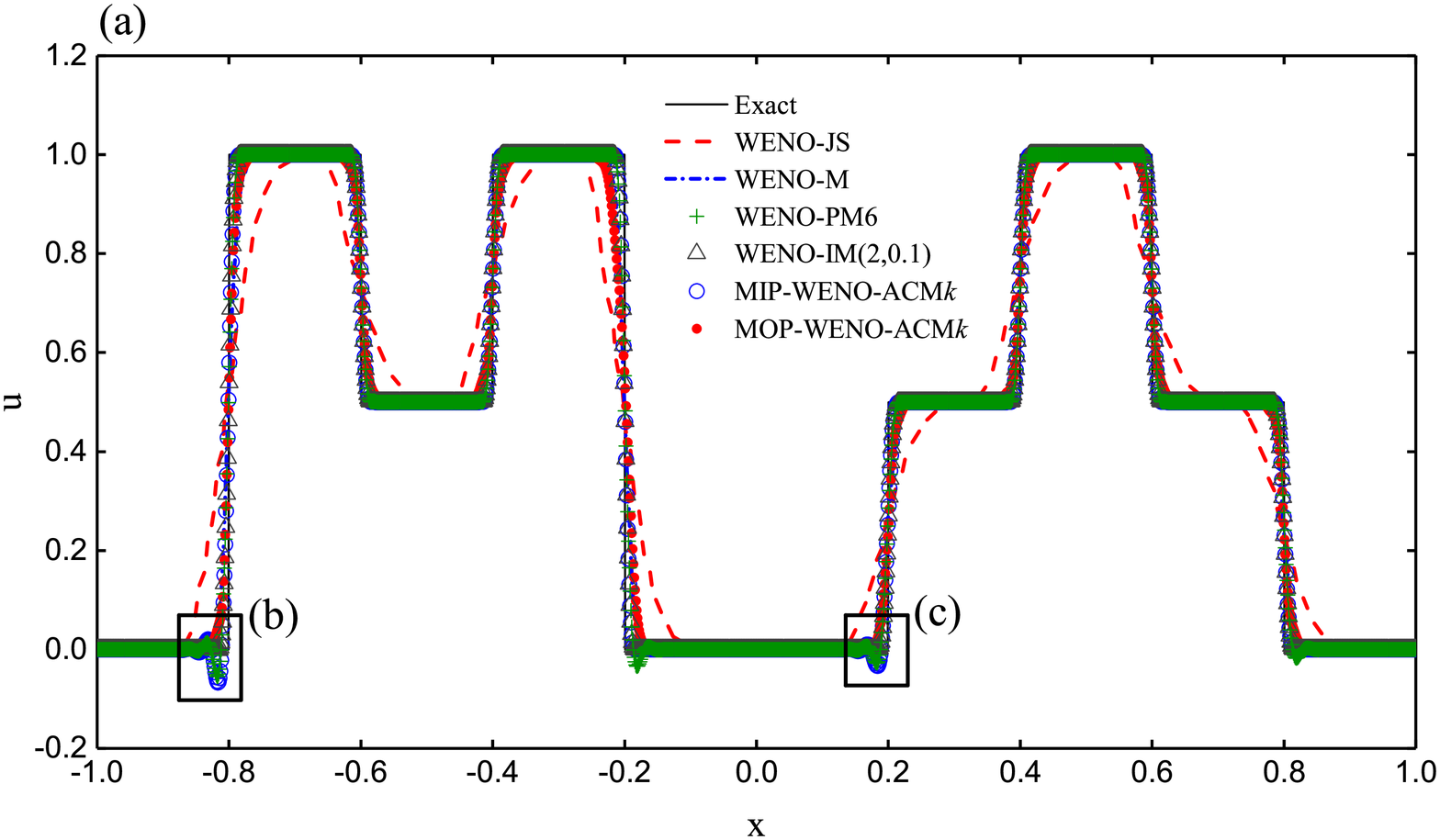}
\includegraphics[height=0.32\textwidth]
{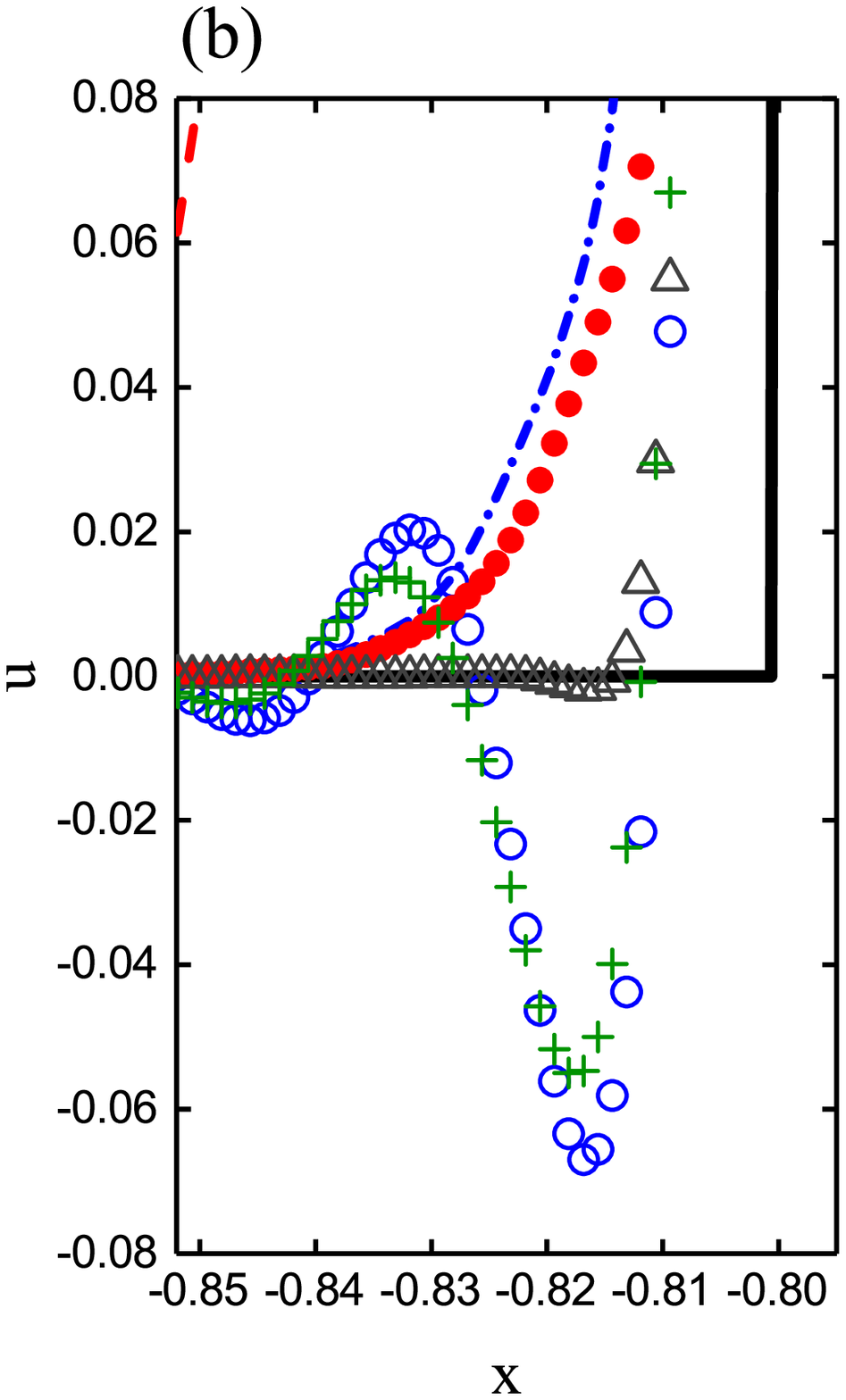}
\includegraphics[height=0.32\textwidth]
{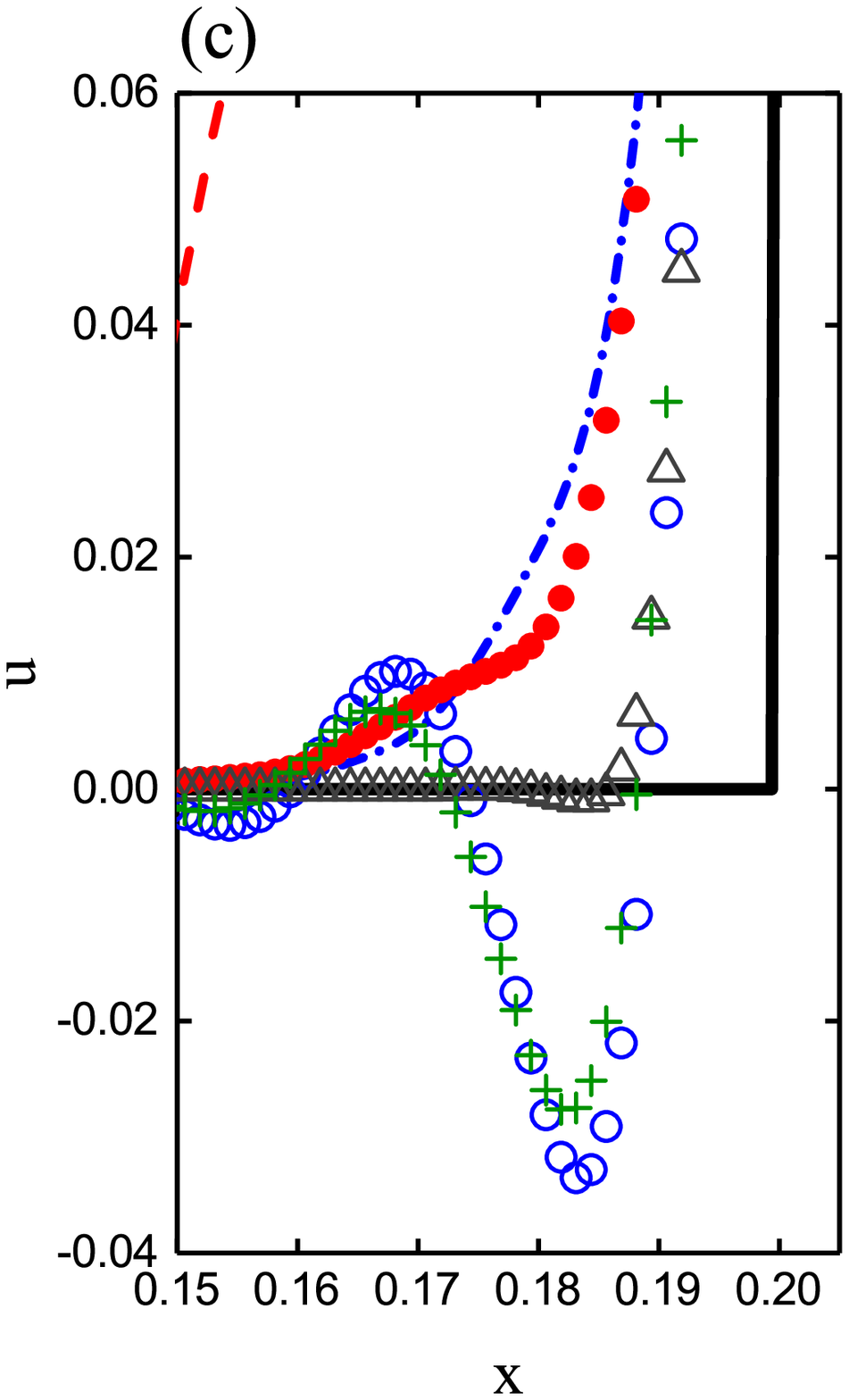}
\caption{Performance of the fifth-order MOP-WENO-ACM$k$, 
MIP-WENO-ACM$k$, WENO-JS, WENO-M, WENO-PM6 and WENO-IM($2,0.1$) 
schemes for the BiCWP with $N=1600$ at long output time $t=200$.}
\label{fig:BiCWP:1600}
\end{figure}

\begin{figure}[ht]
\centering
\includegraphics[height=0.32\textwidth]
{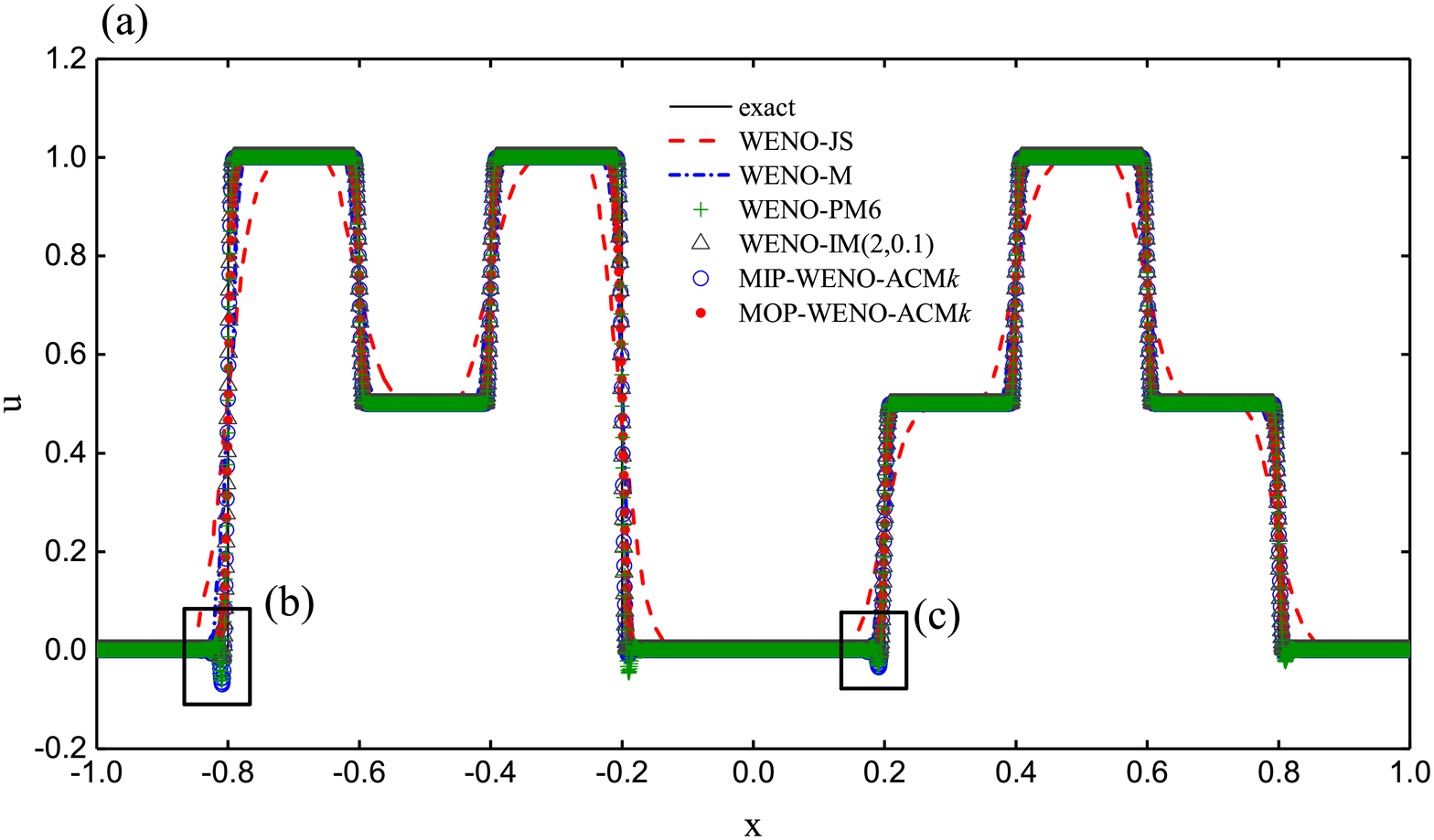}
\includegraphics[height=0.32\textwidth]
{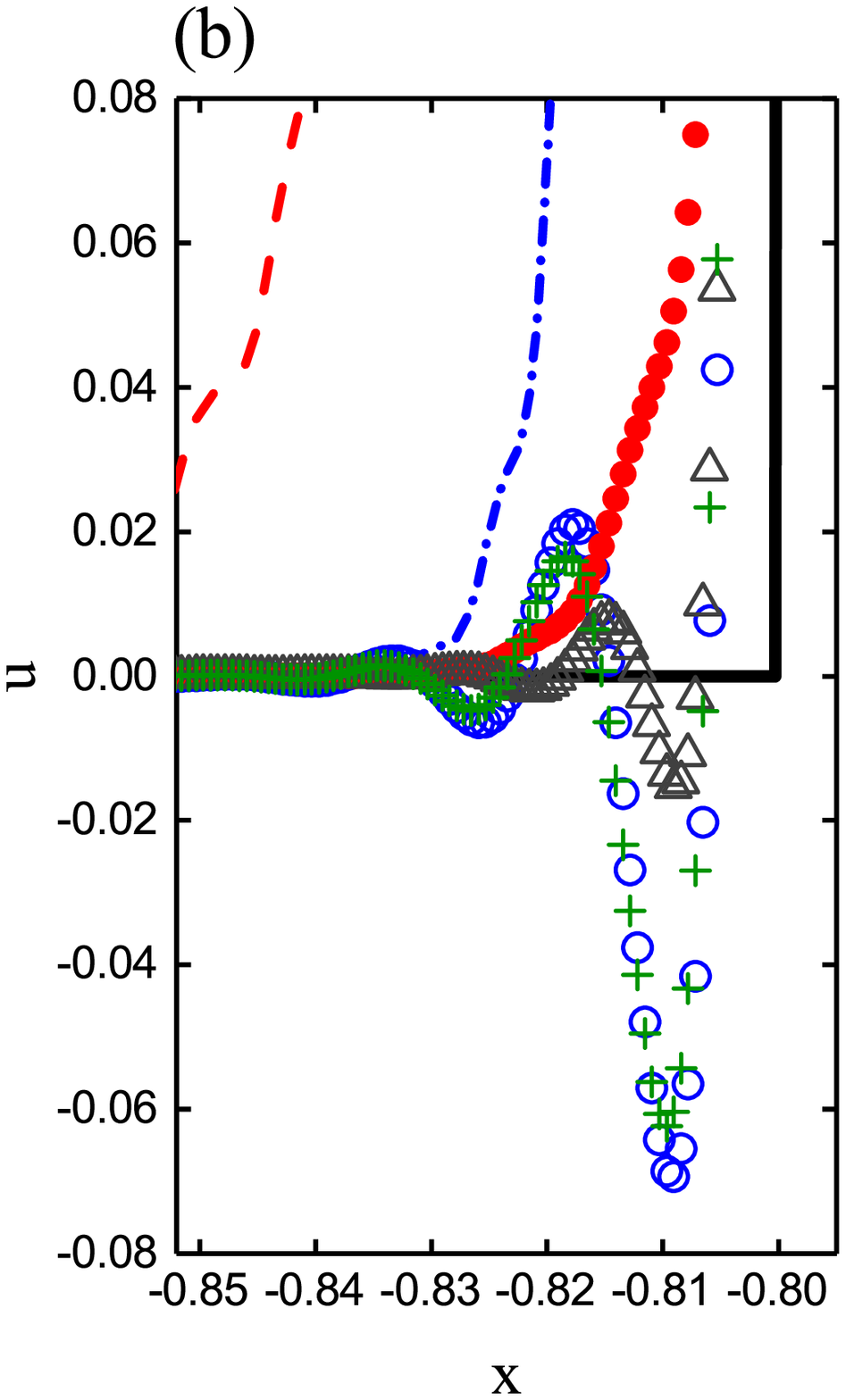}
\includegraphics[height=0.32\textwidth]
{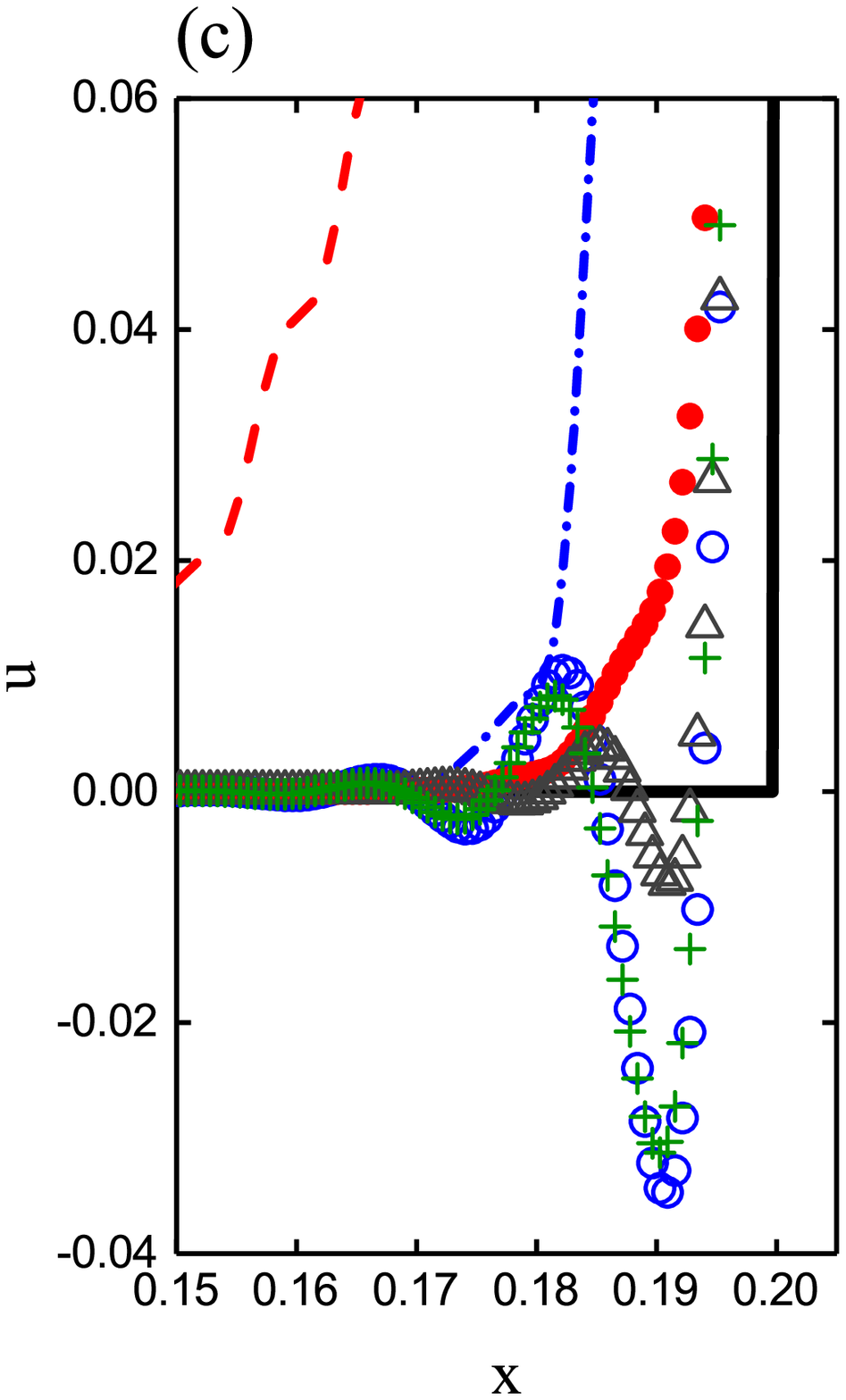}
\caption{Performance of the fifth-order MOP-WENO-ACM$k$, 
MIP-WENO-ACM$k$, WENO-JS, WENO-M, WENO-PM6 and WENO-IM($2,0.1$) 
schemes for the BiCWP with $N=3200$ at long output time $t=200$.}
\label{fig:BiCWP:3200}
\end{figure}

\begin{figure}[ht]
\centering
\includegraphics[height=0.32\textwidth]
{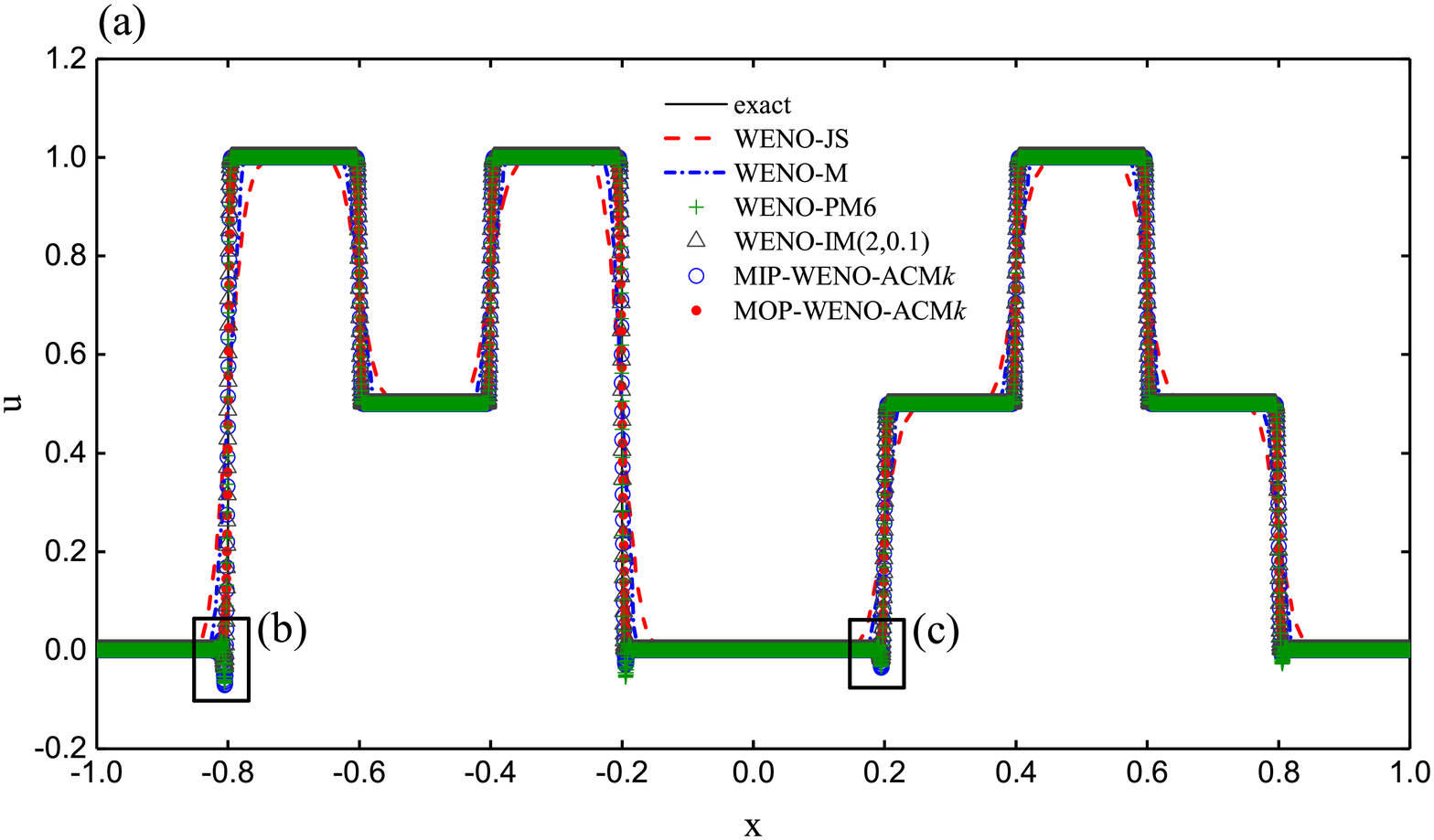}
\includegraphics[height=0.32\textwidth]
{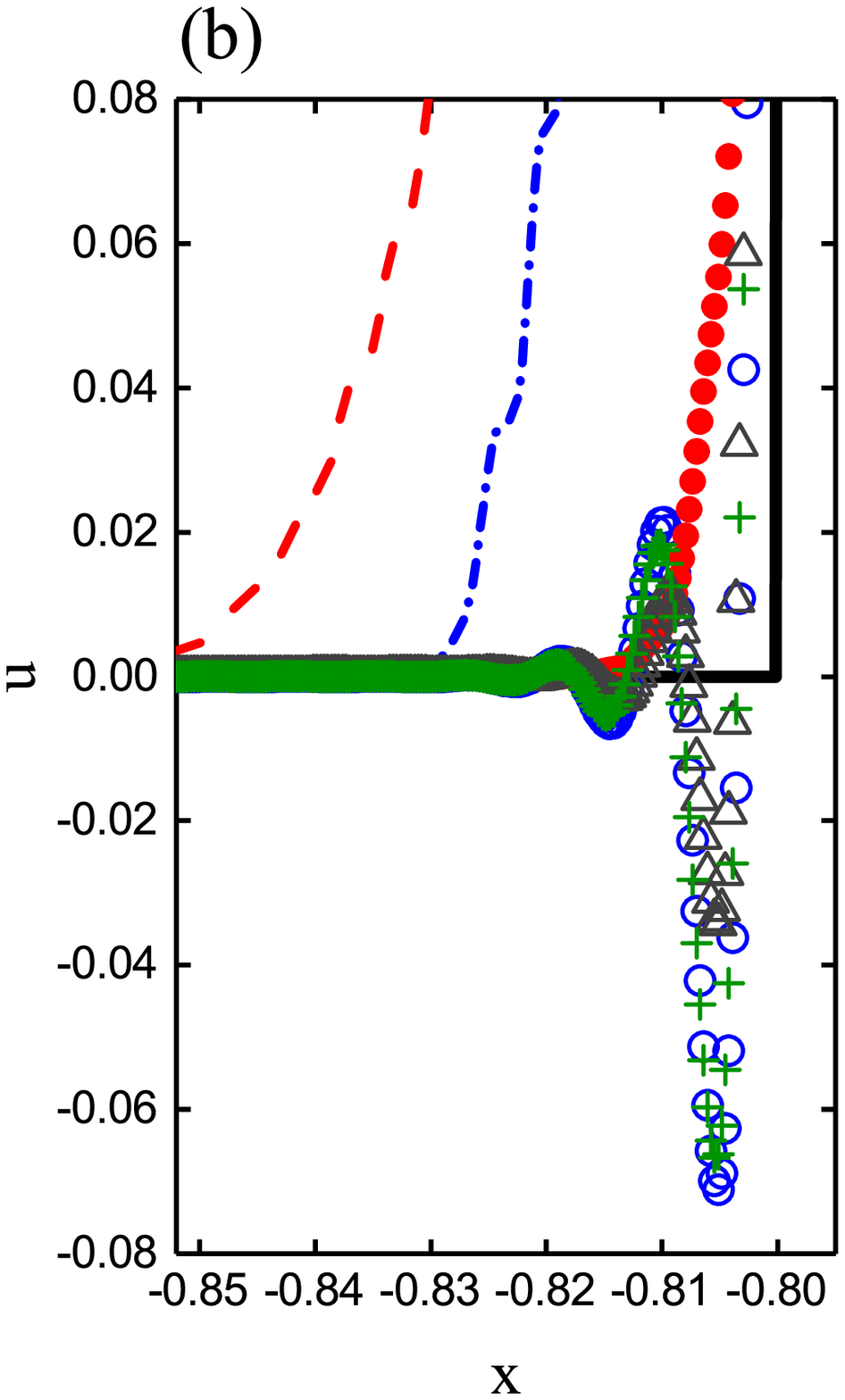}
\includegraphics[height=0.32\textwidth]
{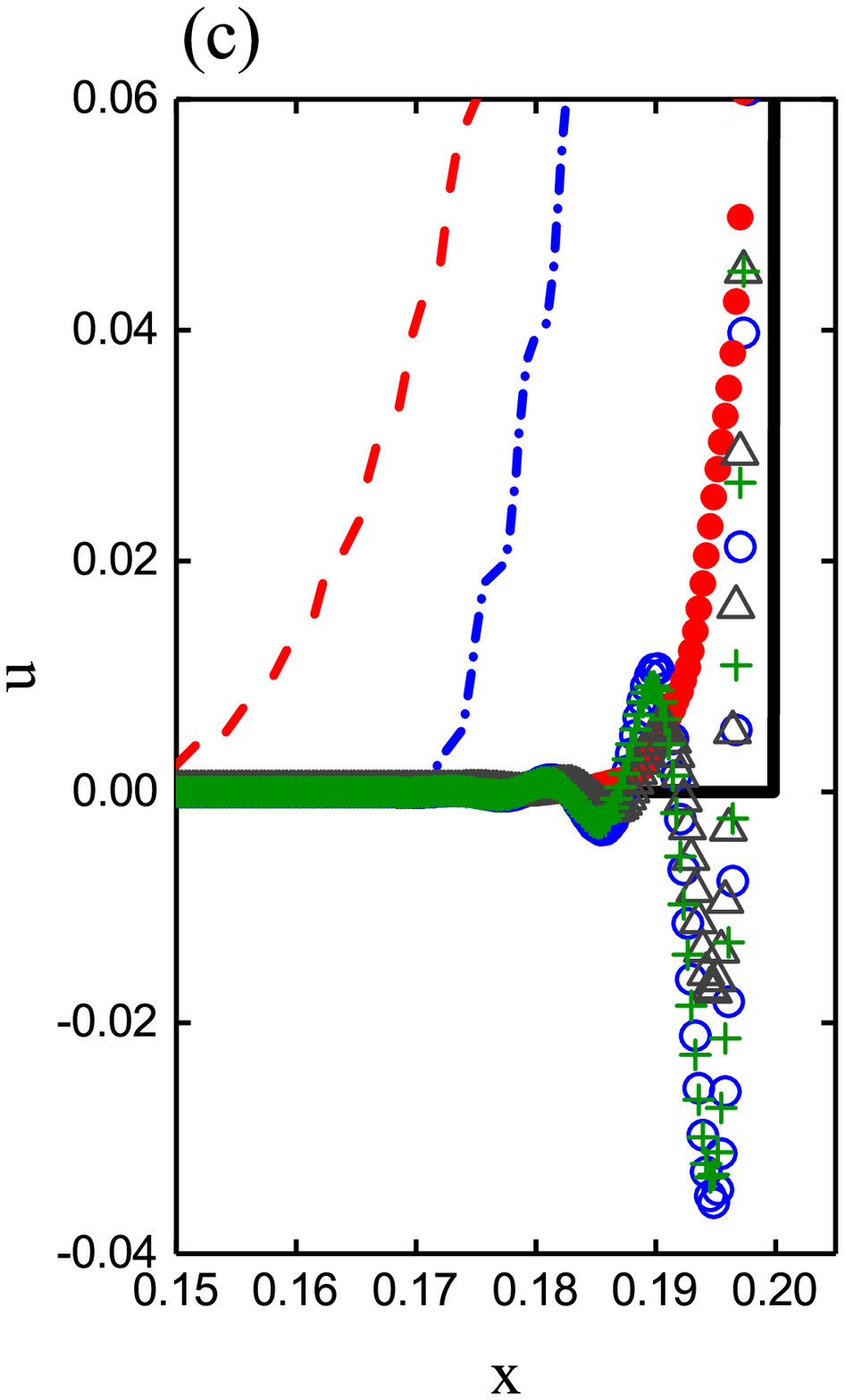}
\caption{Performance of the fifth-order MOP-WENO-ACM$k$, 
MIP-WENO-ACM$k$, WENO-JS, WENO-M, WENO-PM6 and WENO-IM($2,0.1$) 
schemes for the BiCWP with $N=6400$ at long output time $t=200$.}
\label{fig:BiCWP:6400}
\end{figure}

In Figs. \ref{fig:SLP:1600}, \ref{fig:SLP:3200}, \ref{fig:SLP:6400} 
and Figs. \ref{fig:BiCWP:1600}, \ref{fig:BiCWP:3200}, 
\ref{fig:BiCWP:6400}, we show the comparison of various schemes when 
$t=200$ and $N=1600, 3200, 6400$ for SLP and BiCWP, respectively. 
From these solutions computed with larger grid numbers and a 
shortened but still long output time, we observe that, as the grid 
number increases: (1) firstly, the WENO-IM(2, 0.1) scheme generates 
spurious oscillations but the MOP-WENO-ACM$k$ scheme does not while 
provides an improved resolution when solving SLP; (2) although the 
resolutions of the results computed by the WENO-JS and WENO-M 
schemes are significantly improved for both SLP and BiCWP, the 
MOP-WENO-ACM$k$ scheme still evidently provides better resolutions 
than those of these two schemes; (3) the spurious oscillations 
generated by the WENO-PM6, WENO-IM(2,0.1) and MIP-WENO-ACM$k$ 
schemes appear to be more evident and more intense when the grid 
number gets larger, while the MOP-WENO-ACM$k$ scheme can still 
prevent generating spurious oscillations but obtain great 
improvement of the resolution, when solving both SLP and BiCWP.

\begin{figure}[ht]
\centering
\includegraphics[height=0.305\textwidth]
{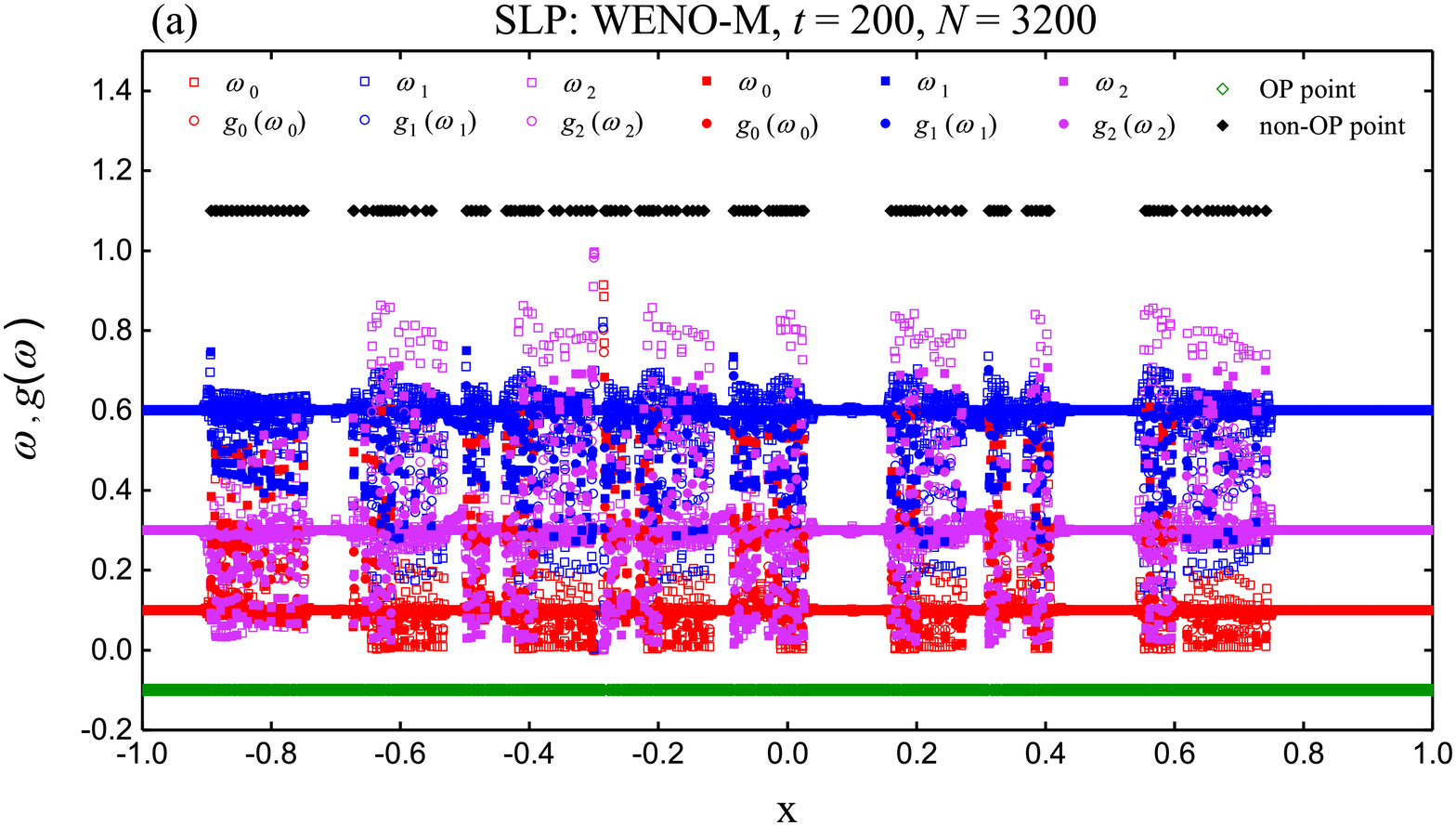}
\includegraphics[height=0.305\textwidth]
{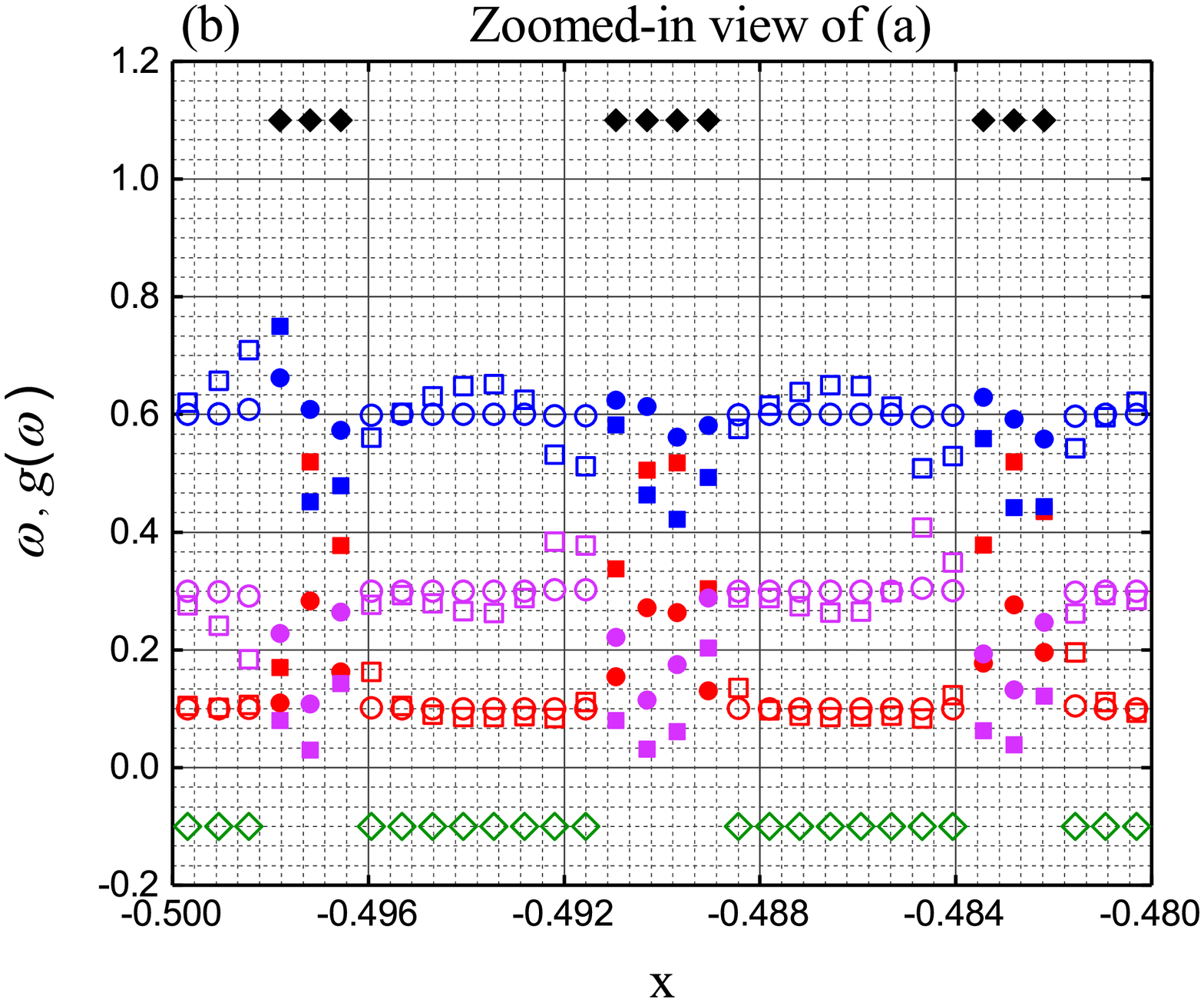}\\
\includegraphics[height=0.305\textwidth]
{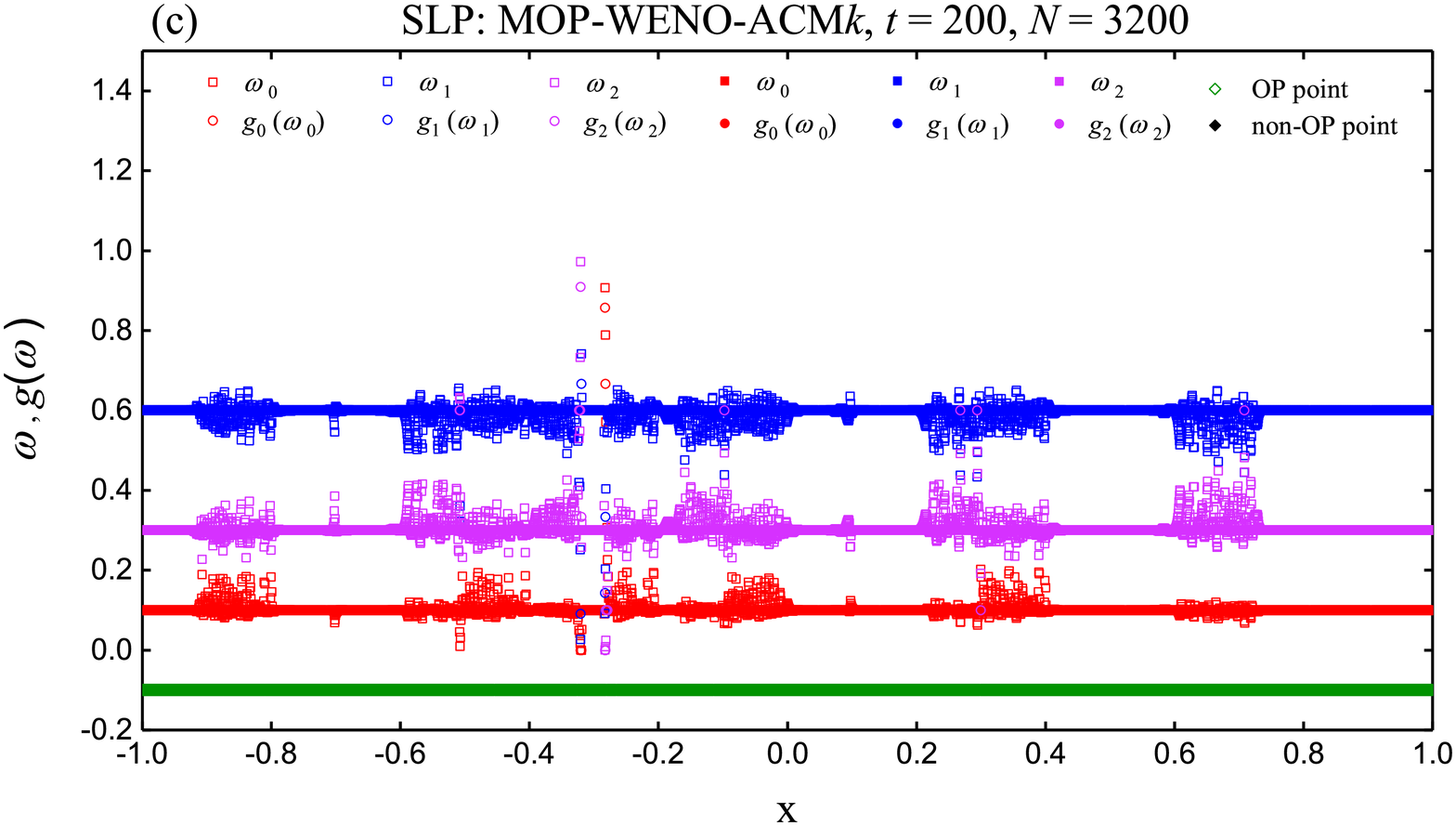}
\includegraphics[height=0.305\textwidth]
{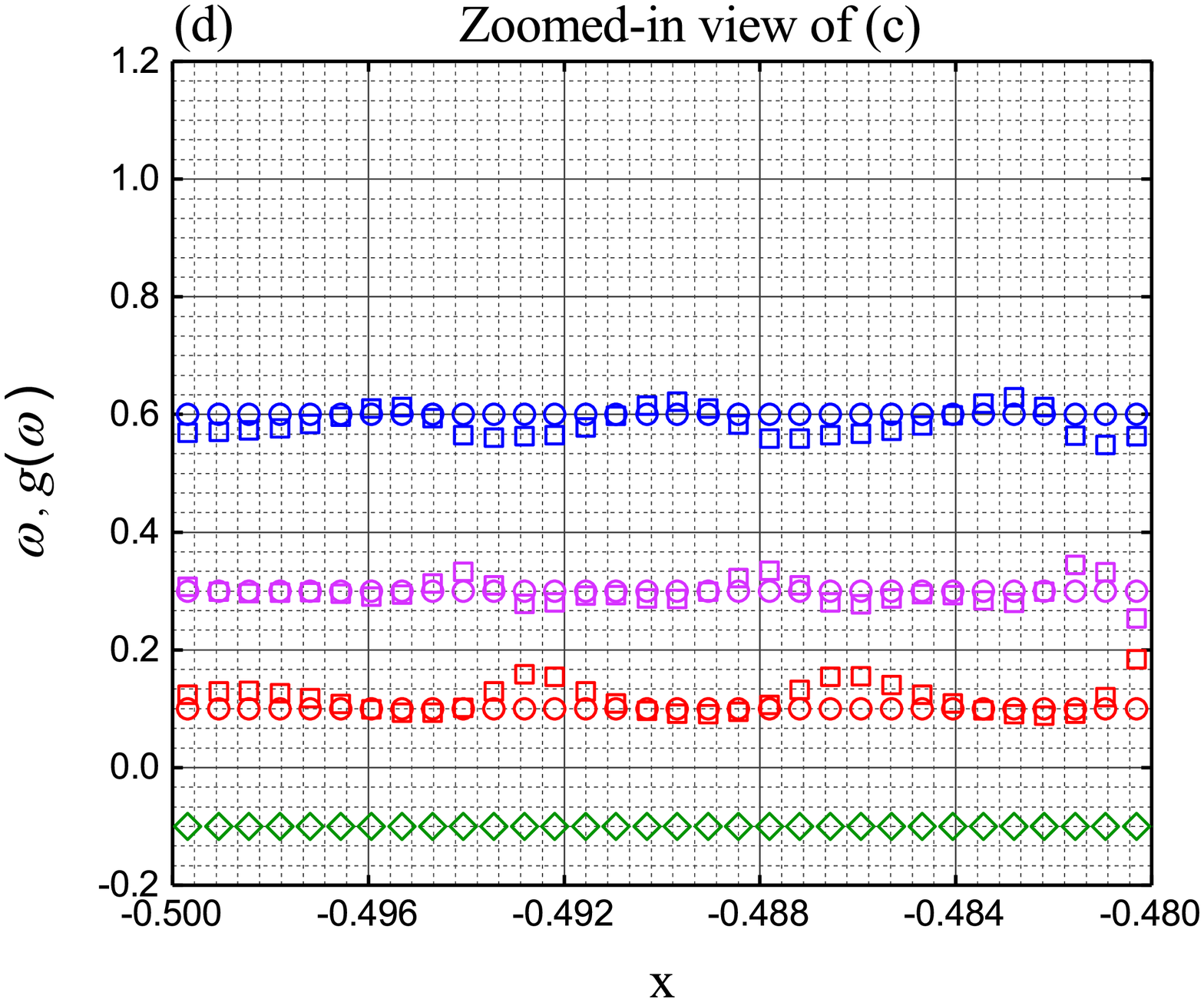}\\
\includegraphics[height=0.305\textwidth]
{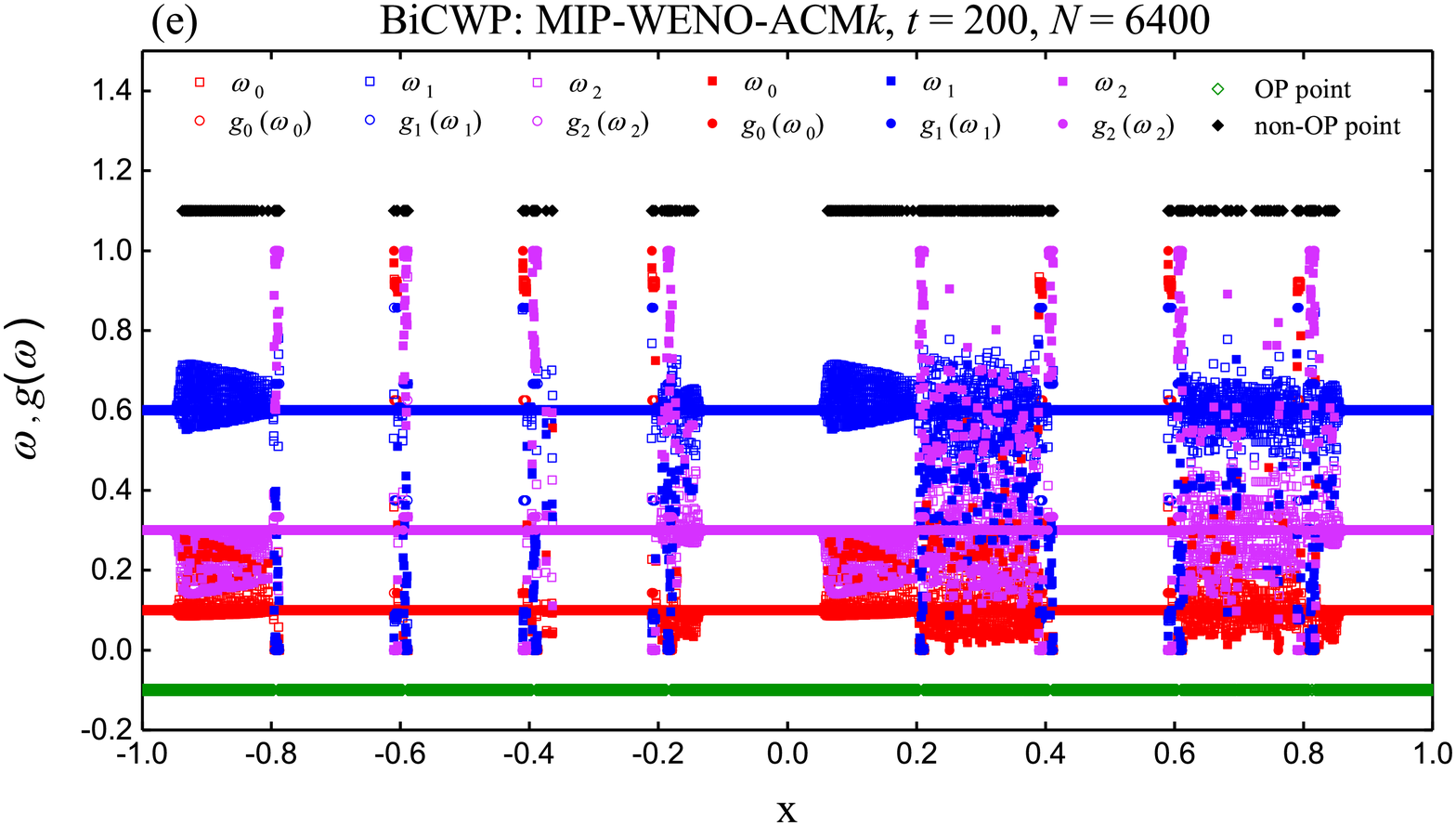}
\includegraphics[height=0.305\textwidth]
{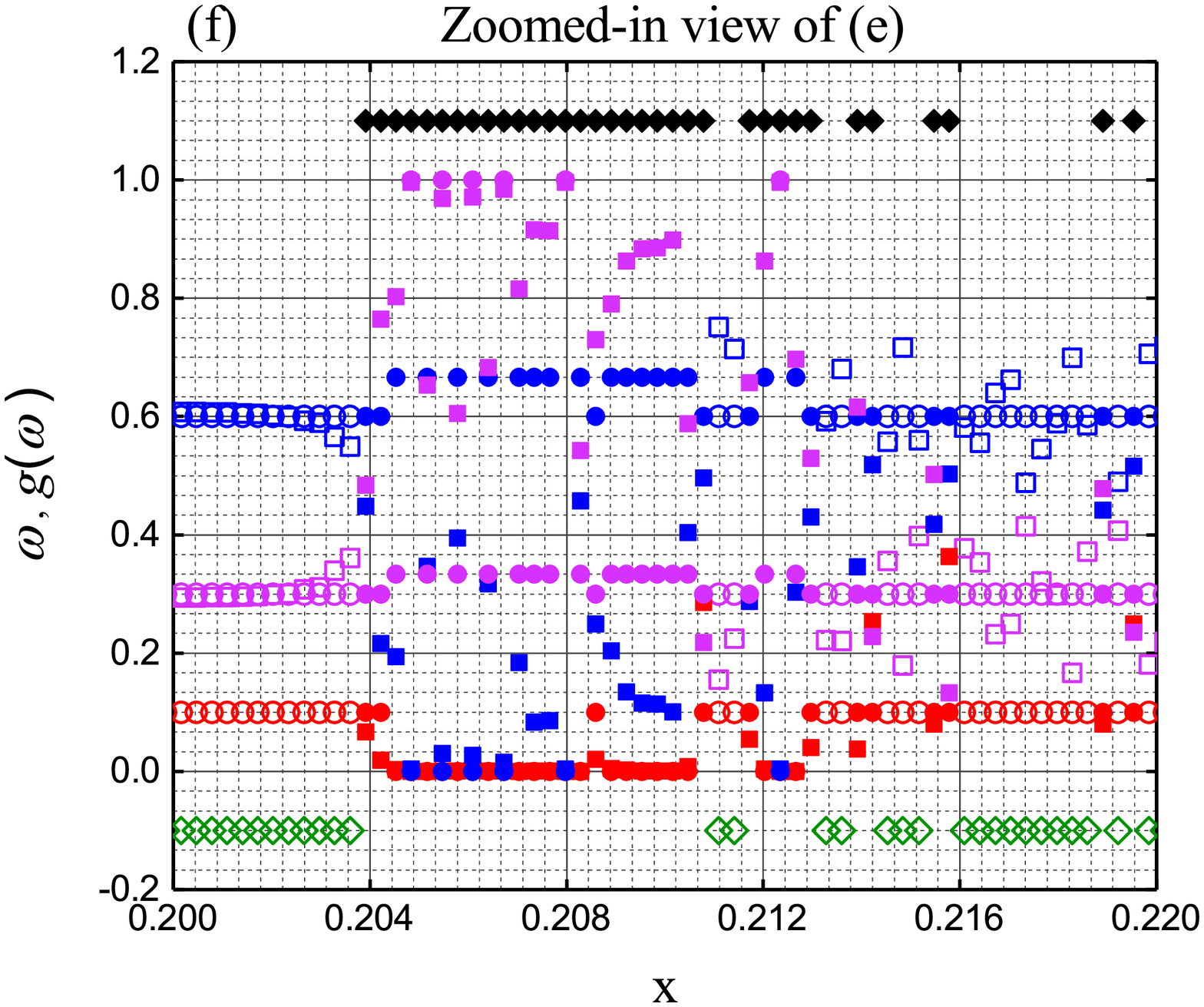}\\
\includegraphics[height=0.305\textwidth]
{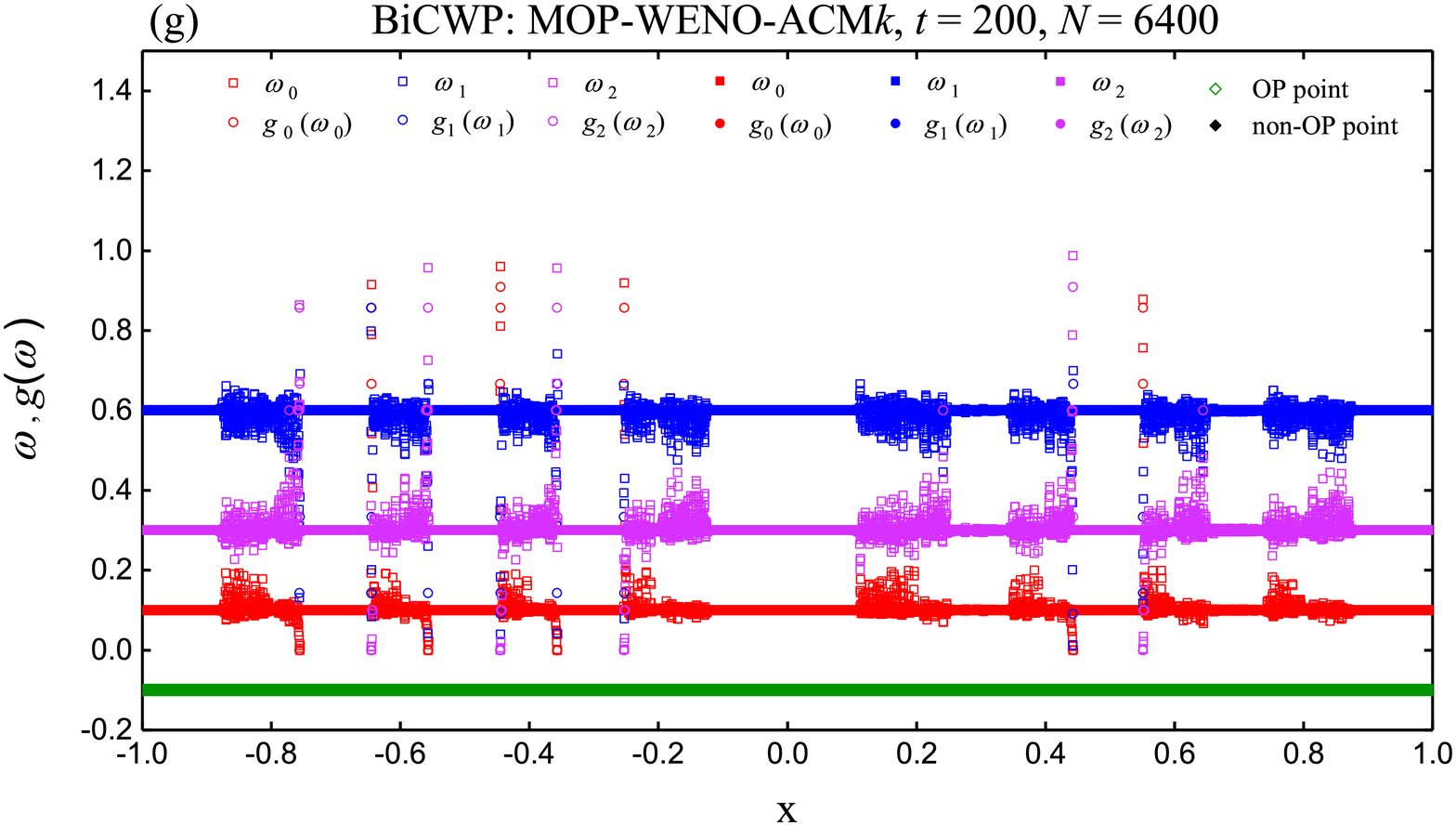}
\includegraphics[height=0.305\textwidth]
{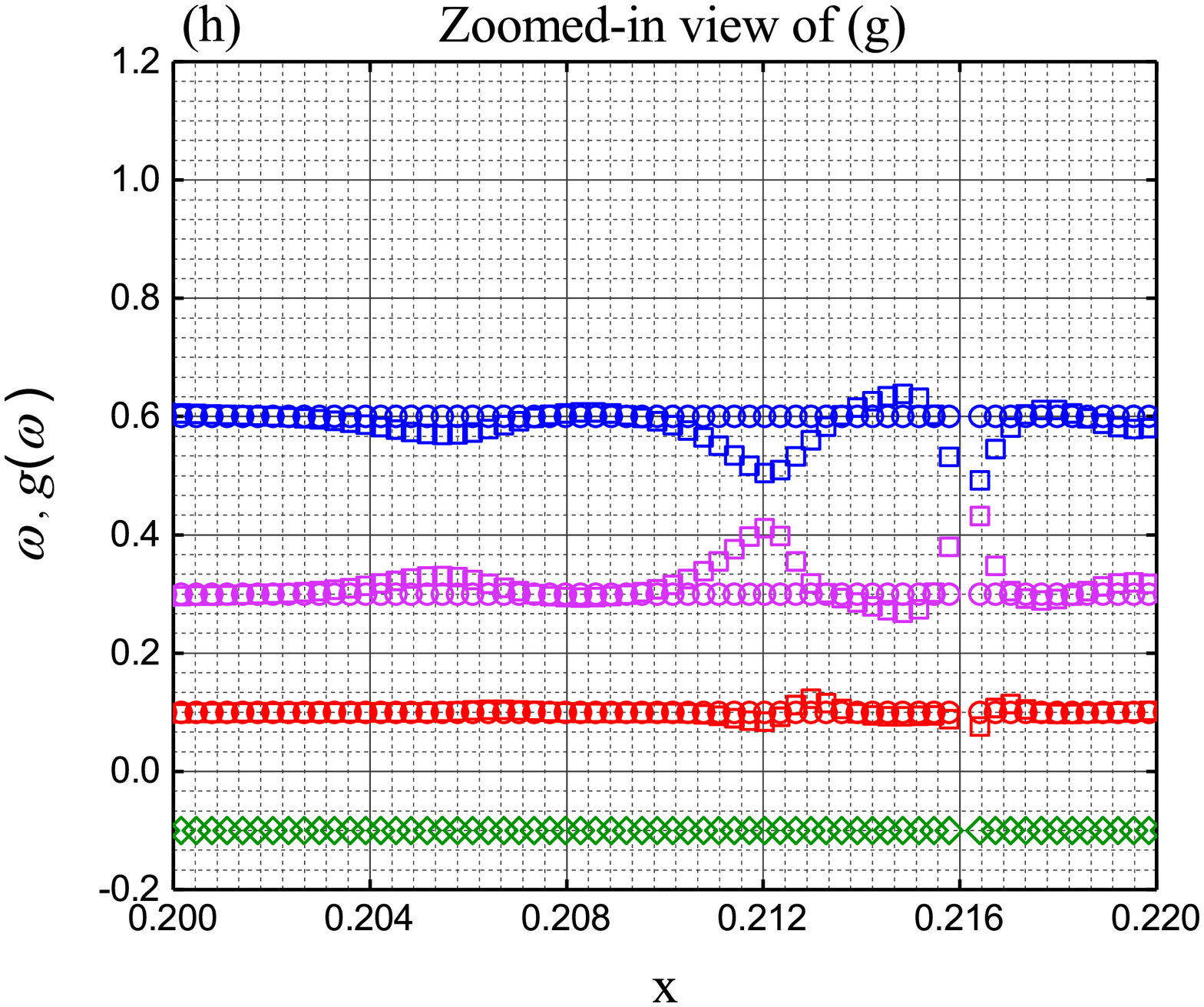}
\caption{The \textit{non-OP points} in the numerical solutions of 
SLP computed by the WENO-M and MOP-WENO-ACM$k$ schemes with $N=3200,
t=200$, and the \textit{non-OP points} in the numerical solutions of 
BiCWP computed by the MIP-WENO-ACM$k$ and MOP-WENO-ACM$k$ schemes 
with $N=6400, t= 200$.}
\label{fig:x-Omega:SLP_BiCWP}
\end{figure}

As examples, in Fig. \ref{fig:x-Omega:SLP_BiCWP}, we present the 
\textit{non-OP points} in the numerical solutions of SLP with a 
uniform mesh size of $N=3200$ computed by the WENO-M and 
MOP-WENO-ACM$k$ schemes, and the \textit{non-OP points} in the 
numerical solutions of BiCWP with a uniform mesh size of $N=6400$ 
computed by the MIP-WENO-ACM$k$ and MOP-WENO-ACM$k$ schemes. We 
can see that there are a great many \textit{non-OP points} in the 
solutions of the SLP computed by the WENO-M scheme and in the 
solutions of the BiCWP computed by the MIP-WENO-ACM$k$ scheme, while 
the numbers of the \textit{non-OP points} in the solutions of these 
two cases computed by the MOP-WENO-ACM$k$ scheme are zero. Actually, 
there are many \textit{non-OP points} for all considered mapped 
schemes whose mapping functions are \textit{non-OP}, like the 
WENO-M, WENO-PM6, WENO-IM(2, 0.1)) and MIP-WENO-ACM$k$ schemes. As 
expected, there are no \textit{non-OP points} for the MOP-WENO-ACM$k$
and WENO-JS schemes for all computing cases here. We do not show the 
results of the \textit{non-OP points} for all computing cases here 
just for the simplicity of illustration. It should be noted that the 
WENO-JS scheme could be treated as a mapped WENO scheme whose 
mapping function is an identical mapping, that is, 
$\big( g^{\mathrm{JS}} \big)_{s}(\omega) = \omega, s = 0,1,2$, and 
it is trivial to verify that the set of mapping functions $\big( 
g^{\mathrm{JS}} \big)_{s}(\omega), s=0,1,2$ is \textit{OP} while its 
optimal weight intervals are zero.

Thus, to sum up, we can conclude that a set of mapping functions which is \textit{OP} can help to improve the resolution of the 
corresponding mapped WENO scheme and prevent it from generating 
spurious oscillations in the simulation of problems with 
discontinuities, especially for long output times. And in the 
upcoming follow-up study of this article, we will provide more 
examples and evidence to further verify this conclusion.

\subsection{Two-dimensional Euler system}
\label{subsec:examples_2D_Euler}
In this subsection, we solve the two-dimensional Euler system of gas 
dynamics. We consider the numerical solutions of the 2D Riemann 
problem \cite{Riemann-2D-01,Riemann2D-02,Riemann2D-03} and the 
shock-vortex interaction problem \cite{Shock-vortex_interaction-1,
Shock-vortex_interaction-2,Shock-vortex_interaction-3}. The 
two-dimensional Euler system is given by the following strong 
conservation form
\begin{equation}
\mathbf{U}_t + \mathbf{F}\big( \mathbf{U} \big)_x +
\mathbf{G}\big( \mathbf{U} \big)_y= 0,
\label{2DEulerEquations}
\end{equation}
where $\mathbf{U} = \big(\rho, \rho u, \rho v, E\big)^{\mathrm{T}}, 
\mathbf{F}\big( \mathbf{U} \big) = \big( \rho u, \rho u^{2} + p, \rho
uv, u(E + p) \big)^{\mathrm{T}}, \mathbf{G}\big( \mathbf{U} \big) = 
\big( \rho v, \rho vu, \rho v^{2} + p, v(E + p) \big)^{\mathrm{T}}$, 
and $\rho, u, v, p$ and $E$ are the density, component of velocity 
in the $x$ and $y$ coordinate directions, pressure and total energy, 
respectively. The relation of the pressure and the total energy, the 
component of velocity in the $x$ and $y$ coordinate directions is 
defined by the equation of state for an ideal polytropic gas, taking 
the form
\begin{equation*}
p = (\gamma - 1)\Big( E - \dfrac{1}{2}\rho (u^{2} + v^{2}) \Big), 
\end{equation*}
where $\gamma$ is the ratio of specific heats and we choose 
$\gamma = 1.4$ here. In all numerical examples of this subsection, 
the CFL number is set to be $0.5$.

\begin{example}
\bf{(2D Riemann problem)} 
\rm{The series of 2D Riemann problems proposed in 
\cite{Riemann-2D-01,Riemann2D-02} has become favorable cases to test 
the resolution of numerical methods \cite{Riemann2D-03,WENO-PPM5,
Riemann2D-04}. It is calculated over a unit square domain $[0,1] 
\times[0,1]$, initially involves the constant states of flow 
variables over each quadrant which is got by dividing the 
computational domain using lines $x = x_{0}$ and $y = y_{0}$. 
Configuration 4 in \cite{Riemann2D-03} is taken here for the test. 
The initial condition of this configuration is given by}
\label{ex:Riemann2D}
\end{example}
\begin{equation*}
\big( \rho, u, v, p \big)(x, y, 0) = \left\{
\begin{aligned}
\begin{array}{ll}
(1.1, 0.0, 0.0, 1.1),         & 0.5 \leq x \leq 1.0, 
                                0.5 \leq y \leq 1.0, \\
(0.5065, 0.8939, 0.0, 0.35),  & 0.0 \leq x \leq 0.5, 
                                0.5 \leq y \leq 1.0, \\
(1.1, 0.8939, 0.8939, 1.1),   & 0.0 \leq x \leq 0.5, 
                                0.0 \leq y \leq 0.5, \\
(0.5065, 0.0, 0.8939, 0.35),  & 0.5 \leq x \leq 1.0, 
                                0.0 \leq y \leq 0.5. \\
\end{array}
\end{aligned}
\right.
\label{eq:initial_Euer2D:Riemann2D}
\end{equation*}
The transmission boundary condition is used on all boundaries. The 
numerical solutions are calculated using considered WENO schemes on 
$800 \times 800$ cells, and the computations proceed to $t = 0.25$.

\begin{figure}[ht]
\centering
  \includegraphics[height=0.29\textwidth]
  {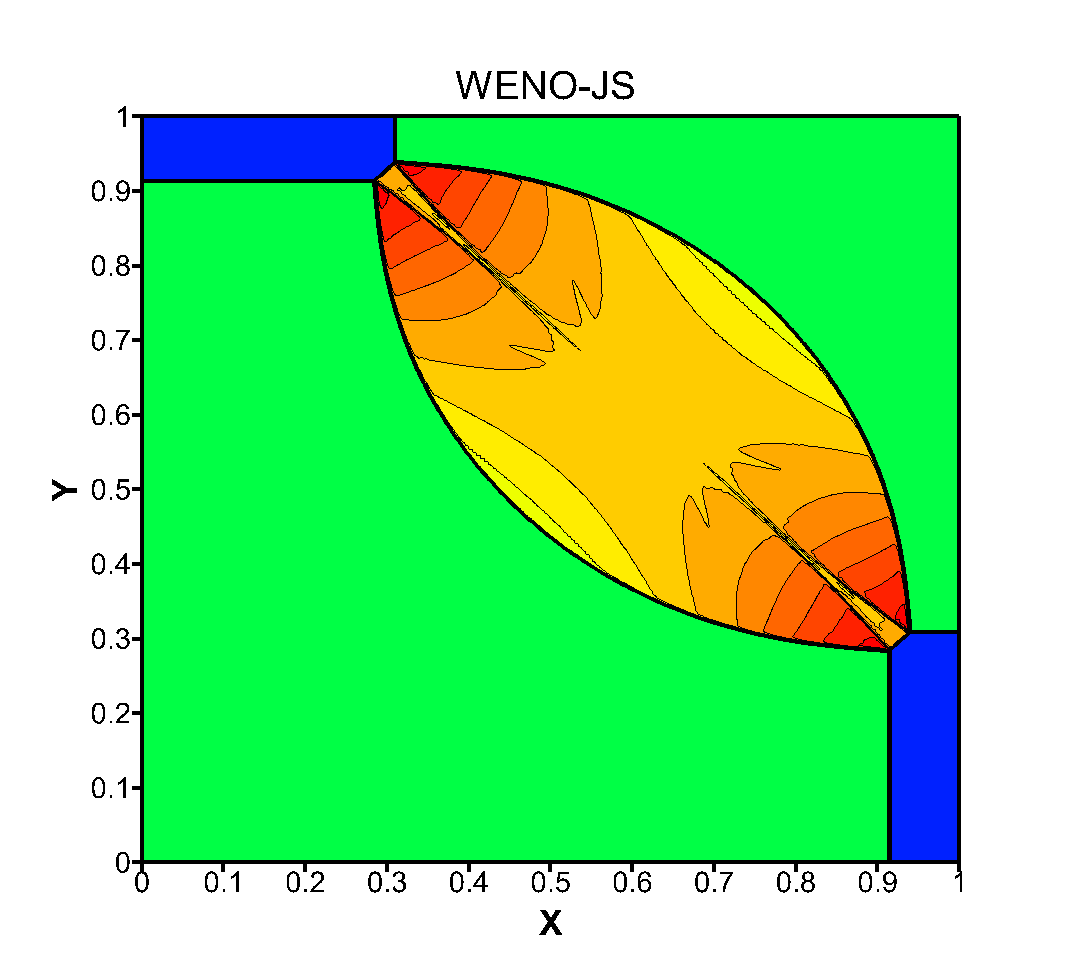}
  \includegraphics[height=0.29\textwidth]
  {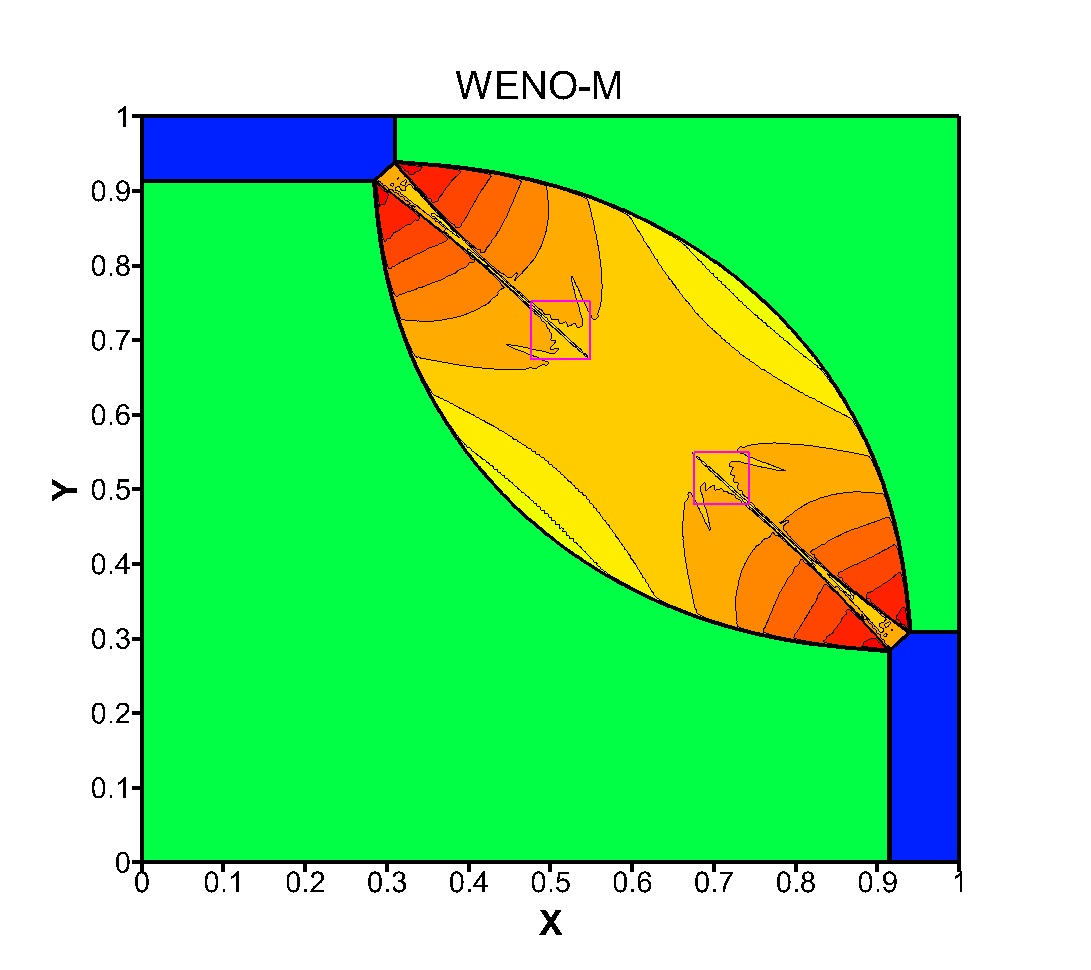}
  \includegraphics[height=0.29\textwidth]
  {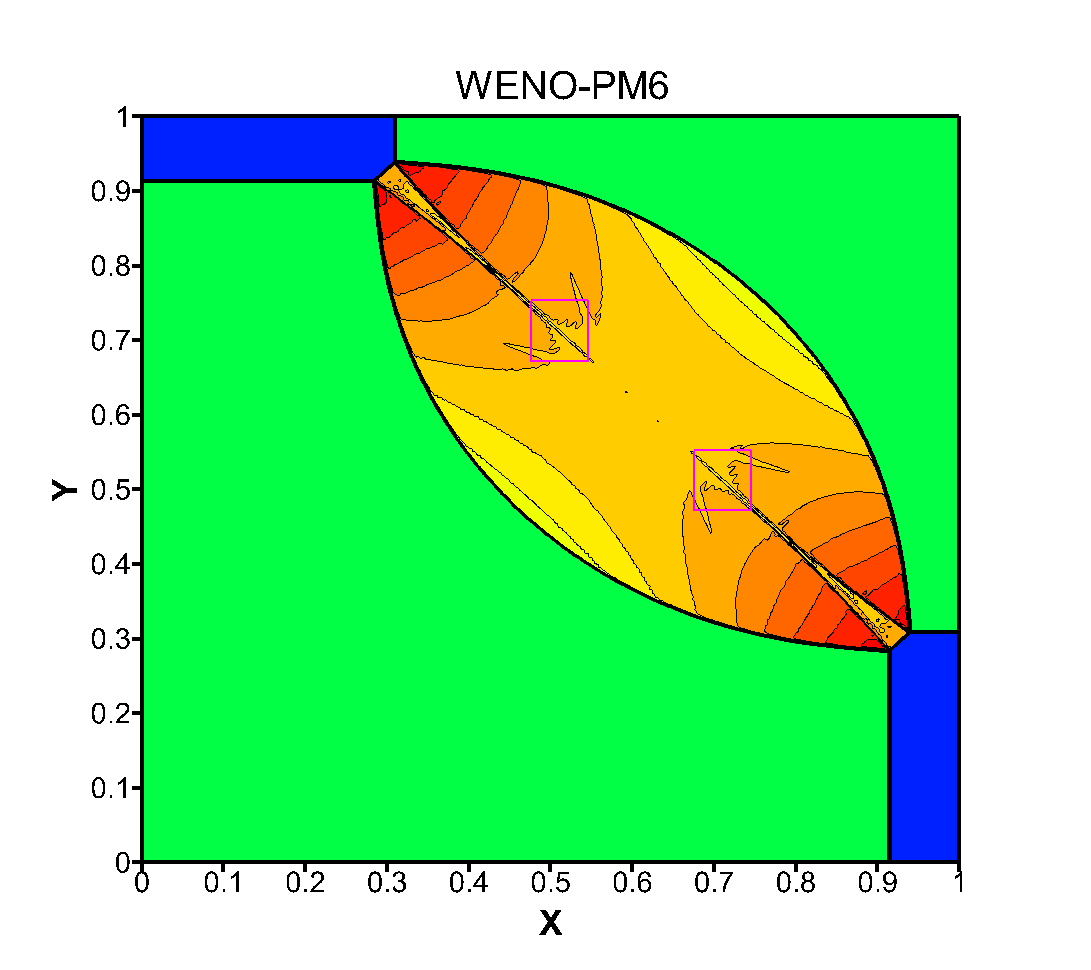}\\
  \includegraphics[height=0.29\textwidth]
  {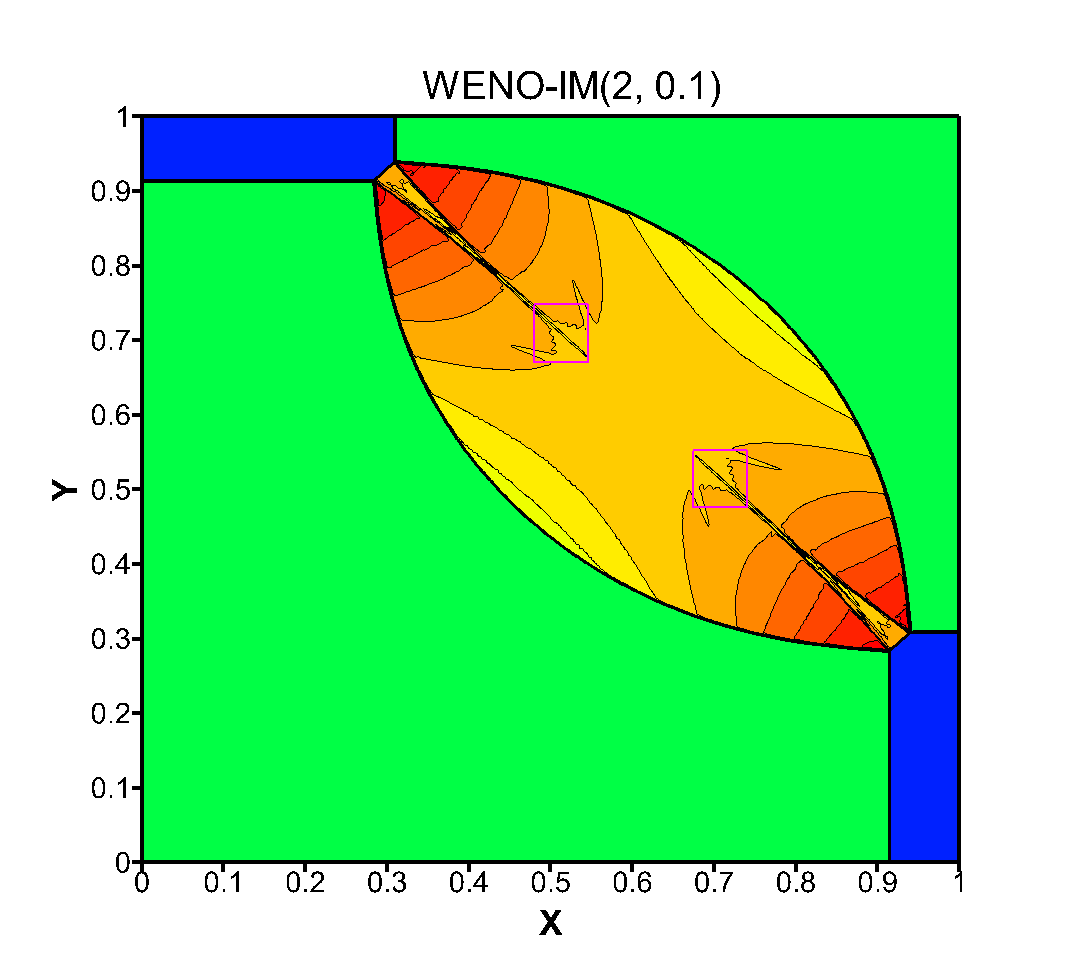}
  \includegraphics[height=0.29\textwidth]
  {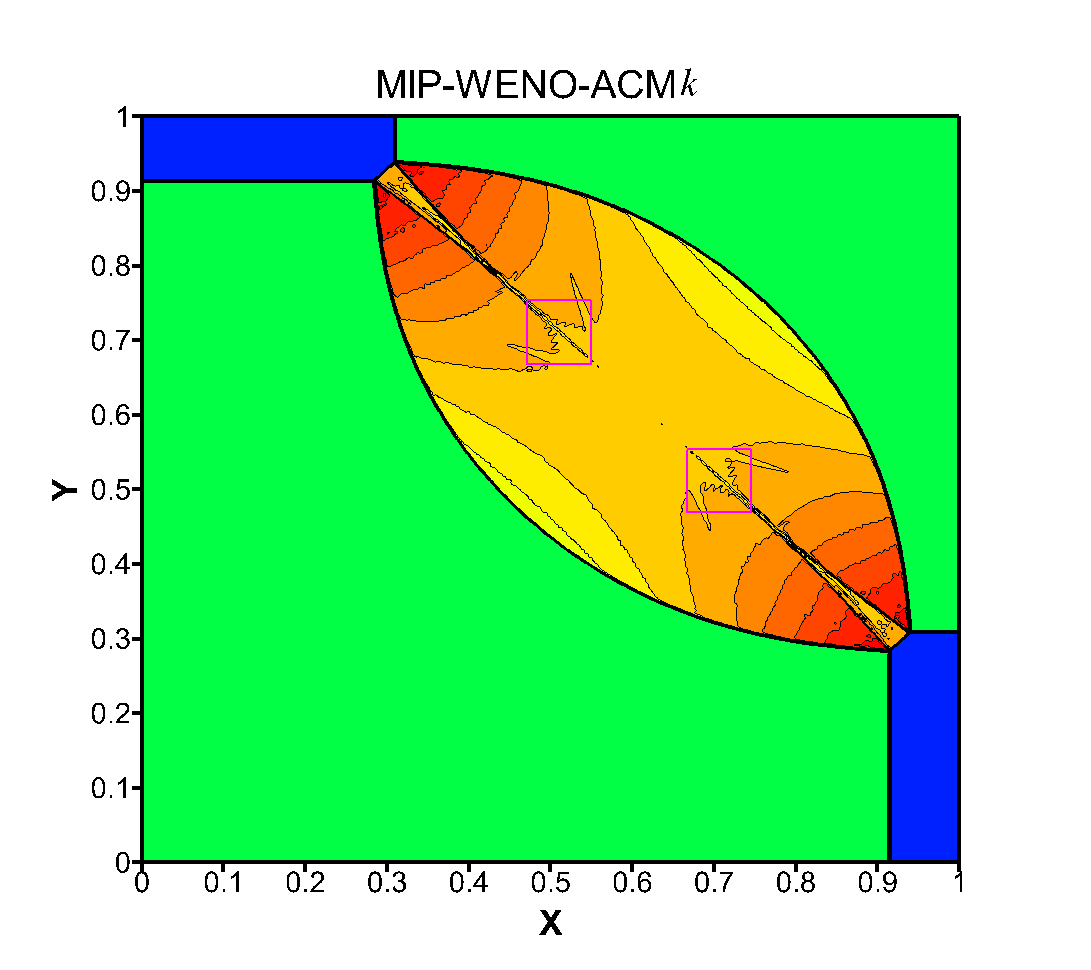}
  \includegraphics[height=0.29\textwidth]
  {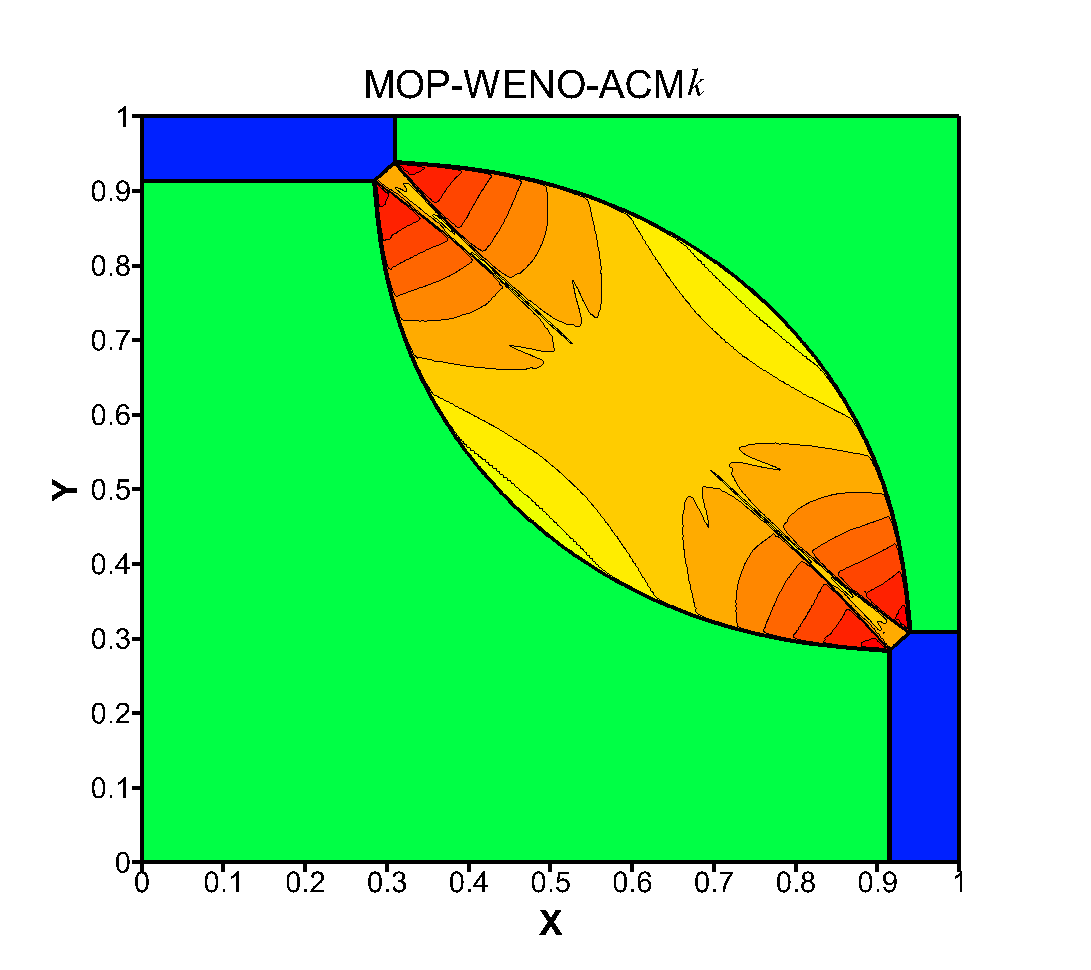}    
\caption{Density plots for the 2D Riemann problem using $30$ contour
lines with range from $0.5$ to $1.9$, computed using the WENO-JS, 
WENO-M, WENO-PM6, WENO-IM($2, .01$), MIP-WENO-ACM$k$ and 
MOP-WENO-ACM$k$ schemes.}
\label{fig:ex:Riemann2D}
\end{figure}

\begin{figure}[ht]
\centering
  \includegraphics[height=0.42\textwidth]
  {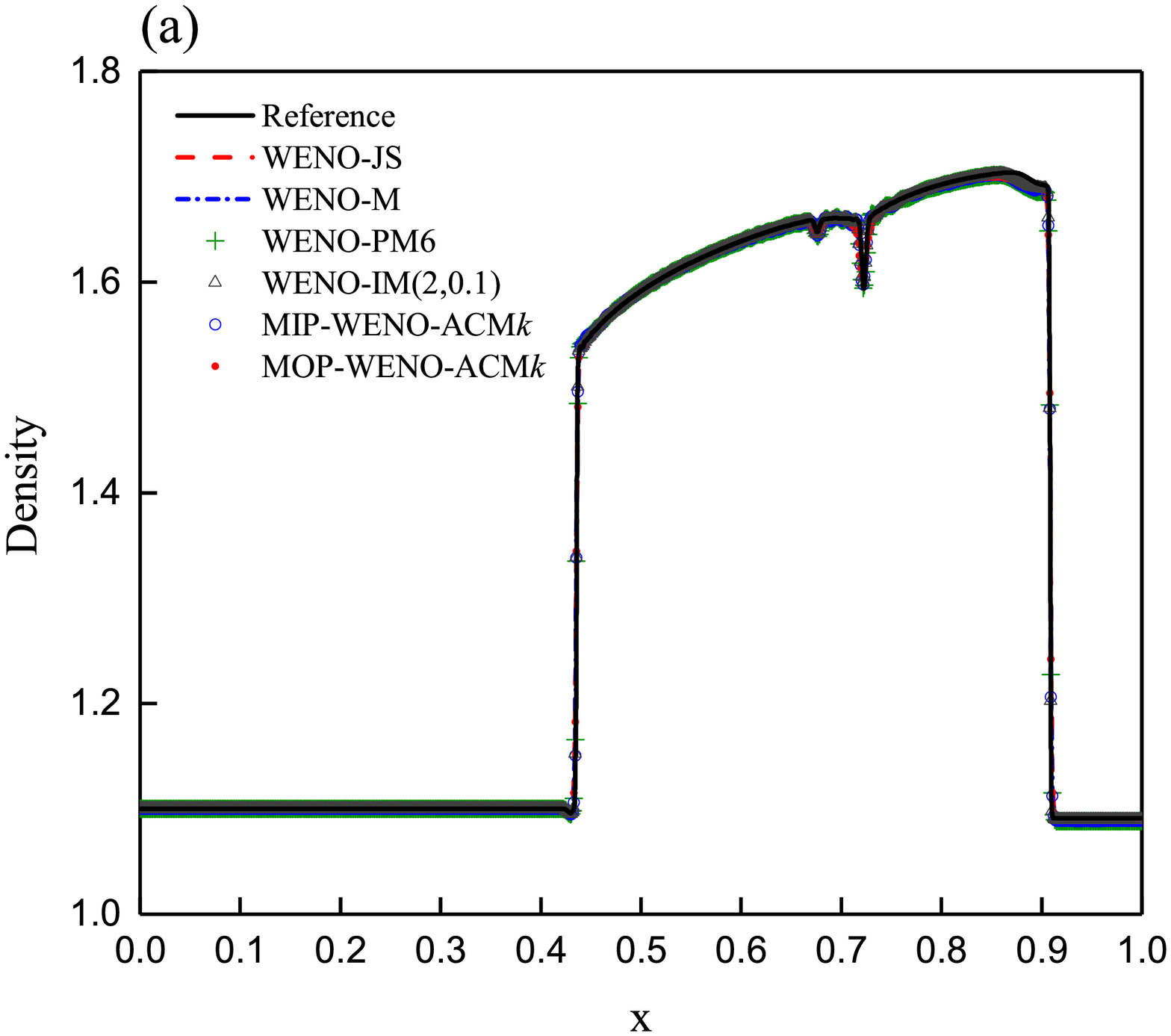}
  \includegraphics[height=0.42\textwidth]
  {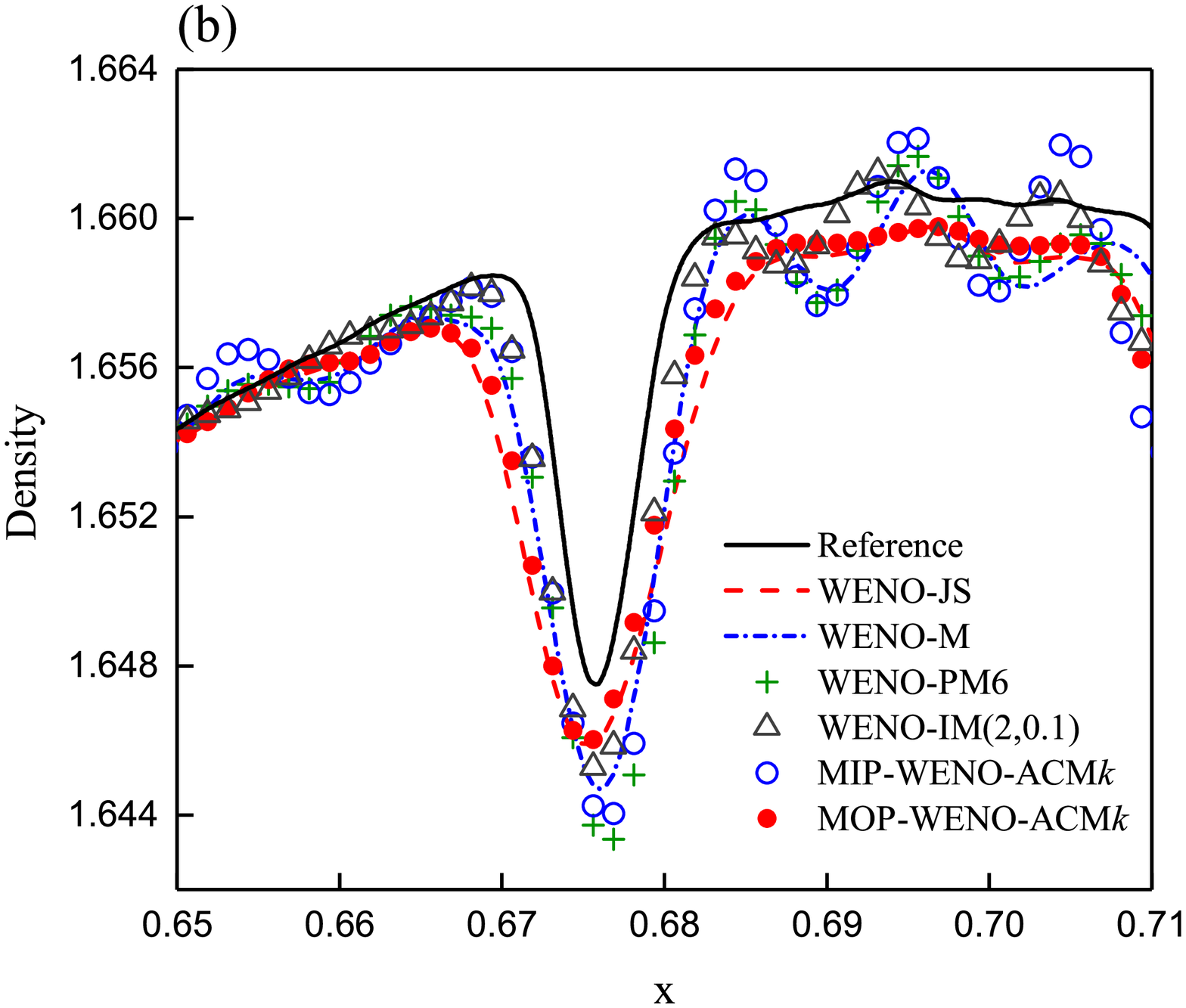}
\caption{The cross-sectional slices of density plot along the plane 
$y = 0.5$, computed using the WENO-JS, WENO-M, WENO-PM6, WENO-IM
($2, 0.1$), MIP-WENO-ACM$k$ and MOP-WENO-ACM$k$ schemes.}
\label{fig:ex:Riemann2D-ZoomedIn}
\end{figure}

In Fig. \ref{fig:ex:Riemann2D}, we have shown the numerical results 
of density obtained by using the WENO-JS, WENO-M, WENO-PM6, 
WENO-IM($2,0.1$), MIP-WENO-ACM$k$ and MOP-WENO-ACM$k$ schemes. We 
can see that all considered schemes can capture the main structure 
of the solution. However, we can also observe that there are obvious 
numerical oscillations (as marked by the pink boxes), which are 
unfavorable for the fidelity of the results, in the solutions of the 
WENO-M, WENO-PM6, WENO-IM(2, 0.1) and MIP-WENO-ACM$k$ schemes. These 
numerical oscillations can be seen more clearly from the 
cross-sectional slices of density profile along the plane $y = 0.5$ 
as presented in Fig. \ref{fig:ex:Riemann2D-ZoomedIn}, where the 
reference solution is obtained by using the WENO-JS scheme with a 
uniform mesh size of $3000\times3000$. Noticeably, there are almost 
no numerical oscillations in the solutions of the MOP-WENO-ACM$k$ 
and WENO-JS schemes, and this should be an advantage of the mapped 
WENO schemes whose mapping functions are \textit{OP}.

\begin{example}
\bf{(Shock-vortex interaction)}
\rm{We solve the shock-vortex interaction problem
\cite{Shock-vortex_interaction-1,Shock-vortex_interaction-2,
Shock-vortex_interaction-3} that consists of the interaction of a 
left moving shock wave with a right moving vortex. The initial 
condition is given by}
\label{ex:shock-vortex}
\end{example}
\begin{equation*}
\big( \rho, u, v, p \big)(x, y, 0) = \left\{
\begin{aligned}
\begin{array}{ll}
\textbf{U}_{\mathrm{L}}, & x < 0.5, \\
\textbf{U}_{\mathrm{R}}, & x \geq 0.5, \\
\end{array}
\end{aligned}
\right.
\label{eq:initial_Euer2D:shock-vortex-interaction}
\end{equation*}
where the left state is taken as $\textbf{U}_{\mathrm{L}} = 
(\rho_{\mathrm{L}}, u_{\mathrm{L}}, v_{\mathrm{L}}, p_{\mathrm{L}}) =
(1, \sqrt{\gamma}, 0, 1)$, and the right state 
$\textbf{U}_{\mathrm{R}} = (\rho_{\mathrm{R}}, u_{\mathrm{R}}, 
v_{\mathrm{R}}, p_{\mathrm{R}})$ is given by
\begin{equation*}
\begin{array}{l}
p_{\mathrm{R}} = 1.3, \rho_{\mathrm{R}} = \rho_{\mathrm{L}}\bigg( 
\dfrac{\gamma - 1 + (\gamma + 1)p_{\mathrm{R}}}{\gamma + 1 + (\gamma
- 1)p_{\mathrm{R}}} \bigg)\\
u_{\mathrm{R}} = u_{\mathrm{L}}\bigg( \dfrac{1 - p_{\mathrm{R}}}{
\sqrt{\gamma - 1 + p_{\mathrm{R}}(\gamma + 1)}}\bigg), v_{\mathrm{R}}
= 0.
\end{array}
\label{eq:Euler2D:shock-vortex-interaction:rightState}
\end{equation*}
A vortex given by the following perturbations is superimposed onto 
the left state $\textbf{U}_{\mathrm{L}}$,
\begin{equation*}
\delta \rho = \dfrac{\rho_{\mathrm{L}}^{2}}{(\gamma - 1)
p_{\mathrm{L}}}\delta T, 
\delta u = \epsilon \dfrac{y - y_{\mathrm{c}}}{r_\mathrm{c}}
\mathrm{e}^{\alpha(1-r^{2})}, 
\delta v = - \epsilon \dfrac{x - x_{\mathrm{c}}}{r_\mathrm{c}}
\mathrm{e}^{\alpha(1-r^{2})}, 
\delta p = \dfrac{\gamma \rho_{\mathrm{L}}^{2}}{(\gamma - 1)
\rho_{\mathrm{L}}}\delta T,
\label{eq:Euler2D:shock-vortex-interactions:Perturbations}
\end{equation*}
where $\epsilon = 0.3, r_{\mathrm{c}} = 0.05, \alpha = 0.204, 
x_{\mathrm{c}} = 0.25, y_{\mathrm{c}} = 0.5,
r = \sqrt{((x - x_{\mathrm{c}})^{2} + (y - y_{\mathrm{c}})^{2})/r_{
\mathrm{c}}^{2}}, \delta T = - (\gamma - 1)\epsilon^{2}\mathrm{e}^{2
\alpha (1 - r^{2})}/(4\alpha \gamma)$.
The transmissive boundary condition is used on all boundaries.

\begin{figure}[ht]
\centering
  \includegraphics[height=0.33\textwidth]
  {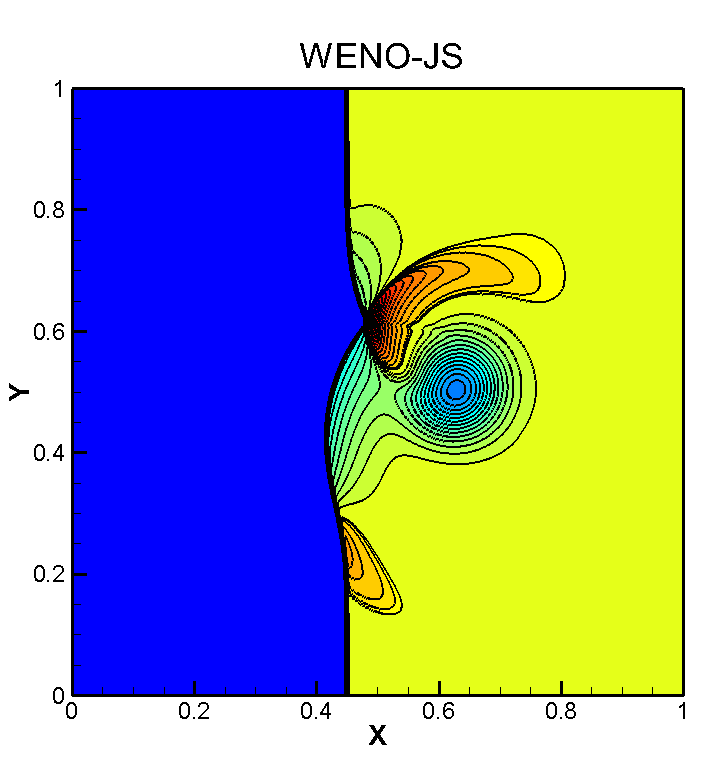}
  \includegraphics[height=0.33\textwidth]
  {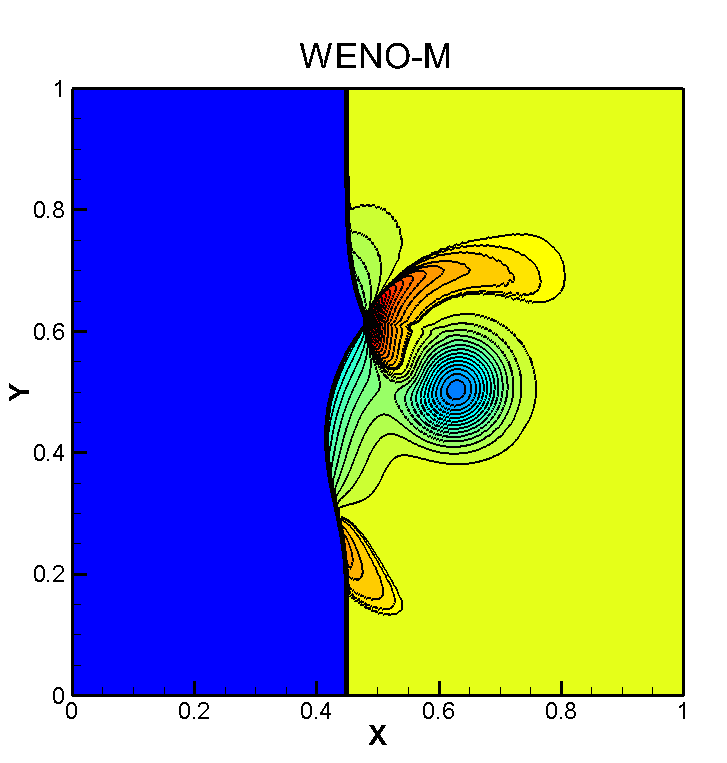}
  \includegraphics[height=0.33\textwidth]
  {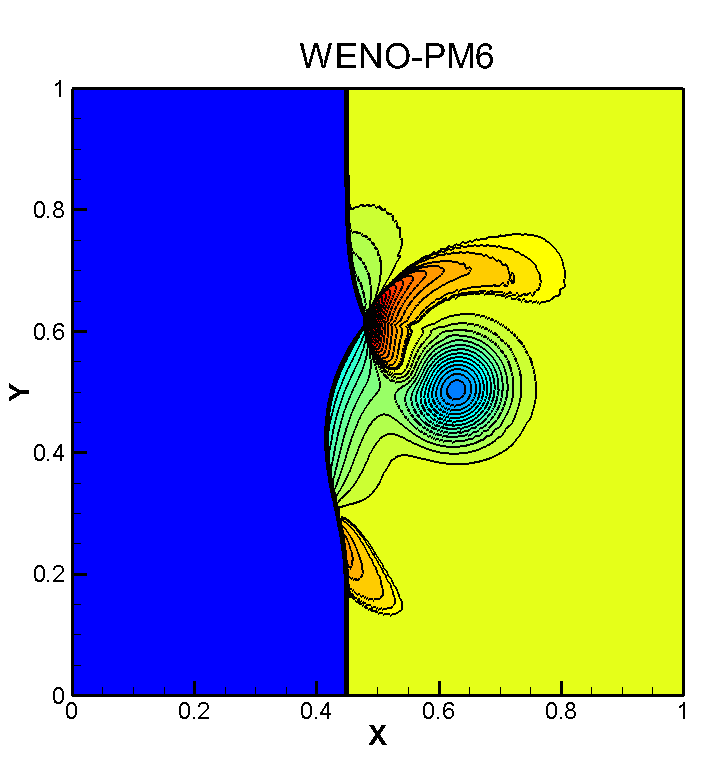}\\
  \includegraphics[height=0.33\textwidth]
  {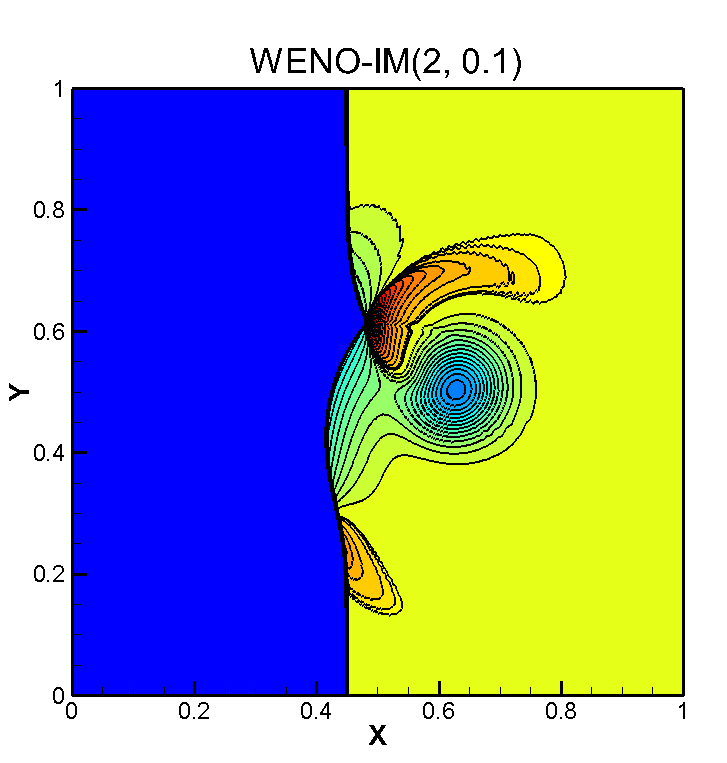}
  \includegraphics[height=0.33\textwidth]
  {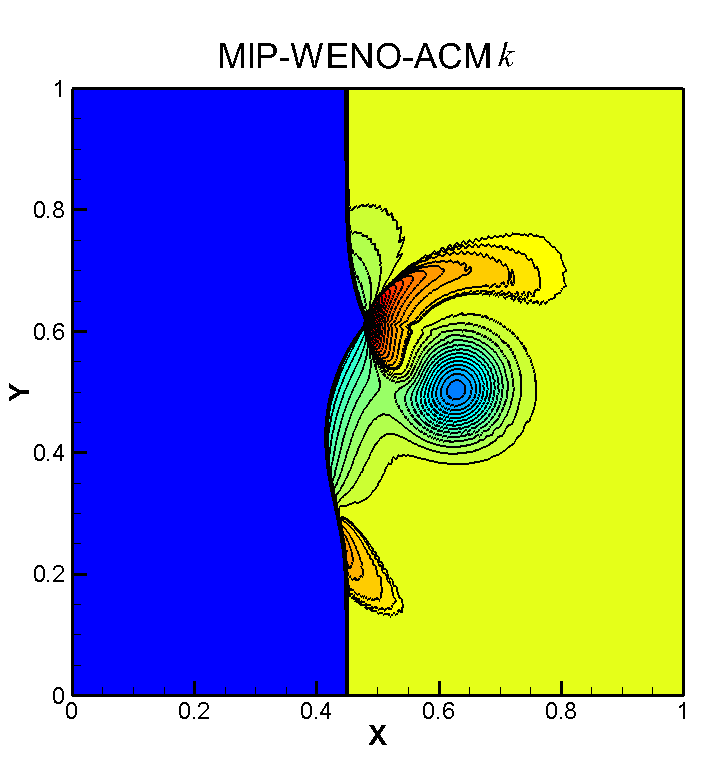}
  \includegraphics[height=0.33\textwidth]
  {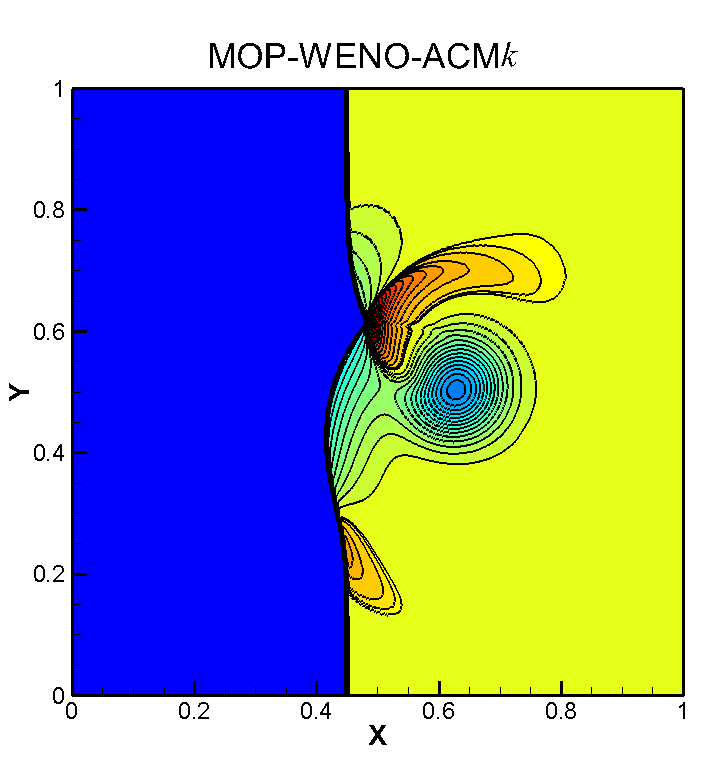}    
\caption{Density plots for the Shock-vortex interaction using $30$ 
contour lines with range from $0.9$ to $1.4$, computed using the 
WENO-JS, WENO-M, WENO-PM6, WENO-IM($2, .01$), MIP-WENO-ACM$k$ and 
MOP-WENO-ACM$k$ schemes.}
\label{fig:ex:SVI}
\end{figure}

\begin{figure}[ht]
\centering
  \includegraphics[height=0.42\textwidth]
  {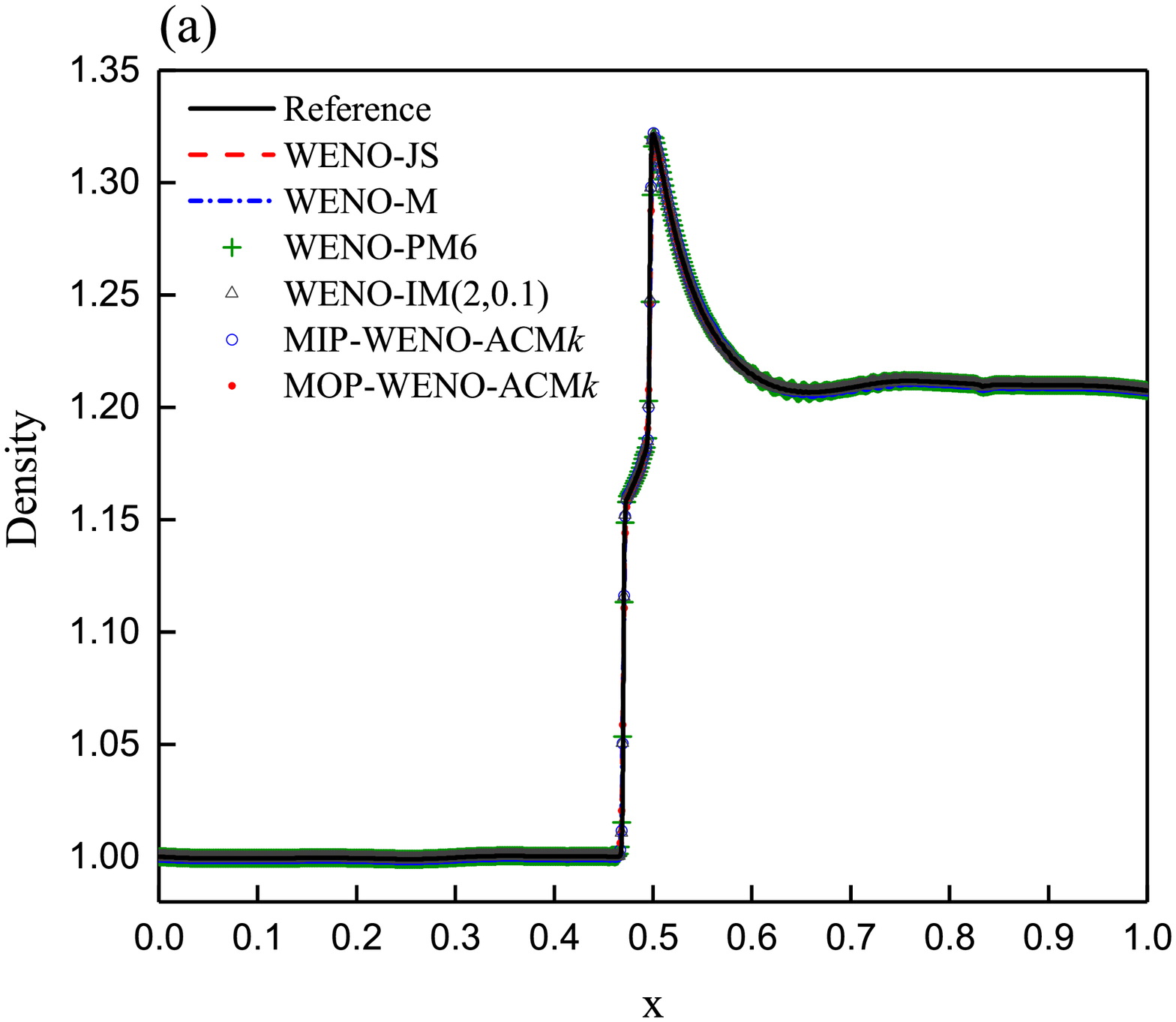}
  \includegraphics[height=0.42\textwidth]
  {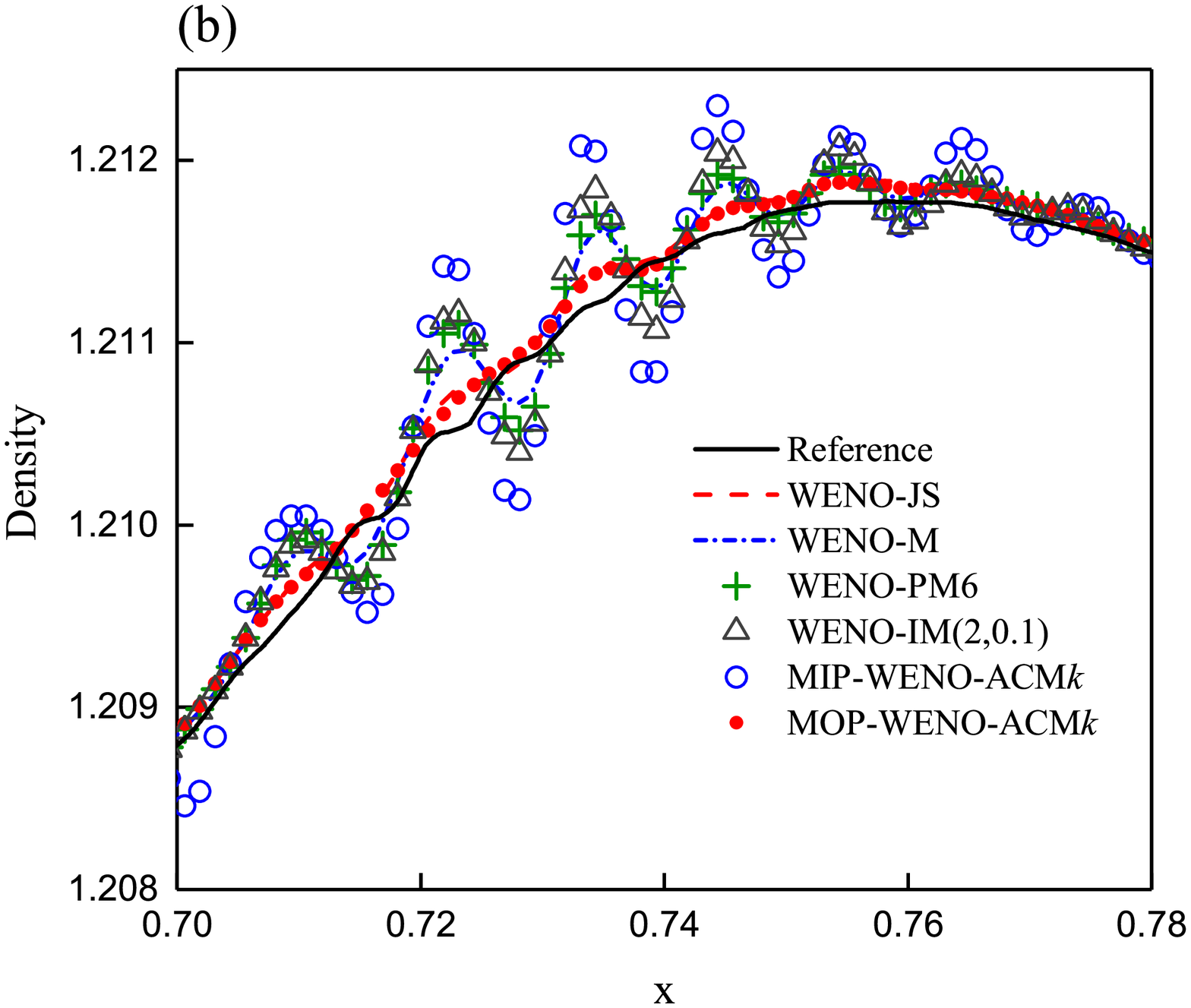}
\caption{The cross-sectional slices of density plot along the plane 
$y = 0.65$, computed using the WENO-JS, WENO-M, WENO-PM6, WENO-IM
($2, 0.1$), MIP-WENO-ACM$k$ and MOP-WENO-ACM$k$ schemes.}
\label{fig:ex:SVI-ZoomedIn}
\end{figure}

The problem has been calculated by the considered WENO schemes with 
a uniform mesh size of $800 \times 800$ and the output time is taken 
as $t = 0.35$. The final structures of the shock and vortex in the 
density profile have been shown in Fig. \ref{fig:ex:SVI}. It is 
observed that all the considered WENO schemes perform well in 
capturing the main structure of the shock and vortex after the 
interaction. We can see that there are clear numerical oscillations 
in the solutions of the WENO-IM(2, 0.1) and MIP-WENO-ACM$k$ schemes, 
and the numerical oscillations can also be observed in the solutions 
of the WENO-M and WENO-PM6 schemes although they are not so severe 
as those of the WENO-IM(2, 0.1) and MIP-WENO-ACM$k$ schemes. However,
in the solutions of the MOP-WENO-ACM$k$ and WENO-JS schemes, we 
almost did not find the numerical oscillations. To further 
demonstrate this, we have plotted the cross sectional slices of 
density profile along the plane $y=0.65$ in Fig.
\ref{fig:ex:SVI-ZoomedIn}. The reference solution is obtained 
using the WENO-JS scheme with a uniform mesh size of 
$1600 \times 1600$. It is evident that the MIP-WENO-ACM$k$ scheme 
produces the numerical oscillations with the biggest amplitudes 
followed by those of the WENO-IM(2, 0.1) scheme. The WENO-PM6 
and WENO-M schemes also generate clear numerical oscillations with 
the amplitudes slightly smaller than that of the WENO-IM(2, 0.1) 
scheme. Obviously, the solutions of the MOP-WENO-ACM$k$ and WENO-JS 
schemes almost generate no numerical oscillations or only generate 
some imperceptible numerical oscillations, and their solutions are 
most close to the reference solution. Again, we argue that this 
should be an advantage of the mapped WENO schemes whose mapping 
functions are \textit{OP}.


\section{Conclusions}
\label{secConclusions} 
This paper has proposed a new mapped weighted essentially 
non-oscillatory scheme named as MOP-WENO-ACM$k$. The motivation to 
design this new scheme is that: (1) the WENO-JS and WENO-M schemes 
generate numerical solutions with very low resolutions when solving 
hyperbolic problems with discontinuities for long output 
times; (2) although various existing improved mapped WENO schemes 
can successfully address the drawback above, as far as we know, 
almost all of them introduce unfavorable spurious oscillations as 
the \textit{non-OP mapping process} occurs in their mappings. By 
introducing a set of mapping functions that is 
\textit{order-preserving (OP)}, the MOP-WENO-ACM$k$ scheme can 
prevent the \textit{non-OP mapping process}, which should be the 
essential cause of the spurious oscillation generation and potential 
loss of accuracy. Therefore, the MOP-WENO-ACM$k$ scheme has a 
significant advantage in that it not only can obtain comparable high 
resolutions but also can prevent generating spurious oscillations, 
when solving problems with discontinuities, especially for long 
output times. Numerical experiments have shown that the proposed 
scheme yields lower dissipation and higher resolution near 
discontinuities than the WENO-JS and WENO-M schemes, especially for 
long output times, and it enjoys better robustness than the 
WENO-PM6, WENO-IM(2, 0.1) and MIP-WENO-ACM$k$ schemes.

\bibliographystyle{model1b-shortjournal-num-names}
\bibliography{refs}

\end{document}